\documentclass[11pt]{amsart}
\usepackage[utf8]{inputenc}
\pdfoutput=1

\usepackage{fullpage}
\usepackage{amsmath,amssymb,amsthm}

\usepackage{tikz}
\usetikzlibrary{calc}

\usepackage[backend=bibtex,url=false]{biblatex}
\usepackage{hyperref}
\usepackage{xurl}
\hypersetup{linkcolor = bluegreen, citecolor = bluegreen, colorlinks = True, breaklinks=true}
\usepackage{cleveref}

\usepackage{color}
\usepackage{enumitem}
\setlist[enumerate,1]{label={(\roman*)}}

\newtheorem{theorem}{Theorem}[section]

\newtheorem{lemma}[theorem]{Lemma}
\newtheorem{corollary}[theorem]{Corollary}
\newtheorem{proposition}[theorem]{Proposition}

\newtheorem{question}[theorem]{Question}

\newtheorem{remark}[theorem]{Remark}
\newtheorem{definition}[theorem]{Definition}

\counterwithin{figure}{section}

\definecolor{green}{RGB}{0,144,0}
\definecolor{bluegreen}{RGB}{17,100,180}

\newcommand*{\abs}[1]{\left\lvert#1\right\rvert}

\def\N{\mathbb{N}}
\def\Z{\mathbb{Z}}
\def\Q{\mathbb{Q}}
\def\R{\mathbb{R}}

\newcommand{\CF}{\mathcal{F}}
\newcommand{\CA}{\mathcal{A}}

\DeclareMathOperator{\dimH}{\mathrm{dim}_{\mathrm{H}}}
\DeclareMathOperator{\dimB}{\mathrm{dim}_{\mathrm{B}}}
\newcommand{\dloc}{d_{\mathrm{loc}}}
\DeclareMathOperator{\diam}{diam}
\DeclareMathOperator{\sizes}{\mathsf{s}}

\title{On the classical Lagrange and Markov spectra: new results on the local dimension and the geometry of the difference set}

\author{Harold Erazo}

\address[Harold Erazo]{IMPA, Estrada Dona Castorina 110, 22460-320, Rio de Janeiro, Brazil}
	\email{harolderaz@gmail.com}

\author{Luke Jeffreys}
\address[Luke Jeffreys]{School of Mathematics , University of Bristol, Fry Building, Woodland Road, Bristol BS8 1UG, UK}
\curraddr{}
\email{luke.jeffreys@bristol.ac.uk}

\author{Carlos Gustavo Moreira}
	\address[Carlos Gustavo Moreira]{SUSTech International Center for Mathematics, Shenzhen, Guangdong, People’s Republic of China; \hfill\break
IMPA, Estrada Dona Castorina 110, 22460-320, Rio de Janeiro, Brazil}
\email{gugu@impa.br}

\thanks{The first author is partially supported by CAPES and FAPERJ. The second author is a Leverhulme Early Career Fellow and thanks the Leverhulme Trust for their support as well as the University of Wisconsin-Madison where he was a Van Vleck Visiting Assistant Professor when this work began. The third author is partially supported by CNPq and FAPERJ}

\date{}


\begin{document}

\begin{abstract}
Let $L$ and $M$ denote the classical Lagrange and Markov spectra, respectively. It is known that $L\subset M$ and that $M\setminus L\neq\varnothing$. Inspired by three questions asked by the third author in previous work investigating the fractal geometric properties of the Lagrange and Markov spectra, we investigate the function $d_{loc}(t)$ that gives the local Hausdorff dimension at a point $t$ of $L'$. Specifically, we construct several intervals (having non-trivial intersection with $L'$) on which $d_{loc}$ is non-decreasing. We also prove that the respective intersections of $M'$ and $M''$ with these intervals coincide. Furthermore, we completely characterize the local dimension of both spectra when restricted to those intervals. Finally, we demonstrate the largest known elements of the difference set $M\setminus L$ and describe two new maximal gaps of $M$ nearby.
\end{abstract}


\maketitle

\section{Introduction}

The classical Lagrange and Markov spectra are two subsets of the real line related to the study of Diophantine approximation. Given a positive real number $\alpha$ we define its \emph{best constant of Diophantine approximation} to be
\[k(\alpha) := \limsup_{p,q\to\infty} \frac{1}{|q(q\alpha - p)|} = \sup\left\{k>0\,:\,\left|\alpha - \frac{p}{q}\right| < \frac{1}{kq^{2}} \text{has infinitely many solutions }\frac{p}{q}\in\Q\right\}.\]
The \emph{Lagrange spectrum} is defined to be the set
\[L:=\{k(\alpha)\,:\,\alpha\in\R\setminus\Q\}.\]
The \emph{Markov spectrum} is related to the approximation of binary quadratic forms and is defined to be
\[M:=\left\{\sup_{(p,q)\in\Z^{2}\setminus\{(0,0)\}}\frac{1}{|ap^{2}+bpq+cq^{2}|}\,:\,ax^{2}+bxy+cy^{2}\text{ real indefinite, }b^{2}-4ac = 1\right\}.\]
It is known that $L\subset M\subset\R^{+}$.

In the late 1800s, Markov~\cite{Markov79,Markov80} determined that
\[L\cap(0,3) = M\cap(0,3) = \left\{\sqrt{5}<2\sqrt{2}<\frac{\sqrt{221}}{5}<\ldots\right\} = \left\{\sqrt{9-\frac{4}{m^{2}}}\,:\,m\,\text{is a Markov number}\right\},\]
where a Markov number is the largest number in a triple of positive integers $(x,y,z)$ satisfying the so-called Markov equation $x^{2}+y^{2}+z^{2} = 3xyz$. In 1975, Freiman~\cite{Fr75} showed that $[c_{F},\infty)\subset L\subset M$ and $(\nu_{F},c_F)\cap M=\varnothing$ where $c_{F} = 4.527829566\ldots$ and $\nu_{F} = 4.527829538\ldots\in M$. Here, $c_{F}$ and $\nu_{F}$ are algebraic numbers, whose explicit expressions as sums of eventually periodic continued fractions are given in Section 1.1. This ray $[c_{F},\infty)$ is known as Hall's ray after earlier work of Hall~\cite{Ha47} (see also the intermediate results of Freiman-Judin~\cite{FJ66}, Hall~\cite{Ha71}, Freiman~\cite{Fr73} and Schecker~\cite{Sch77}) and $c_{F}$ is now known as Freiman's constant. Hall's ray is currently the only known continuous part of the spectra. It is a long-standing and wide-open conjecture of Berstein~\cite{BerConj} that $[4.1,4.52]\subset L\subset M$.

Between 3 and $c_{F}$, both $L$ and $M$ have interesting fractal structure. Indeed, the third author proved~\cite{geometricproperties} that
\[d(t):=\dimH(L\cap(-\infty,t)) = \dimH(M\cap(-\infty,t)),\]
is a continuous function with $d(3+\epsilon)>0$ for any $\epsilon>0$ and $d(\sqrt{12}) = 1$. In the same paper, it was also proved that
\[D(t) := \dimH(k^{-1}(-\infty,t]) = \dimH(k^{-1}(-\infty,t)),\]
is also a continuous function and that $d(t)=\min\{1,2\cdot D(t)\}$.

More recently, the third author in joint work with Matheus, Pollicott and Vytnova~\cite{MMPV} determined that
\[t_1 := \min\{t\in\R : d(t) = 1\} = 3.334384....\]
One can also define the local dimension function $d_{loc}:L'\to [0,1]$ by
\[d_{loc}(t) := \lim_{\epsilon\to 0}\dimH(L\cap(t-\epsilon,t+\epsilon)).\]

As mentioned above, it is known that $L\subset M$. Freiman~\cite{Fr68} also showed that $M\setminus L \neq \varnothing$. Recently, in the same work of the third author, Matheus, Pollicott and Vytnova discussed above, it was shown that the Hausdorff dimension $\dimH(M\setminus L)$ of $M\setminus L$ satisfies
\[0.537152 < \dimH(M\setminus L) <  0.796445.\]
This lower bound was heuristically improved to 0.593 in recent work of the second and third authors with Matheus~\cite{JMM}. 

In~\cite{geometricproperties}, the third author raised the following questions for further investigation:
\begin{itemize}
\item[\textbf{1)}] Is the function $d_{loc}$ non-decreasing?
\item[\textbf{2)}] What is the geometric structure of the difference set $M\setminus L$?
\item[\textbf{3)}] Is $M'' = M'$?
\end{itemize}

In this work, we will investigate these questions and will answer them for certain subsets of the spectra lying within intervals of the real line that we call ``good intervals" (see Definition~\ref{def:good_interval}). We will show that we have a complete understanding of the local dimension of both spectra restricted to these intervals. As a by-product we also improve the upper bound on $\dimH(M\setminus L)$ in these regions.

The above questions remain open in general.

\subsection{Continued fractions and shift space dynamics}

Before discussing the main results, we must introduce the modern notation used to describe the Lagrange and Markov spectra.

The values in the Lagrange and Markov spectra can be calculated using the theory of continued fractions thanks to the work of Perron~\cite{Pe21} who proved that if we have
\[\alpha = [a_{0};a_{1},a_{2},\ldots] := a_{0} + \frac{1}{a_{1} + \frac{1}{a_{2} + \frac{1}{\dots}}},\]
then
\[k(\alpha) = \limsup_{n\to\infty}\,( [a_{n};a_{n-1},\ldots,a_{1}] + [0;a_{n+1},a_{n+2},\ldots]).\]
This allows us to study the structure of the Lagrange spectrum using the dynamics of the bi-infinite shift space $\Sigma := \{1,2,3,\ldots\}^{\Z}$. That is, for $(a_{i})_{i\in\Z}\in\Sigma$ we define
\[\lambda_{0}((a_{i})_{i\in\Z}) := [a_{0};a_{1},a_{2},\ldots] + [0;a_{-1},a_{-2},\ldots],\]
and, for $j\in\Z$,
\[\lambda_{j}((a_{i})_{i\in\Z}) := \lambda_{0}(\sigma^{j}((a_{i})_{i\in\Z})) = \lambda_{0}((a_{i+j})_{i\in\Z}),\]
where $\sigma:\Sigma\to\Sigma$ is the left-shift sending $(a_{i})_{i\in\Z}$ to $(a_{i+1})_{i\in\Z}$. This shift map is related to the classical Gauss map
\[g([a_{0};a_{1},a_{2},\ldots]) = [a_{1};a_{2},a_{3},\ldots].\]
The Lagrange spectrum can then be equivalently defined as
\[L:= \{\ell(\underline{a}):= \limsup_{j\to\infty}\lambda_{j}(\underline{a})\,\mid\,\underline{a}\in\Sigma\}.\]
The Markov spectrum also permits a dynamical definition as
\[M:= \{m(\underline{a}):=\sup_{n\in\Z}\lambda_{n}(\underline{a}) \,\mid\, \underline{a}\in\Sigma\}.\]

Given an infinite sequence $(a_{i})_{i\in\Z}\in\Sigma$, we will often express it as $\ldots a_{-2}a_{-1}a_{0}^{*}a_{1}a_{2}\ldots$ where the asterisk denotes the 0th position. We will also use an overline to denote periodicity; e.g., $\overline{1^{*}23} = \ldots1231231^{*}23123123\ldots$. This notation should be clear from the context as we will typically restrict to the subshift $\{1,2,3,4\}^{\Z}$ so, in particular, all $a_{i}$ will be single digits. Moreover, when writing an equation of the form $t=m(\dots a_{-1}a_0^*a_1\dots)$, we always mean that the Markov value is being attained at the marked position, that is $m(\dots a_{-1}a_0^*a_1\dots)=\lambda_0(\dots a_{-1}a_0^{*}a_1\dots)$. With this notation, Freiman's gap can be described as
\[c_{F} = m(\overline{121313}22344^{*}3211\overline{313121}),\]
\[\nu_{F} = m(\overline{323444}313134^{*}313121133\overline{313121}).\]

Furthermore, given a finite word $\alpha\in(\N_{>0})^{\Z\cap [-b,a]}$ for some nonnegative integers $a, b$, inequalities of the form $\lambda_{0}(\alpha)>x$ will mean that we have $\lambda_{0}(w)>x$ for all bi-infinite sequences $w$ that are obtained by extending the finite sequence $\alpha$ on both sides.

Given an infinite sequence $\underline{a}=(a_i)_{i\in\Z}\in\Sigma$, the sequence $\underline{a}^T:=(a_{-i})_{i\in\Z}$ is called the \emph{transpose} of $\underline{a}$. Analogously, for finite words and for one-sided infinite sequences, that is for $a=a_1\dots a_n\in(\N_{>0})^n$ and $\underline{a}=a_0a_1\dots\in(\N_{>0})^\N$, the transpose are $a^T=a_n\dots a_1$ and $\underline{a}^T=\dots a_1a_0$, respectively.

\subsection{Local dimension}

Our first result investigates the third author's question about when the function $d_{loc}$ is non-decreasing.

\begin{theorem}\label{thm:dloc}
Consider the intervals
\begin{itemize}
\item $[3.05082, 3.122183)$
\item $[3.1299, 3.285441)$
\item $[3.28603, 3.28729)$
\item $[3.29296, 3.29335)$
\item $[3.33396, 3.33475)$
\item $[3.359, 3.423)$
\item $[\sqrt{12},3.8465)$
\item $[3.873,3.930691)$
\item $[3.93616, 3.943767)$
\item $[3.944054, 3.971606)$
\item $[3.97995, 3.9857)$
\item $[4.520781, 4.523103)$
\item $[4.5251, 4.5279)$.
\end{itemize}
For each interval $I$ above, $I\cap L'\neq\varnothing$ and for all $t\in I\cap L^\prime$ we have
\begin{equation*}
    \dloc(t)=\lim_{\epsilon\to 0}\dimH(L\cap(t-\epsilon,t+\epsilon))=\lim_{\epsilon\to 0}\dimH(M\cap(t-\epsilon,t+\epsilon))=d(t).
\end{equation*}
In particular $d_{loc}|_{I\cap L'}$ is non-decreasing.
\end{theorem}

To prove this, we will define the notion of a \emph{good interval} (see Definition~\ref{def:good_interval}). Such intervals $[\nu,\mu)$ will have the property that $d_{loc}$ is non-decreasing on $L'\cap[\nu,\mu)$. To prove that an interval is good, we will be required to understand how two shift spaces related to this interval are combinatorially related. Indeed, good intervals have strong transitivity properties for two naturally related subshifts. We analyse the properties of good intervals abstractly in Section~\ref{sec:goodprops}. We then prove that each of the above intervals are good intervals in Section~\ref{sec:good_ints}. These intervals are shown in blue in Figure~\ref{fig:goodints}.

\begin{figure}[ht!]
\centering
\begin{tikzpicture}[scale=7.5]

%Berstein conjecture
\fill[opacity = 0.4, red] (4.1,0.01)--(4.52,0.01)--(4.52,-0.01)--(4.1,-0.01)--cycle;

% Halls ray
\fill[opacity = 0.4, gray] (4.52782956616,0.015) --
(4.7,0.015) -- 
(4.7,-0.015) -- 
(4.52782956616,-0.015);

% Interval 1
\fill[opacity = 0.2, blue] (3.05082,-.03) -- (3.122183, -.03) -- (3.122183, .03) -- (3.05082,.03) -- cycle;

% Interval 2
\fill[opacity = 0.2, blue] (3.1299,-.03) -- (3.285441, -.03) -- (3.285441, .03) -- (3.1299,.03) -- cycle;

% Interval 3
\fill[opacity = 0.2, blue] (3.28603,-.03) -- (3.28729, -.03) -- (3.28729, .03) -- (3.28603,.03) -- cycle;

% Interval 3.5
\fill[opacity = 0.2, blue] (3.29296,-.03) -- (3.29335, -.03) -- (3.29335, .03) -- (3.29296,.03) -- cycle;

% Interval 3.75
\fill[opacity = 0.2, blue] (3.33396,-.03) -- (3.33475, -.03) -- (3.33475, .03) -- (3.33396,.03) -- cycle;

% Interval 4
\fill[opacity = 0.2, blue] (3.359,-.03) -- (3.423, -.03) -- (3.423, .03) -- (3.359,.03) -- cycle;

% Interval 5
\fill[opacity = 0.2, blue] (3.464,-.03) -- (3.84, -.03) -- (3.84, .03) -- (3.464,.03) -- cycle;

% Interval 6
\fill[opacity = 0.2, blue] (3.873,-.03) -- (3.930691, -.03) -- (3.930691, .03) -- (3.873,.03) -- cycle;

% Interval 7
\fill[opacity = 0.2, blue] (3.93616,-.03) -- (3.943767, -.03) -- (3.943767, .03) -- (3.93616,.03) -- cycle;

% Interval 8
\fill[opacity = 0.2, blue] (3.944054,-.03) -- (3.971606, -.03) -- (3.971606, .03) -- (3.944054,.03) -- cycle;

% Interval 9
\fill[opacity = 0.2, blue] (3.97995,-.03) -- (3.9857, -.03) -- (3.9857, .03) -- (3.97995,.03) -- cycle;

% Interval 10
\fill[opacity = 0.2, blue] (4.520781,-.03) -- (4.523103, -.03) -- (4.523103, .03) -- (4.520781,.03) -- cycle;

% Interval 11
(4.5251, 4.5279)
\fill[opacity = 0.2, blue] (4.5251,-.03) -- (4.5279, -.03) -- (4.5279, .03) -- (4.5251,.03) -- cycle;

% Line
\draw[-, thick] (2.6,0) -- (2.675,0);
\draw[thick, densely dotted] (2.675,0) -- (2.76,0);
\draw[->, thick] (2.76,0) -- (4.7,0);

% Cuts
\draw[thick] (3,0.02) -- (3,-0.02) node[below] {3};
\draw[thick] (3.334384,0.02) -- (3.334384,-0.02) node[below] {$t_1$};
\draw[thick] (4,0.02) -- (4,-0.02) node[below] {4};
\draw[thick] (4.52782956616,0.02) -- (4.52782956616,-0.02);

\draw[thick] (3.9420011599,-0.02) -- (3.9420011599,0.02) node[above] {$m_{1}$};

\node[below] at (4.52782956616,-0.03) {$c_{F}$};

% points
\foreach \x in {2.62,2.82842,2.97321,2.99605}
	\filldraw (\x, 0) circle[radius=0.15pt];
\node[below] at (2.62,0) {$\sqrt{5}$};
\node[below] at (2.82842,0) {$2\sqrt{2}$};

\end{tikzpicture}
\caption{The intervals of Theorem~\ref{thm:dloc} are depicted in blue. Hall's ray is depicted in grey. The interval $[4.1,4.52]$ of Berstein's conjecture is depicted in red. The largest known value of $M\setminus L$, denoted in the figure as $m_1$, is given by \Cref{thm:M-L}.}\label{fig:goodints}
\end{figure}

\subsection{Geometry of $\boldsymbol{M\setminus L}$}

Recall that a \emph{semi-symmetric word} is a finite word that is either a palindrome or a concatenation of two palindromes. As explained in \cite{Fractalgeometryofcomplement}, this concept is crucial for the construction of all known examples of $M\setminus L$. We will see in \Cref{subsec:berprops}, that restricted to good intervals, any element of $M\setminus L$ can be obtained by the construction outlined in \cite{Fractalgeometryofcomplement}.

So far all known regions of $M\setminus L$ have the structure depicted in Figure~\ref{fig:MminusL} where the spaces represent gaps in $M$, $j_1\in L^\prime$, $j_0$ is isolated in $M$, and $(M\setminus L)\cap(j_0,j_1)$ is given as a disjoint union of the form
\begin{equation*}
    (M\setminus L)\cap(j_0,j_1)=C\sqcup D
\end{equation*}
where $C$ is a Cantor set and $D$ is an infinite set of discrete points in $M$. The pre-images by the function $m(\underline a)$ of all elements of $(M\setminus L)\cap(j_0,j_1)$ connect in the past or in the future with some periodic not semi-symmetric word $\overline{w}$ of odd period. That is, the elements are given by sequences that (up to transposition) are eventually periodic on the left or on the right with period $w$. The fact that $j_1\in L^\prime$ for all known regions is a consequence of \cite[Theorem 2]{Flahive77}.

\begin{figure}[h]
\centering
\begin{tikzpicture}[scale=4, ]
\node[] at (-0.5,0) {(}; 
\node[] at (-0.01,0) {)}; 
\node[] at (0.01,0) {(}; 
\draw[dotted] (1.5,0) -- (0.5,0) node {)};
\draw[dotted] (0.5,0) -- (1.5,0) node {(};
\draw[-,dotted] (2.5,0) -- (2,0) node {)};
\node[] at  (0,-0.1) {$j_0\in L$};
\node[] at  (2,-0.1) {$j_1\in L^\prime$};
\node[] at  (0.5,-0.1) {$m_0$};
\node[] at  (1.5,-0.1) {$m_1$};
\node[] at  (1,-0.1) {$M\setminus L$};
\end{tikzpicture}
\caption{Structure of known regions of $M\setminus L$}\label{fig:MminusL}
\end{figure}

Since $(j_0,m_0)$ and $(m_1,j_1)$ are gaps of $M$, we know that the intersection $(M\setminus L)\cap(j_0,j_1)$ is a closed set. Recently, it was shown by Lima-Matheus-Moreira-Vieira in \cite{MminusLisnotclosed}  that $M\setminus L$ is not a closed set. More precisely they constructed a decreasing sequence $\widetilde{m}_{k}\in M\setminus L$  lying in distinct maximal gaps of $L$ such that $\lim_{k\to\infty}\widetilde{m}_{k}=1+\frac{3}{\sqrt{2}}\in L$. However, it is not known if there is a maximal gap $(\nu,\mu)$ of $L$ such that $(M\setminus L)\cap(\nu,\mu)$ is not closed. We do not yet have examples of elements in $L\cap\overline{(M\setminus L)}$ greater than $\sqrt{12}$, in other words such that the associated sequences contain at least one 3 or 4.

It is discussed in Cusick-Flahive~\cite[Chapter 3]{Cusick-Flahive}, that Berstein \cite{Berstein} was trying to find necessary and sufficient conditions for a Markov value to not be contained in the Lagrange spectrum. Berstein gave a list of conditions and constructed some intervals within which these conditions were both necessary and sufficient for a point $m \in M$ to be outside of $L$. 

In \Cref{subsec:berprops}, we will demonstrate that these conditions also hold for points $m\in(M\setminus L)\cap[\nu,\mu)$ for a good interval $[\nu,\mu)$. In particular, we will establish that any such element $m$ always connects to some non-semi-symmetric word $\overline{w}$. As a consequence, we conclude that all elements of $(M\setminus L) \cap [\nu,\mu)$ can be constructed using the method discovered by Freiman. This method and its relation with dynamical systems is very well explained in \cite{Fractalgeometryofcomplement}.

It is important to note that Berstein's intervals are constructed in an analogous way to our good intervals. That is, Berstein's proofs of the necessity and sufficiency of the conditions for $m\in M$ being in $L$ rely on proving that one theorem~\cite[Theorem 2]{Berstein} is true for each such interval. However, and importantly so, Berstein is not considering the transitivity of two subshifts related to the interval. These subshifts and their transitivity are crucial to our proofs on local dimension and perfectness of $M^\prime$.

In spite of the fact that the first digit of $\dimH(M\setminus L)$ is still unknown, we are able to characterize completely the Hausdorff dimension of $M\setminus L$ when restricted to good intervals.

\begin{theorem}
Let $[\nu,\mu)\subset\R$ be an interval. If $[\nu,\mu)$ is good, then \begin{enumerate}
\item If $(\ell_1,\ell_2)$ is a maximal gap of $L$ contained in $[\nu,\mu)$ and such that $M^\prime\cap(M\setminus L)\cap(\ell_1,\ell_2)\neq\varnothing$, then 
\begin{equation}\label{eq:precise_dimension}
    \dimH((M\setminus L)\cap(\ell_1,\ell_2))=D(\ell_1)=D(\ell_2).
\end{equation}
\item For all $t\in M^\prime\cap(M\setminus L)\cap[\nu,\mu)$ we have
\begin{equation*}
    \lim_{\varepsilon\to 0}\dimH((M\setminus L)\cap(t-\varepsilon,t+\varepsilon))=D(t).
\end{equation*}
\end{enumerate}
\end{theorem}

On the other hand, in recent work of the second and third authors with Matheus~\cite{JMM}, the largest known elements of $M\setminus L$ were found near to 3.938. Using \eqref{eq:precise_dimension} we can confirm rigorously the lower bound of \cite{JMM}: we have that $\dimH((M\setminus L)\cap(3.94,3.943))\geq D(3.94)>0.594561$ (see \Cref{subsec:lower_bound_dimM-L} for details).

The values of $M\setminus L$ mentioned in that paper were discovered via a computer-assisted investigation of gaps in the Lagrange spectrum suggested by numerical approximations in the work of Delecroix, Matheus and the third author~\cite[Figure 5]{DMM}. The algorithm used in the computer-assisted investigation is described in an appendix of the above work of Matheus and the second and third authors~\cite[Appendix B]{JMM}, but we also describe how the algorithm works in Subsection~\ref{subsec:M-Lalg}. In~\cite{JMM}, it is also mentioned that the computer investigations suggested that there should also be a region of $M\setminus L$ near to 3.942 (see~\cite[Appendix A]{JMM}). It was (heuristically) determined that new values of $M\setminus L$ in this range would not lead to an appreciable improvement on the lower bound of $\dimH(M\setminus L)$ and so this region was not investigated further in that work. Our second result is a full investigation of the structure of $M\setminus L$ in this region near 3.942.

Let
\begin{equation*}
    j_0=m(\overline{121112333^*11133232}) = 3.942001159911341469213548\dots\in L.
\end{equation*}

Let us denote the finite word $w=12111233311133232$ and $w^*=121112333^*11133232$. It is better to write the period of $j_0$ in this way because we have the same number of digits at each side of $3^*$. 

We have the following result.

\begin{theorem}\label{thm:M-L}
The intersection of $M\setminus L$ with $(3.942,3.943)$ is non-empty. The largest known element of $M\setminus L$ is
\begin{equation*}
m_1=\lambda_0(\overline{w}w^*w1211111\overline{23})=3.94200115991134146921437465\dots\approx j_0+8.26\cdot 10^{-22}.    
\end{equation*}
\end{theorem}

The full description of $M\setminus L$ in this region is given in \Cref{thm:characterization}. The value $m_{1}$ is also indicated on Figure~\ref{fig:goodints}.

We should highlight that in the computer approximations given by Delecroix-Matheus-Moreira~\cite[Figure 5]{DMM}, the last visible gap in the approximation of $L$ is near 3.942. 

In this direction, we prove the following result that demonstrates the existence of two new maximal gaps of $M$ near to this largest known value of $M\setminus L$. These intervals could also be contributing to the last visible gap in the computer approximations.

\begin{theorem}\label{thm:gaps}
Let
\[\nu_{1} = m(\overline{23}331113^*3\overline{32})=3.94254\dots,\]
\[\mu_{1} = m(\overline{23}3311133113212311333^*111331113331132123113311133\overline{32})=3.943304\dots,\]
\[\nu_{2} = m(\overline{12}31133311133111331113^*33113\overline{21})=3.94330534\dots, \]
and
\[\mu_{2} = m(\overline{21}11331113^*331132123113311133\overline{32})=3.94330716\dots.\]
The intervals $(\nu_{i},\mu_{i})$, $i = 1,2$, are maximal gaps of $M$; i.e., $(\nu_{i},\mu_{i})\cap M = \varnothing$ and $\nu_{i},\mu_{i}\in M$.
\end{theorem}

We prove this theorem in Section~\ref{sec:gaps}.

\subsection{$\boldsymbol{M'}$ and $\boldsymbol{M''}$}

In Section~\ref{sec:goodprops}, we will prove the following theorem addressing the question of whether $M' = M''$ in the setting of good intervals.

\begin{theorem}\label{thm:perfect}
Let $[\nu,\mu)\subset\R$ be an interval. If $[\nu,\mu)$ is good in the sense of Definition~\ref{def:good_interval}, then $M^\prime\cap[\nu,\mu)=M^{\prime\prime}\cap[\nu,\mu)$.
\end{theorem}

This will be proved as one part of Theorem~\ref{thm:sec2perfect}.

\subsection{Properties of $D(t)$}

Recall that $D(t) = \dimH(k^{-1}(-\infty,t])$. In~\cite{MV}, the third author and Villamil showed a concentration of dimension result, i.e., $D(t) = \dimH(k^{-1}(t))$, at points of the closure of the interior of the spectra, i.e., for $t\in\overline{\textrm{int}(L)} = \overline{\textrm{int}(M)}$. Moreover, they proved that $D(t)$ is strictly increasing across subintervals of the closure of the interior.

Using similar techniques and ideas, we will prove that the above results also hold on good intervals.  

\begin{proposition}
Let $[\nu,\mu)\subset\R$ be an interval. If $[\nu,\mu)$ is good, then for all $t\in L^\prime\cap[\nu,\mu)$ we have
\begin{equation*}
    \dimH (k^{-1}(t))=D(t).
\end{equation*}

Moreover if $s_1<t<s_2$ are such that $t\in L^\prime\cap[\nu,\mu)$, then
\begin{equation*}
    D(s_1)<D(s_2).
\end{equation*}
\end{proposition}

Recall that it is still a wide-open conjecture whether $\textrm{int}(L\cap(-\infty,c_F])\neq\varnothing$.

\subsection{Open questions}

We remind the reader that the questions of the third author discussed above remain open in general; i.e., outside of good intervals.

Motivated by Figure~\ref{fig:MminusL}, we also propose the following questions.

\begin{question}
Given a maximal gap $(\nu,\mu)$ of $L$, it is true that $M\cap (\nu,\mu)$ is closed in $\R$?
\end{question}

\begin{question}
Is $M^\prime\cap L=L^\prime$?
\end{question}

Recalling the maximal gap $(\nu_{F},c_{F})$ before Hall's ray, we also ask:

\begin{question}
Is $(M\setminus L)\cap (3.945, \nu_{F})\neq\varnothing$? Are there points of $M\setminus L$ close to $\nu_{F}$?
\end{question}

\begin{question}
What can be said about the structure of the sets $\partial L\cap[4.52,\nu_{F}]$ and $\partial M\cap [4.52,\nu_{F}]$? Are they uncountable?
\end{question}

\noindent {\bf Acknowledgements:} We would like to thank Carlos Matheus for mentioning the relation of good intervals with the previous work of Berstein. We thank the referee for their thorough and careful reading of the paper. For the purpose of open access, the authors have applied a Creative Commons Attribution (CC BY) licence to any Author Accepted Manuscript version arising from this submission.

\section{Monotonicity of local dimension and good intervals}\label{sec:goodprops}

In this section, we will investigate the local dimension function $d_{loc}:L'\to[0,1]$, define good intervals, and demonstrate that on such intervals $d_{loc}$ is non-decreasing.

\subsection{Good intervals}

Here we introduce the notion of a ``good interval." These intervals will be the settings in which we are able to address the questions of the third author discussed in the introduction. Berstein~\cite{Berstein} produced intervals with very similar properties in his study of $M\setminus L$. We discuss how our good intervals are related to Berstein's intervals in Subsection~\ref{subsec:berprops}.

The following two lemmas are from \cite{MMPV} and \cite{MMV} respectively. The first will be used heavily in Section~\ref{sec:good_ints} in order to bound Markov values.

Given a subshift  $\Sigma\subset(\N_{>0})^\Z$ we denote $\Sigma^+=\{(a_n)_{n\geq 0}:(a_n)_{n\in\Z}\in\Sigma\}$ and by $K(\Sigma)=\{[0;a_1,a_2,\dots]:(a_n)_{n\in\Z}\in\Sigma\}$. When $\Sigma$ is of finite type, $K(\Sigma)$ is called the stable Gauss-Cantor set associated to $\Sigma$. 

\begin{lemma} 
    \label{lem:tech}
    Let $\Sigma(C)$ be a transitive symmetric
    subshift.  
    Assume that three half-infinite sequences $v^1, v^2, v^3
    \in\Sigma^{+}(C)$ are such that 
    $[0;v^1]>[0;v^3]>[0;v^2]$. Then for all $\underline{a} \in \Sigma(C)$ 
    and for all $j \le n+1$
    \begin{multline*}
    \lambda_0\left(\sigma^j(\dots  a_{-2} a_{-1} a_0\dots
     a_{n} v^3)\right) \le 
    \max\left(m(\dots  a_{-2} a_{-1} a_0\dots  a_{n} v^1),
    m(\dots  a_{-2} a_{-1} a_0\dots  a_{n} v^2)\right) .
    \end{multline*} 
\end{lemma}

We number the sequences in the order above (i.e.,  $[0;v^1]>[0;v^3]>[0;v^2]$) as it resembles the proofs in Section~\ref{sec:good_ints} where we will bound a sequence $v_{B}$ between two sequences $v^{1}$ and $v^{2}$.

The second lemma will be used in Section~\ref{sec:good_ints} to determine the minimum Markov values given by sequences in certain subshifts. These minima will be the right endpoints of our good intervals.

\begin{lemma}
    \label{lem:minimizing} 
    Let $\beta$ be a finite word and $\theta$ be a symmetric finite word of even size on $\{1,2,\dots,A\}$. Assume that, under some conditions (e.g., after forbidding a finite list of finite strings), the Markov value of an infinite word of the type $\gamma=\omega_2^T \beta^T A\theta A\beta \omega_1$ is attained at one of the $A$ next to $\theta$, where $\omega_1, \omega_2$ are infinite words (on $\{1,2,\dots,A\}$). Then this Markov value is minimised when $\omega_1=\omega_2=:\omega$ and $[0;\beta,\omega]$ is minimal (under these conditions).
\end{lemma}

We will need a simple criterion to determine when two subshifts $\Sigma(B)$ and $\Sigma(C)$ are transitive.

Let $\CA=\{1,2,\dots,A\}$ be an alphabet. Denote by $\CA^*$ the set of finite words in the alphabet $\CA$. Given $w\in\CA^*$, we write $w=w^{-}a$ where $a\in\CA$ is the last letter of $w$. We will assume that $\CF$ is a finite set of words in $\CA^*$ that is symmetric, i.e. it contains all the transposes. Furthermore we assume that no word $w\in\CF$ is a subword of another word $\widetilde{w}\in\CF, w\neq\widetilde{w}$. Let
\[
\Sigma=\{\underline{a} \in 
\CA^{\mathbb Z} \mid \underline{a} \hbox{ has no substring from } \CF\}
\]
be a symmetric finite type subshift. 

Given $\underline{a}=(a_n)_{n\in\Z}\in\Sigma$ and $N\in\Z$ we denote $a_{-\infty,N}=\dots a_{N-2}a_{N-1} a_N$ and $a_{N,\infty}= a_{N}a_{N+1} a_{N+2}\dots$ the tails of $\underline{a}$.

\begin{proposition}\label{prop:mixing_criteria}
Suppose there is an $a\in\CA$ with $\overline{a}\in\Sigma$ and such that for all $w\in\CF$ there is a finite word $\tau_w\in\CA^\ast$ such that $w^{-}\tau_w\overline{a}$ does not contain words from $\CF$. Then the subshift $\Sigma$ is transitive.
\end{proposition}

\begin{proof}
Let $\underline{a}=(a_n)_{n\in\Z}\in\Sigma$ and $N\in\N_{>0}$. We want to find a finite $\tau$ in $\CA$ such that $a_{-\infty,N}\tau\overline{a}\in\Sigma$. Let $L$ be the length of the longest word of $\CF$. We start by considering $\underline{b}^{(1)}=a_{-\infty,N+L}\overline{a}$. If it does not contain words from $\CF$, then we are done, otherwise it must contain a word $w_1\in\CF$. That is, there are $n_{1}^{-},n_{1}^{+}\in\N$ with $N<n_{1}^{-}\leq N+L<n_{1}^{+}$ such that $w_1=b_{n_{1}^{-}}^{(1)}\dotsb b_{n_{1}^{+}}^{(1)}$. We take $w_1$ with $n_{1}^{-}$ minimal. By hypothesis there are $\tau_{w_1}$ a finite word in $\CA$ such that $w_1^{-}\tau_{w_1}\overline{a}$ does not contain words from $\CF$. Consider the bi--infinite sequence
\begin{equation*}
    \underline{b}^{(2)}=b^{(1)}_{-\infty, n_1^{+}-1}\tau_{w_1}\overline{a}=\dots a_{n_1^{-}-2}a_{n_1^{-}-1}w_1^{-}\tau_{w_1}\overline{a}.
\end{equation*}
If this bi--infinite sequence does not contain words from $\CF$ we are done, otherwise it contains a word $w_2\in\CF$, that is there are $n_2^{-},n_2^{+}\in\N$ such that $w_2=b_{n_{2}^{-}}^{(2)}\dotsb b_{n_{2}^{+}}^{(2)}$. Take $n_2^{-}$ minimal. We claim that this subword starts before $w_1^{-}$ and ends after $w_1^{-}$, that is $n_2^{-}<n_1^{-}<n_1^{+}\leq n_2^{+}$. Indeed, first note if $n_1^{-}\leq n_2^{-}$ then $w_2$ would be a subword of $b_{n_{1}^{-},\infty}^{(2)}=w_1^{-}\tau_{w_1}\overline{a}$ which contradicts the hypothesis. If $n_2^{+}< n_1^{+}$, then $w_2$ would be a subword of $b_{-\infty,n_{1}^{+}-1}^{(2)}=b_{-\infty,n_{1}^{+}-1}^{(1)}$, hence $w_2$ would be a forbidden subword of $\underline{b}^{(1)}$ starting at $n_2^{-}<n_1^{-}$, but this contradicts the minimality of $n_1^{-}$. 

Now that the claim is proved, in particular we obtain that $w_2$ is longer $|w_2|=n_2^{+}-n_2^{-}+1\geq n_1^{+}-n_1^{-}+2=|w_1|+1$. Now consider the bi--infinite sequence
\begin{equation*}
    \underline{b}^{(3)}=b^{(2)}_{-\infty, n_2^{+}-1}\tau_{w_2}\overline{a}=\dots a_{n_2^{-}-2}a_{n_2^{-}-1}w_2^{-}\tau_{w_2}\overline{a}.
\end{equation*}

Inductively we find a sequence of possible continuations $\underline{b}^{(1)},\dots,\underline{b}^{(r)}$ each one associated with a word $|w_i|\geq|w_{i-1}|+1$ and such that $n_i^{-}<n_{i-1}^{-}<n_{i-1}^{+}\leq n_{i}^{+}$. Since $\CF$ is finite, this algorithm must finish for some $r\leq L$ and so 
\begin{equation*}
    \underline{b}^{(r+1)}=b^{(r)}_{-\infty, n_r^{+}-1}\tau_{w_r}\overline{a}=\dots a_{n_r^{-}-2}a_{n_r^{-}-1}w_r^{-}\tau_{w_r}\overline{a}\in\Sigma.
\end{equation*}
Finally this completes our initial aim because $n_r^{-}\geq n_r^{+}-L\geq n_1^{+}-L>N$.    
\end{proof}

Now we introduce the notion of when an orbit connects to a transitive subshift of finite type. The first definition is from \cite[Definition 3.2]{MMPV} and the second one is reminiscent of \cite[Definition 3.1]{MV}. 

\begin{definition}
    \label{def:cplus}
    Consider two transitive and symmetric subshifts of finite type $\Sigma(B)
    \subset \Sigma(C)$. 
     Let $\underline{a}\in \Sigma(C)$ be a sequence with
$m(\underline{a}) = \lambda_0(\underline{a}) = m \in M$. 
\begin{itemize}
\item We say that $\underline{a}$ {\it
connects positively to~$B$} if, for every $k\in\mathbb{N}$, there exists
a finite sequence~$\tau$ and an infinite sequence $ v\in\Sigma^{+}(B)$ such that for
$\widetilde{\underline{a}}:=\dots  a_{-2} a_{-1} a_0\dots  a_{k}
\tau v$ we have 
\begin{equation}
    \label{eq:m<m}
m(\widetilde{\underline{a}})<m(\underline{a})+2^{-k}.
\end{equation}
We say that $\underline{a}$
{\it connects negatively} to~$B$ if the reversed sequence $\underline{a}^T$ connects
positively to~$B$. 
\item We say that $\underline{a}$ {\it
connects positively to~$B$ before $t$} if, for every $k\in\mathbb{N}$, there exists
a finite sequence~$\tau$ and an infinite sequence $ v\in\Sigma^{+}(B)$ such that for
$\widetilde{\underline{a}}:=\dots  a_{-2} a_{-1} a_0\dots  a_{k}
\tau v$ we have 
\begin{equation*}
m(\widetilde{\underline{a}})<t+2^{-k}.
\end{equation*}
We say that $\underline{a}$
{\it connects negatively to~$B$ before $t$} if the reversed sequence $\underline{a}^T$ connects
positively to~$B$ before $t$. 
\end{itemize}
\end{definition}

\begin{remark}
    \label{rem:MtoL}
    Assume that $\Sigma(B) \subset \Sigma(C)$ are two transitive symmetric
    subshifts and let $x$ be such that $m(\underline{b}) \le x$ for all $\underline{b} \in
    \Sigma(B)$. Consider a sequence $\underline{c} \in \Sigma(C)$ with $m(\underline{c}) =
    \lambda_0(\underline{c}) = m > x$. 
    Then for any finite sequence $\tau$ and half-infinite sequence $ v \in
    \Sigma^{+}(B)$ directly from definition of Lagrange and Markov
    numbers we get
$$ 
 \limsup_{j\to+\infty}\lambda_0(\sigma^j(\dots \gamma_{-2}\gamma_{-1}\gamma_0\dots
\gamma_{k} \tau v))<m,
$$
Thus, if we want to get that $m(\dots \gamma_{-2}\gamma_{-1}\gamma_0\dots
\gamma_{k} \tau v)<m+2^{-k}$, then it suffices to check that  
\begin{equation*}
\lambda_0(\sigma^j(\dots \gamma_{-2}\gamma_{-1}\gamma_0\dots \gamma_{k}
\tau v))\!<\!m+2^{-k} 
\end{equation*}
for finitely many values of $j$, namely, for all $0\le j\le k+|\tau|+l$ where $l$ is 
sufficiently large (so that $2^{1-l}\le m-x+2^{-k}$). This is because, for two continued fractions $\alpha = [0;a_{1},a_{2},\ldots]$ and $\beta = [0;b_{1},b_{2},\ldots]$ with $a_{i} = b_{i}, 1\leq i\leq N$, and $a_{N+1}\neq b_{N+1}$, we have $|\alpha - \beta|<2^{1-N}$.
\end{remark}

The important consequence of connecting to both sides to a transitive subshift with smaller Markov values is given by \cite[Lemma 3.4]{MMPV}. 

\begin{lemma}\label{lem:key} Consider two transitive and symmetric subshifts of finite type $\Sigma(B)\subset \Sigma(C)$. Let~$x$ be such that $m(\underline{b}) \le x$ for all $\underline{b} \in \Sigma(B)$. Suppose that a sequence $\gamma \in\Sigma(C)$ satisfies $m(\gamma) =\lambda_0(\gamma) = m \geq x$ and connects positively and negatively to~$B$. Then $m \in L$. 
\end{lemma}

\begin{remark}\label{rem:sup_transitive_subshifts}
In general, if $\Sigma$ is a transitive symmetric subshift with infinite cardinality then $\sup\{m(\gamma):\gamma\in\Sigma\}\in L^\prime$.
\end{remark}

In this direction, we will associate to a suitable interval $[x,y)$ two symmetric transitive subshifts of finite type $\Sigma(B) = \Sigma(B_x)\subset
\Sigma(C) = \Sigma(C_y)\subset(\mathbb{N}^*)^{\mathbb{Z}}$ such that 
\begin{itemize} 
    \item $m(\underline{b})\leq x$ for all $\underline{b} \in \Sigma( B)$; and
    \item for all $\underline{c}\in(\mathbb{N}^*)^{\mathbb{Z}}$ such that $m(\underline{c})<y$
        we have $\underline{c} \in \Sigma(C)$.  
\end{itemize}
We will always choose $x:=\sup\{m(\underline{b}):\underline{b}\in\Sigma(B)\}$. Note that if $\Sigma(B)$ has infinite cardinality then Remark~\ref{rem:sup_transitive_subshifts} gives us that $x\in L'$.

\begin{definition}\label{def:good_interval}
Let $[\nu,\mu)$ be an interval. We say that $[\nu,\mu)$ is good if it can be covered by intervals $[x,y)$ as above such that for any $m(\underline{a})=\lambda_0(\underline{a})<y$ it holds: given any position $N\in\N$ such that there are two continuations $v^1, v^2 \in
\Sigma^+(C)$ with \emph{different} first term, then there are a (perhaps empty) finite word $\tau$ and $v_B\in\Sigma^+(B)$ such that \[m(\dots a_{-1}\alpha_0\dots a_N\tau v_B)\leq \max\{m(\dots a_{-1}a_0\dots a_N\upsilon^1),m(\dots a_{-1}a_0\dots a_N\upsilon^2),x\}.\]     
\end{definition}

Note that a connected finite union of good intervals is again a good interval.

\subsection{Local dimension on good intervals}

Define $d(t):\R\to[0,1]$ by
\begin{equation*}
    d(t)=\dimH(L\cap(-\infty,t))=\dimH(M\cap(-\infty,t)).
\end{equation*}

    Recall that Moreira~\cite{geometricproperties}  proved the following formula:
    \[
        d(t)=\min\{1,2\cdot D(t)\},
    \]
    where $D(t)=\dimH(K_t)=\overline{\dimB}(K_t)$ is a continuous function and 
    \begin{multline*}
        K_t=\{[0;c_1,\dots,c_n,\dots]\ \mid\ \text{there exists $(c_{-n})_{n\geq 0}\in(\N_{>0})^{\N}$ such that} \\ 
 [c_k;c_{k+1},\dots,]+[0;c_{k-1},c_{k-2},\dots]\leq t,\forall k\in\Z\}.
\end{multline*}
Note that, using the notation from the introduction,
\[D(t) = \dimH(k^{-1}(-\infty,t]) = \dimH(k^{-1}(-\infty,t)).\]

Define the function $\dloc:L^\prime\to[0,1]$
\begin{equation*}
    \dloc(t)=\lim_{\varepsilon\to 0}\dimH(L\cap(t-\varepsilon,t+\varepsilon)).
\end{equation*}

By using the fact that $d$ is continuous and the monotonicity of the Hausdorff dimension, we always have $\dloc\leq d$.

\begin{proposition}\label{prop:conjecture_characterization}
The function $\dloc:L^\prime\to[0,1]$ is non-decreasing if and only if $\dloc=d$ on $L^\prime$.
\end{proposition}
\begin{proof}
Suppose $\dloc$ is non-decreasing. If $\dloc(t)<d(t)$ for some $t$, let $t_0=\min\{s:d(s)=d(t)\}$. By definition we have that $d(t_0-\varepsilon)<d(t_0)$ for all $\varepsilon>0$, so we must have $\dloc(t_0)=d(t_0)$. However since $t_0\leq t$ by hypothesis  $d(t)=d(t_0)=\dloc(t_0)\leq \dloc(t)<d(t)$, a contradiction. 
\end{proof}

 Recall that the notation $\Sigma_{t}$ for some $t\in\R$ denotes the set
\[\Sigma_{t} := \{\underline{a}\in(\N^{*})^{\Z}\,|\,m(\underline{a})\leq t\}.\]

\begin{theorem}
Let $[\nu,\mu)\subset\R$ be an interval. If $[\nu,\mu)$ is good, then $\dloc(t)=d(t)$ for all $t\in L^\prime\cap[\nu,\mu)$. In particular $\dloc$ is non-decreasing in $[\nu,\mu)$. Moreover for all $t\in L^\prime\cap[\nu,\mu)$ we have
\begin{equation}\label{eq:localdimM}
    \dloc(t)=\lim_{\epsilon\to 0}\dimH(L\cap(t-\epsilon,t+\epsilon))=\lim_{\epsilon\to 0}\dimH(M\cap(t-\epsilon,t+\epsilon))=d(t).
\end{equation}

\end{theorem}

\begin{proof}
Let $t\in L^\prime\cap[x,y)$ with $x,y$ as in \Cref{def:good_interval}. We want to show $\dloc(t)\geq d(t)$ by connecting $\Sigma_{s}\to\Sigma(B)\to\Sigma_{t^\prime}\to\Sigma(B)\to\Sigma_{s}$ for values of $s\in M$, $s<t$ close to $t$ and values $t^\prime\in L$ close to $t$. In other terms, given $\varepsilon>0$ small, we would like to find a sequence $\gamma=(\gamma_i)_{i\in\Z}$ with $m(\gamma)=\lambda_0(\gamma)=t^\prime$ close to $t$, finite subwords $\beta_r,\beta_l$ of some sequences in $\Sigma(B)$, finite words $\underline{\tau},\tilde{\tau},\tau,\underline{\tilde{\tau}}$ and a complete subshift $\Sigma(A)\subset\Sigma_s$ with $\dimH(K(A))$ close to $\dimH(K_s)$, such that
\begin{equation*}
    \lvert{m(v_l^T\underline{\tau}\beta_l\tilde{\tau}\gamma_{-k}\dots\gamma_k\tau\beta_r\underline{\tilde{\tau}} v_r)-t\rvert}<\varepsilon, \qquad \text{for all } v_l,v_r\in\Sigma^+(A). 
\end{equation*}

More specifically, we will prove the following: for all $\varepsilon>0$ it holds
\begin{equation}\label{eq:2}
    \dim_H(L\cap(t-\varepsilon, t+\varepsilon))\geq d(t)-\widetilde{\delta},
\end{equation}
where $\widetilde{\delta}\to0$ as $\varepsilon\to0$. In particular, since $L\subset M$ one has
\begin{equation*}
    \dimH(L\cap(t-\epsilon,t+\epsilon))\leq \dimH(M\cap(t-\epsilon,t+\epsilon))\leq d(t+\epsilon)
\end{equation*}
so letting $\varepsilon\to0$ in \eqref{eq:2} will show \eqref{eq:localdimM}.

Now we want to show \eqref{eq:2}. Since $t\in L^\prime$, by the proof of \cite[Theorem 3]{geometricproperties} we know that for any $\delta^\prime>0$ there are regular Cantor sets $K(A_1),K(A_2)$ defined by iterates of the Gauss map such that $\diam\left(t-\big(K(A_1)+K(A_2)\big)\right)\leq\delta^\prime$ and $K(A_1)+K(A_2)\subset L$. We claim that for all $m(\gamma)=\lambda_0(\gamma)\in K(A_1)+K(A_2)$ we have that $\gamma$ connects positively and negatively with $B$ before $\max\{m(\gamma),x\}$. Indeed, given any $k\in\N$, since $K(A_2)$ is a regular Cantor set, we can find some large $N\in\N, k\leq N$ and continuations $ v^1, v^2$ with different first term and such that $m(\dots\gamma_{N} v^i)<m(\gamma)+2^{-N}$ for $i=1,2$. Since the interval is good we know there are a finite word $\tilde\tau$ and $v_B\in\Sigma^+(B)$ such that $m(\dots\gamma_{N}\tilde\tau v_B)<\max\{m(\gamma),x\}+2^{-N}$. So by taking $\tau=\gamma_{k+1}\dots\gamma_N\tilde\tau$ and $v=v_B\in\Sigma^+(B)$ we are done. The same proof applies to $\gamma^T$ and $K(A_1)$. If $t>x$ we impose that $0<\delta^\prime<t-x$ and if $t=x$ we choose $\gamma\in\Sigma(B)$ such that $\lambda_0(\gamma)=m(\gamma)=x$.

We will use the following facts that hold for any $t>3$. By \cite[Lemma 2]{geometricproperties}, for any $\eta\in(0,1)$ there is $\delta>0$ and a complete shift $\Sigma(A)$, where $A$ is a finite set of finite words such that $\Sigma(A)\subset\Sigma_{t-\delta}$ and $\dimH(K(A))>(1-\eta)\dimH(K_t)$. Suppose that the maximum of $m(\Sigma(A))$ is attained at
\begin{equation}\label{eq:3}
    \underline{\theta}=(\dots,\widetilde{\alpha}_1,\widetilde{\alpha}_0,\widetilde{\alpha}_1,\dots), \quad\widetilde{\alpha}_i\in A,\forall i\in\Z.
\end{equation}
in a position belonging to the word $\widetilde{\alpha}_0$. In particular $m(\underline{\theta})\leq t-\delta$. By the same proof given above we have that $\underline{\theta}$ connects positively and negatively with $B$ before $\max\{t-\delta,x\}$. 

Since $\underline{\theta}$ connects positively and negatively to~$B$ before $\max\{t-\delta,x\}$, given any $k\in\N_{>0}$, there exist finite
sequences $\underline{\tau}, \underline{\widetilde\tau}$ and infinite sequences $\underline{ v}, \underline{\widetilde v} \in\Sigma^+(B)$ such that 
\begin{multline}\label{eq:123}
m(\dots \widetilde{\alpha}_{-2}\widetilde{\alpha}_{-1}\widetilde{\alpha}_0\dots
\widetilde{\alpha}_{k} \underline{\tau}\underline{ v})<\max\{t-\delta,x\}+2^{-k}\mbox{  and  }\\ m(\underline{\widetilde  v}^T
\underline{\widetilde\tau} \widetilde{\alpha}_{-k}\dots \widetilde{\alpha}_0 \widetilde{\alpha}_1 \widetilde{\alpha}_2 \dots)<\max\{t-\delta,x\}+2^{-k}.    
\end{multline}
Similarly, since $\gamma$ connects positively and negatively to~$B$, there exist finite
sequences $\tau, \widetilde\tau$ and infinite sequences $ v, \widetilde v \in\Sigma^+(B)$ such that 

\begin{equation}\label{eq:4}
m(\dots \gamma_{-2}\gamma_{-1}\gamma_0\dots
\gamma_{k} \tau v)<m(\gamma)+2^{-k}\mbox{  and  } m(\widetilde  v^T
\widetilde\tau \gamma_{-k}\dots \gamma_0 \gamma_1 \gamma_2 \dots)<m(\gamma)+2^{-k}.
\end{equation}

In the particular case where $m(\gamma)=\lambda_0(\gamma)=x$ we connect to $B$ just using $\tau\upsilon=\gamma_{k+1}\gamma_{k+2}\dots$ and $\tilde{\upsilon}^T\tilde{\tau}=\dots\gamma_{-k-2}\gamma_{-k-1}$.

The Markov value of either of the continuations \eqref{eq:4} is attained at the position $n=0$ if $m(\gamma)=\lambda_0(\gamma)=x$. Otherwise $m(\gamma)=\lambda_0(\gamma)>x$ and the Markov value of the leftmost, resp. rightmost continuation in \eqref{eq:4} is attained at some position $n\leq|\tau|+|\widetilde{\tau}|+k+l$, resp. $n\geq-|\tau|-|\widetilde{\tau}|-k-l$, where $l$ is such that $2^{1-l}< m(\gamma)-x+2^{-k}$. In particular
\begin{equation}\label{eq:first_connection}
    |m(\widetilde{ v}^T\widetilde{\tau}\gamma_{-k}\dots\gamma_{-1}\gamma_0\gamma_1\dots\gamma_k\tau v)-m(\gamma)|<2^{1-k},
\end{equation}
and the Markov value is attained in a position $|n|\leq|\tau|+|\widetilde{\tau}|+k+l$.

Let $ v^{l+k}:= v_1\dots v_{l+k}$, 
$\widetilde v^{l+k}: = \widetilde  v_{l+k} \dots \widetilde
 v_1$, $\underline{ v}^{l+k}:=\underline{ v}_1\dots\underline{ v}_{l+k}$, 
$\underline{\widetilde v}^{l+k}: = \underline{\widetilde  v}_{l+k} \dots \underline{\widetilde
 v}_1$ be the initial segments of $ v,\widetilde v^T,\underline{ v},\underline{\widetilde v}^T$ respectively. By transitivity of $\Sigma(B)$, there exists $\beta \in \Sigma(B)$ which contains non-overlapping occurrences of the strings $\underline{ v}^{l+k}$ and $\widetilde v^{l+k}$ in this order and analogously the same holds for $ v^{l+k}$ and $\underline{\widetilde v}^{l+k}$. Let us denote by $(\underline{ v}^{l+k}\ast\widetilde v^{l+k})$ a finite substring of $\beta$ which begins with $\underline{ v}^{l+k}$ and terminates with $\widetilde  v^{l+k}$ and analogously define $( v^{l+k}\ast\underline{\widetilde{ v}}^{l+k})$.

Finally, using all the above we obtain that for $j\geq k$
\begin{multline}\label{eq:all_connected}
    m(\dots \alpha_{-j-2} \alpha_{-j-1}\widetilde{\alpha}_{-j}\dots\widetilde{\alpha}_0\dots\widetilde{\alpha}_k\underline{\tau}(\underline{ v}^{l+k}\ast\widetilde{ v}^{l+k})\widetilde{\tau}\gamma_{-k}\dots\gamma_0\dots \\
    \dots\gamma_{k}\tau( v^{l+k}\ast\underline{\widetilde{ v}}^{l+k})\underline{\widetilde{\tau}}\widetilde{\alpha}_k\dots\widetilde{\alpha}_0\dots\widetilde{\alpha}_j \alpha_{j+1} \alpha_{j+2}\dots)<t+2^{3-k}
\end{multline}
for any choice of $ \alpha_{|i|}\in A$ for $|i|\geq j+1$, where we used \eqref{eq:first_connection} and $|m(\gamma)-t|=|\lambda_0(\gamma)-t|<\delta^\prime<\delta$. By the same reason the Markov value of the above bi-infinite sequence is attained in some position $|n|\leq|\tau|+|\widetilde{\tau}|+k+l$ inside $(\underline{ v}^{l+k}\ast\widetilde{ v}^{l+k})\widetilde{\tau}\gamma_{-k}\dots\gamma_0^*\dots\gamma_{k}\tau( v^{l+k}\ast\underline{\widetilde{ v}}^{l+k})$ where $\gamma_0^*$ denotes the zero position. By \eqref{eq:first_connection}, we will also have that the above Markov value is at least $t-\delta^\prime-2^{3-k}$. 

Let us choose $\delta^\prime=\delta/2$ and $j=k$ with $k$ such that $2^{5-k}<\delta$. We claim that in fact, given any fixed choice of the $\alpha_{|i|}$ for $|i|\geq k+1$, the Markov value of \eqref{eq:all_connected} belongs to $L$. We will use the fact that Markov values of periodic sequences are dense in $L$ (see \cite[Theorem 2, Chapter 3]{Cusick-Flahive}). Since the Markov values of the complete shift $\Sigma(A)$ are less than or equal to $t-\delta$, and since the Markov value in \eqref{eq:all_connected} is at least $t-\delta/2-2^{3-k}>t-\delta+2^{2-k}$, we have that for any $n> k$ the Markov value of the bi-infinite periodic sequence with period
\begin{equation*}
    \alpha_{-n} \dots \alpha_{-k-1}\widetilde{\alpha}_{-k}\dots\widetilde{\alpha}_0\dots\widetilde{\alpha}_k\underline{\tau}(\underline{ v}^{l+k}\ast\widetilde{ v}^{l+k})\widetilde{\tau}\gamma_{-k}\dots\gamma_0\dots\gamma_{k}\tau( v^{l+k}\ast\underline{\widetilde{ v}}^{l+k})\underline{\widetilde{\tau}}\widetilde{\alpha}_k\dots\widetilde{\alpha}_0\dots\widetilde{\alpha}_k \alpha_{k+1} \dots \alpha_{n}
\end{equation*}
is again being attained inside $(\underline{ v}^{l+k}\ast\widetilde{ v}^{l+k})\widetilde{\tau}\gamma_{-k}\dots\gamma_0\dots\gamma_{k}\tau( v^{l+k}\ast\underline{\widetilde{ v}}^{l+k})$. Therefore these periodic Markov values are converging to the Markov value of \eqref{eq:all_connected} and thus it belongs to $L$.

By all of the above, we finally obtain, using the Markov values in \eqref{eq:all_connected}, that $L\cap(t-\varepsilon,t+\varepsilon)$ contains a sum of the form $K+K^\prime$ where both $K$ and $K^\prime$ are Gauss-Cantor sets that are bi-Lipschitz equivalent to $K(A)$. By using Moreira's dimension formula \cite{CartesianProduct} we get
\begin{align*}
    \dimH(L\cap(t-\varepsilon,t+\varepsilon))\geq\dimH(K+K^\prime)&=\min\{2\dimH(K(A)),1\} \\
    &\geq\min\{2(1-\eta)\dimH(K_t),1\}\geq d(t)-\widetilde{\delta}.
\end{align*}
\end{proof}

\subsection{The topology and geometry of $\boldsymbol{M'}$ and the Hausdorff dimension of $M\setminus L$}

Here we discuss the topological and geometric properties of $M'$ in the setting of good intervals.

As discussed in the introduction, an open question asked by the third author and concerning the topology of the spectra is whether the set $M^\prime$ is perfect; that is, is $M^\prime=M^{\prime\prime}$?

We will prove in the following theorem that these two sets do coincide within good intervals. 

\begin{theorem}\label{thm:sec2perfect}
Let $[\nu,\mu)\subset\R$ be an interval. If $[\nu,\mu)$ is good, then \begin{enumerate}
\item If $(\ell_1,\ell_2)$ is a maximal gap of $L$ contained in $[\nu,\mu)$ and such that $M^\prime\cap(M\setminus L)\cap(\ell_1,\ell_2)\neq\varnothing$, then 
\begin{equation}\label{eq:regions_dimension}
    \dimH((M\setminus L)\cap(\ell_1,\ell_2))=D(\ell_1)=D(\ell_2).
\end{equation}
\item For all $t\in M^\prime\cap(M\setminus L)\cap[\nu,\mu)$ we have
\begin{equation}\label{eq:MminusL_has_lowdim}
    \lim_{\varepsilon\to 0}\dimH((M\setminus L)\cap(t-\varepsilon,t+\varepsilon))=\dimH(K_t)=D(t).
\end{equation}
\item We have
\begin{equation}\label{eq:Mprimeperfect}
    M^\prime\cap[\nu,\mu)=M^{\prime\prime}\cap[\nu,\mu).
\end{equation}

\end{enumerate}
\end{theorem}

\begin{proof}
Let $t\in (M\setminus L)\cap(x,y)$ with $x,y$ as in \Cref{def:good_interval} (note that $x\in L^\prime$ by \Cref{rem:sup_transitive_subshifts}). We first show that for all
$\varepsilon>0$ small enough
\begin{equation}\label{eq:1}
    \dimH\left((M\setminus L)\cap(t-\varepsilon, t+\varepsilon)\right)\leq\dimH(K_{t+\varepsilon}),
\end{equation}

We know that $\{\underline{a}\in(\N_{>0})^\Z:m(\underline{a})=\lambda_0(\underline{a})=t\}$ is finite: this was proved by
Bernstein \cite[Theorem 3, Page 25]{Berstein} and, for the sake of completeness, we give a short proof in \Cref{prop:Berstein_properties} below. Let $\{\underline{a}^{(1)},\dots,\underline{a}^{(r)}\}$ be the sequences such that $m(\underline{a}^{(i)})=\lambda_0(\underline{a}^{(i)})=t$. We claim that for any $n\in\N$, there is $\varepsilon>0$ such that if $|m(\underline{a})-t|<\varepsilon$ with $m(\underline{a})=\lambda_0(\underline{a})$, then $a_{-n,n}=a^{(j)}_{-n,n}$ for some $1\leq j\leq r$. Indeed, by contradiction there will be a sequence $\underline{b}^{(i)}$ with $m(\underline{b}^{(i)})=\lambda_0(\underline{b}^{(i)})$, $|m(\underline{b}^{(i)})-t|<1/i$ but $b^{(i)}_{-n,n}\neq a^{(j)}_{-n,n}$ for all $i,j$. Therefore, passing to a subsequence of indexes, we conclude that the $\underline{b}^{(i)}$ converge to a sequence $\underline{b}$ with $m(\underline{b})=\lambda_0(\underline{b})=t$ but with $b_{-n,n}\neq a^{(j)}_{-n,n}$ for all $1\leq j\leq r$, a contradiction. 

As a consequence, for any $\varepsilon>0$ small enough, it holds that given any $\underline{a}$ such that $|m(\underline{a})-t|<\varepsilon$ with $m(\underline{a})=\lambda_0(\underline{a})$, we will have that either $a_{0,\infty}=a^{(j)}_{0,\infty}$ or $a_{-\infty,0}=a^{(j)}_{-\infty,0}$ for some $1\leq j\leq r$. To see this, we can assume, without loss of generality, that each $\underline{a}^{(j)}$ does not connect positively to $B$. Thus there is a $k\in\N$, such that for any $N\geq k$ we will have $m(a^{(j)}_{-\infty,N}\tau v_B)\geq t+2^{-k}$ for any finite word $\tau$, any $v_B\in\Sigma^+(B)$ and all $1\leq j\leq r$. Suppose $\varepsilon<2^{-k}$ is small enough to guarantee that $a_{-N,N}=a^{(j)}_{-N,N}$ for some $j$ (up to transposition of $\underline{a}$). Hence if $N_1\geq N$ is such that $a_{N_1}\neq a^{(j)}_{N_1}$, since the interval is good there is a finite $\tau$ and $v_B\in\Sigma^+(B)$ such that
\begin{equation*}
    m(a^{(j)}_{-\infty,N_1-1}\tau v_B)\leq\max\{m(\underline{a}^{(j)}),m(\underline{a})\}<t+\varepsilon<t+2^{-k}
\end{equation*}
which is impossible.

Given $\varepsilon>0$ small, by the previous argument we know that up to transposition there is a finite set $X$ of possible continuations $\upsilon$ and $N\in\N$ such that given any $\lambda_0(\gamma)=m(\gamma)\in (M\setminus L)\cap(t-\varepsilon,t+\varepsilon)$, we have $\gamma_{N}\gamma_{N+1}\dots=\upsilon$ for some $\upsilon\in X$. In particular, we will have that 
\begin{equation*}
    (M\setminus L)\cap(t-\varepsilon,t+\varepsilon)\subset\bigcup_{\upsilon\in X}\bigcup_{\substack{1\leq a\leq T \\ i=1,2,\dots}}(a+g^i([0;\upsilon])+K_{t+\varepsilon}),
\end{equation*}
where $g$ is the Gauss map.

Since the Hausdorff dimension of a countable union is the supremum of the Hausdorff dimensions of each set, we get $\dimH\left((M\setminus L)\cap(t-\varepsilon,t+\varepsilon)\right)\leq\dimH(K_{t+\varepsilon})$ for all $\varepsilon$ small enough. This finishes the proof of \eqref{eq:1}.

Since $D(t)=\dimH(K_t)$ is constant on gaps of $L$ (because $D(t)=\dimH(k^{-1}(-\infty,t))=\dimH(k^{-1}(-\infty,t])$), observe that \eqref{eq:1} already gives one side of the equalities in \eqref{eq:regions_dimension} and \eqref{eq:MminusL_has_lowdim}. To prove the converse inequalities, the idea is to connect $\Sigma_{t^\prime}\to\Sigma(B)\to\Sigma_{s}$ where $s<t$ and $t^\prime$ are very close to $t$. More precisely, for any $\eta\in(0,1)$ and $\varepsilon>0$, we will show that there is a Gauss-Cantor set $K$ such that $K\subset M^\prime$, $\diam(t-K)<\varepsilon$ and $\dimH(K)>(1-\eta)\dimH(K_t)$. In particular, if $t\in M^\prime\cap[\nu,\mu)$ then $t\in M^{\prime\prime}$ as well, and if $t\in M^\prime\cap(M\setminus L)\cap[\nu,\mu)$ then $K\subset (M\setminus L)\cap(t-\varepsilon,t+\varepsilon)$ for $\varepsilon$ small.

Let $t \in M^\prime\cap[x,y)$ with $x,y$ as in \Cref{def:good_interval}. Since we already know that $x\in L^\prime=L^{\prime\prime}$, we can assume that $t\in(x,y)$. Consider a sequence $t_n$ converging to $t$, $t_n \in M$, $t_n \ne t$. Choose $\underline{\theta}^{(n)} \in \Sigma$
such that $t_n = m(\underline{\theta}^{(n)})=\lambda_0(\underline{\theta}^{(n)})$. Let $\underline{\theta}^{(n)}=(b_\ell^{(n)})_{\ell\in\Z}$ and assume $b_\ell^{(n)} \le 4$, $\forall\, \ell$, $\forall\, n$ (which is possible since we may assume that the $t_n$ are not in Hall's ray). Given $\varepsilon>0$, there exists $n_0\in\N$ large such that $n\geq n_0\implies|m(\underline{\theta}^{(n)})-t|<\varepsilon/2$. Let $N =
\lceil \varepsilon^{-1} \rceil+4$. There is a sequence $S$ such that for infinitely many values of $n$, $S$ appears infinitely many times as $(b_{-N}^{(n)},b_{-N+1}^{(n)},\dots,b_{N}^{(n)})$, i.e. there are $n_0<n_1<n_2<\dots$ such that $S=(b_{-N}^{(n_i)},b_{-N+1}^{(n_i)},\dots,b_{N}^{(n_i)})$, $\forall i\geq 1$. Since $t_{n_i}\neq t$ and $t_{n_i}\to t$, we have (up to transposition) that for some $N_1\geq N$ there are two continuations $ v^1, v^2$ with different first term such that $ v^1=(b_{N_1+1}^{(n_i)},b_{N_1+2}^{(n_i)},\dots)$ and $ v^2=(b_{N_1+1}^{(n_j)},b_{N_1+2}^{(n_j)},\dots)$ for some $n_i\neq n_j$ and $(b_0^{(n_i)},\dots,b_{N_1}^{(n_i)})=(b_0^{(n_j)},\dots,b_{N_1}^{(n_j)})$. Define, for some $h\in\Z$, $\underline{\theta}^{(n_i)}_{-\infty,h}=\dots b^{(n_i)}_{h-2}b^{(n_i)}_{h-1} b^{(n_i)}_h$ the left tail of $\underline{\theta}^{(n_i)}$.  Since the interval is good, there is a finite $\tau$ and continuation $v_B\in\Sigma^+(B)$ such that 
\begin{equation}\label{eq:5}
    |m(\underline{\theta}^{(n_i)}_{-\infty,N_1}\tau v_B)-t_{n_i}|=|m(\dots b_{-N_1-1}^{(n_i)}\dots b_{N_1}^{(n_i)}\tau v_B)-t_{n_i}|<2^{-N_1}.
\end{equation}
Now we will use the constructions of \eqref{eq:3} and \eqref{eq:123} with some fixed $k\geq\max\{\lceil \delta^{-1} \rceil,N_1\}$. Let $\underline{\widetilde v}^{k}:=\underline{\widetilde v}_k\dots\underline{\widetilde v}_{1}$ and $v_B^k=(v_B)_1\dots (v_B)_k$ be the initial segments of $\underline{\widetilde v}^T$ and $v_B$, respectively. By transitivity of $\Sigma(B)$, there exists $\beta \in \Sigma(B)$ which contains non-overlapping occurrences of the strings  $v_B^k$ and $\underline{\widetilde v}^{k}$ in this order. Let us denote by $(v_B^k\ast\underline{\widetilde v}^{k})$ a finite substring of $\beta$ which begins with $v_B^k$ and and terminates with $\underline{\widetilde  v}^{k}$. Using all the above we obtain that
\begin{equation*}
    |m(\underline{\theta}^{(n_i)}_{-\infty,N_1}\tau(v_B^k\ast\underline{\widetilde v}^{k})\underline{\widetilde{\tau}}\widetilde{\alpha}_k\dots\widetilde{\alpha}_0\dots\widetilde{\alpha}_k \alpha_{k+1} \alpha_{k+2}\dots)-t|<2^{3-k}+\varepsilon/2<\varepsilon,
\end{equation*}
for any choice of $ \alpha_{i}\in A$ for $i\geq k+1$. This finishes the proof.

\end{proof}

It is a natural question to ask whether the local dimension of $M^\prime$ is also monotone. More precisely, define the function $\dloc^{M}:M^\prime\to[0,1]$ by
\begin{equation*}
    \dloc^{M}(t):=\lim_{\varepsilon\to 0}\dimH(M\cap(t-\varepsilon,t+\varepsilon)).
\end{equation*}
By using the explicit construction of Cantor sets in $M\setminus L$ carried out in \cite{MminusL0353}, one can see that the answer is that $\dloc^{M}$ is not monotone. Recall from~\cite{MminusL0353} the two numbers
\[b_\infty = [2; \overline{1,1,2,2,2,1,2}] + [0; \overline{1, 2,2,2, 1,1, 2}] = 3.2930442439\ldots\]
and
\begin{gather*}
B_\infty = [2; 1, \overline{1, 2,2,2, 1, 2, 1,1, 2, 1,1, 2}] + [0; 1, 2,2,2, 1,1, 2, 1, 2,2,2,\\ 1,1, 2, 1, 2,2, \overline{1, 2,2,2, 1, 2, 1,1, 2, 1,1, 2}] = 3.2930444814\ldots.
\end{gather*}
We have that $(b_\infty,B_\infty)$ is a maximal gap of $L$, hence $D(t)=\dimH(k^{-1}(-\infty,t))$ satisfies $D(b_\infty)=D(B_\infty)$, however by the explicit description of the Cantor set done in \cite{MminusL0353} one has that $0<\dimH((M\setminus L)\cap(b_\infty,B_\infty))\leq D(B_\infty)<d(B_\infty)=d(b_\infty)$ which implies that there is $s\leq b_\infty$ such that for all $t\in(M\setminus L)\cap(b_\infty,B_\infty)$ one has $\dloc^{M}(t)<d(B_\infty)$ but $\dloc^{M}(s)=d(b_\infty)$. 

From the fact that the local dimension of $M^\prime\cap (M\setminus L)$ is equal to $D(t)$ in good intervals, we can obtain an alternative proof of the fact that $\dloc^{M}$ is not monotone and obtain a stronger conclusion:
\begin{proposition}
Let $[\nu,\mu)$ be a good interval. For any  $t\in M^\prime\cap (M\setminus L)\cap[\nu,\mu)$, there is $s<t$ such that $\dloc^{M}(s)>\dloc^{M}(t)$. 
\end{proposition}

\begin{proof}
Let $z_1<t<z_2$ with $z_1,z_2\in L$ be a maximal gap containing $t$. Then we have that $D(z_1)=D(z_2)$ and $\dloc^{M}(t)=D(t)<d(z_2)=d(z_1)$, but on the other hand  by using the same proof of \Cref{prop:conjecture_characterization} one can show that there exists $s\leq z_1$ such that $\dloc^{M}(s)=d(z_1)$.
\end{proof}

A weaker question could be whether $\dloc^{M}:M^\prime\to[0,1]$ is a continuous function. Using the fact that $M\setminus L$ is not closed, established in \cite{MminusLisnotclosed}, we answer this question negatively. Recall from \cite{MminusLisnotclosed} that there are values $\widetilde{m}_k\in M\setminus L$ and $\widetilde{\ell}_k\in L$ such that $\widetilde{\ell}_k<\widetilde{m}_k<\widetilde{\ell}_{k-1}$ and $\lim_{k\to\infty}\widetilde{m}_{k}=1+\frac{3}{\sqrt{2}}\in L^\prime$.

\begin{theorem}
The function $\dloc^{M}:M^\prime\to[0,1]$ is not continuous.
\end{theorem}

\begin{proof}
Let $t=1+\frac{3}{\sqrt{2}}\simeq3.1213\dots$  and notice that it belongs to the good interval (3.0508, 3.1221) (cf. \Cref{thm:good_intervals_list} to
be proved in the next section). We claim that $\widetilde{m}_k\in M^\prime$. This implies that $\dloc^{M}(\widetilde{m}_k)=D(\widetilde{m}_k)$ while $\dloc^{M}(t)=2\cdot D(t)$. Since $\widetilde{m}_k\to t$ and $D$ is continuous, this will finish the proof. 

In fact we will exhibit a small Cantor set contained in $M\setminus L$ that contains $\widetilde{m}_k$. For $k\in\N$ let $2_{k}$ denote a string of 2s of length $k$. Recall from \cite{MminusLisnotclosed} that 
\begin{equation*}
    \widetilde{m}_k=m(\overline{w_k}w_k^{*}w_k2_{k-2}12_{2k-1}11\overline{2})=\lambda_0(\overline{w_k}w_k^{*}w_k2_{k-2}12_{2k-1}11\overline{2}), \quad \widetilde{\ell}_k=m(\overline{w_{k}}),
\end{equation*}
where $w_k=2_{k-2}12_{2k+1}12_{2k-1}12_{k+2}$ and $w^{*}_k=2_{k-2}12_{2k+1}12^*2_{2k-2}12_{k+2}$. Moreover by the proof of \cite[Proposition 1]{MminusLisnotclosed}, we know that the Markov value is only attained at 0, in fact there is $\delta_k>0$ such that $\lambda_i(\overline{w_k}w_k^{*}w_k2_{k-2}12_{2k-1}11\overline{2})\leq\widetilde{m}_k-\delta_k$ for any $i\neq 0$. Let $r_k\in\N_{>0},r_k\geq 2$ such that $2^{-r_k}<\delta_k$. In particular we will have that
\begin{equation*}
    \{\lambda_0(\overline{w_k}w_k^{*}w_k2_{k-2}12_{2k-1}112_{r_k}\gamma_1\gamma_2\dots):\gamma_i\in\{2211,11\},\forall i\geq 1\}, 
\end{equation*}
is a Cantor set contained in $(M\setminus L)\cap(\widetilde{\ell}_k,\widetilde{\ell}_{k-1})$ (we used that $2_{2k-2}12_{k+2}w_k2_{k-2}12_{2k-1}11$ is even and that $\lambda_0(22^*11)<3.06$). 
\end{proof}

\subsection{Good intervals and Berstein's investigation of $\boldsymbol{M\setminus L}$}\label{subsec:berprops}\hfill

Here we discuss how our good intervals are related to the intervals previously considered by Berstein in his work on $M\setminus L$~\cite{Berstein}.

It is also discussed in Cusick-Flahive~\cite[Chapter 3]{Cusick-Flahive}, that Berstein was trying to find necessary and sufficient conditions for a Markov value to also be contained in the Lagrange spectrum. Berstein gave a list of conditions~\cite[Theorems 2-8]{Berstein} (see~\cite[Chapter 3, Theorem 6 and Theorem 7]{Cusick-Flahive} for a summary) and constructed 23 intervals (listed in Appendix~\ref{app:ber}) on which the conditions are both necessary and sufficient for a point $m = m(\underline{a})$ to be contained in $L$.

Berstein provided a list of conditions that must hold for points in $M\setminus L$. We will prove that these conditions also hold for points $m\in(M\setminus L)\cap[\nu,\mu)$ for a good interval $[\nu,\mu)$.

\begin{definition}[$\varepsilon$-property]
A bi-infinite sequence $\underline{a}=(a_n)_{n\in\Z}$ is said to have the right (or, respectively left) $\varepsilon$-property if there is an integer $N$ such that for any different $\underline{a^\prime}=(a^\prime_n)_{n\in\Z}$ which has $a_n=a^\prime_n$ for all $n\leq N$ (respectively, $n\geq N$) we have $m(\underline{a^\prime})>m(\underline{a})+\varepsilon$.
\end{definition}

\begin{proposition}
Let $[\nu,\mu)$ be a good interval. Then for all $t=\lambda_0(\underline{a})=m(\underline{a})\in(M\setminus L)\cap[\nu,\mu)$ we have that $\underline{a}$ satisfies the right or left $\varepsilon$-property for some $N\in\N_{>0}$.
\end{proposition}

\begin{proof}
Consider $t\in[x,y)$ as in \Cref{def:good_interval}. By \Cref{rem:sup_transitive_subshifts} we have that $t\in(x,y)$. Since $m(\underline{a})\in M\setminus L$, without loss of generality we can assume that $\underline{a}$ does not connect positively to $B$ by \Cref{lem:key}. We will prove that in this case, $\underline{a}$ has the right $\varepsilon$-property. By the \Cref{def:cplus}, there is $k\in\N$ such that:
\begin{equation}\label{eq:no_positive_connection}
    \forall N\geq k+2, \forall v_B\in\Sigma^+(B), \forall\tau\in(\N_{>0})^* \quad\text{it holds}\quad m(a_{-\infty,N}\tau v_B)\geq t+2^{-k},
\end{equation}
where $(\N_{>0})^*$ is the set of all finite words in the alphabet $\N_{>0}$. In particular for any continuation $ v =  v_{N+1} v_{N+2}\dots$ with $ v_{N+1}\neq a_{N+1}$ we must have $m(a_{-\infty,N} v)\geq t+2^{-k}$, otherwise since the interval is good there will be a finite $\tau$ and $v_B\in\Sigma^+(B)$ such that $m(a_{-\infty,N}\tau v_B)\leq \max\{m(\underline{a}),m(a_{-\infty,N} v)\}<t+2^{-k}$.    
\end{proof}

Although the next proposition was already proved by Berstein (\cite[Theorems 2-8]{Berstein}), we will give here a proof for completeness.

\begin{proposition}\label{prop:Berstein_properties}
Suppose a bi-infinite sequence $\underline{a}$ satisfies the right or left $\varepsilon$-property for some $N\in\N_{>0}$ and $m(\underline{a})\in M\setminus L$. Then
\begin{enumerate}
\item $\underline{a}$ is periodic on at least one side,
\item if $p=p_1\dots p_l$ is the period of the previous item, then $\overline{p}$ has the $\varepsilon$-property to the same side as $\underline{a}$ for some $\varepsilon>m(\underline{a})-m(\overline{p})$,
\item the period $p$ is not semi-symmetric,
\item if $t=m(\underline{a})$, then $m^{-1}(t)$ is the union of finitely many orbits and any $\underline{b}\in m^{-1}(t)$ satisfies $m(\underline{b})=\lambda_i(\underline{b})$ for some $i$.
\end{enumerate}
\end{proposition}

\begin{proof}\hfill
\begin{enumerate}
\item 
This is \cite[Chapter 3, Lemma 3]{Cusick-Flahive}, which states that if $\underline{a}$ has the $\varepsilon$-property on one side, then it will be periodic on that respective side. 

\item It is not difficult to see that in fact this statement is equivalent to $\underline{a}$ having the $\varepsilon$-property.

\item Suppose $\underline{a}$ has the right $\varepsilon$-property and write $\underline{a}=\dots a_{n-1}a_n\overline{p}$ with $a_n\neq p_l$ (since $m(\underline{a})\not\in L$). If $p=rs$ where $r^T=r$ and $s^T=s$ are palindromes, then we would have that $\dots a_{n-1}a_{n}p^kra_{n}a_{n-1}\dots$ has arbitrarily close Markov value with $\underline{a}$ if $k$ is large, which contradicts the right $\varepsilon$-property for $\underline{a}$.

\item Suppose there is an infinite sequence $\underline{a}^{(i)}\in \{1,\dots,\lfloor t\rfloor\}^\Z$ with all terms distinct such that $m(\underline{a}^{(i)})=\lambda_0(\underline{a}^{(i)})=t$. Since the entries are bounded, we can assume that they converge to a limit, that is, there is a sequence $\underline{b}$ such that for any $n\in\N$ there is $i_0 = i_0(n)$ such that $a^{(i)}_{-n,n}=b_{-n,n}$ for all $i\geq i_0$. In particular we have $m(\underline{b})=\lambda_0(\underline{b})=t$. By the first item, we can assume that $\underline{b}$ has the right $\varepsilon$-property for some $N\in\N_{>0}$. Fix some $j\geq i_0=i_0(N_1)$ where $N_1\geq N$ is such that $2^{-N_1}<\varepsilon$. Since $\underline{a}^{(i)}\neq \underline{b}$ for all $i$, we must have $a^{(j)}_{N_1+1,\infty}\neq b_{N_1+1,\infty}$, otherwise it will contradict the fact that $\lambda_0(\underline{a}^{(j)})=\lambda_0(\underline{b})$ but $\underline{a}^{(j)}\neq\underline{b}$. Hence, using that $a^{(j)}_{-N_1,N_1}=b_{-N_1,N_1}$, we have for $i\geq 0$ 
\begin{equation*}
    \lambda_i(b_{-\infty,N_1}a^{(j)}_{N_1+1,\infty})\leq m(\underline{a}^{(j)})+2^{-N_1} < t+\varepsilon,
\end{equation*} 
while for $i<0$
\begin{equation*}
    \lambda_i(b_{-\infty,N_1}a^{(j)}_{N_1+1,\infty})\leq m(\underline{b})+2^{-N_1} < t+\varepsilon,
\end{equation*} 
which contradicts the fact that $\underline{b}$ satisfies the right $\varepsilon$-property. 

Hence we have finitely many bi-infinite sequences $\underline{a}^{(1)},\dots,\underline{a}^{(r)}$ such that $m(\underline{a}^{(i)})=\lambda_0(\underline{a}^{(i)})=t$. Let $\underline{c}\in \{1,\dots,\lfloor t\rfloor\}^\Z$ be any other bi-infinite sequence such that $m(\underline{c})=t$ but $\lambda_i(\underline{c})<t$ for all $i\in\Z$. By definition of Markov value there is a subsequence say $n_1<n_2<\dots$ such that $\lambda_{n_i}(\underline{c})\to t$. By compactness this means that, restricting to a further subsequence if necessary , that $\sigma^{n_i}(\underline{c})$ converges to a sequence which has Markov value $t$ and that attains it precisely at the index $0$. Since there are only finitely such bi-infinite sequences say it converges to to $\underline{a}^{(1)}$. However, since $\underline{a}^{(1)}$ has the $\varepsilon$-property to one side, the same proof above yields that $\sigma^{n_i}(\underline{c})=\underline{a}^{(1)}$ for some $i$, contrary to the hypothesis.    
\end{enumerate}

\end{proof}

\begin{remark}
As mentioned in the introduction, it is important to note that Berstein's intervals are constructed in a similar way to our good intervals. That is, Berstein's proofs of the necessity and sufficiency of the conditions for $m(\underline{a})$ being in $L$ depend on establishing that one theorem~\cite[Theorem 2]{Berstein} is true for each such interval. In the proof of this theorem, Berstein considers two distinct continuations of a particular sequence and constructs a third continuation giving rise to a Markov value that is no greater than the values for the original two continuations. We see that this is very reminiscent of our definition of a good interval. However, and importantly so, Berstein is not considering the transitivity of two subshifts related to the interval. These subshifts are crucial to our proofs on local dimension, etc., given above and below.
\end{remark}

\subsection{Properties of $D(t)$ on good intervals}

We finish this section by showing that the main results of  \cite{MV} also hold on good intervals. Recall that $D(t) = \dimH(K_{t}) = \dimH(k^{-1}(-\infty,t])$.

\begin{proposition}\label{prop:concentration}
Let $[\nu,\mu)\subset\R$ be an interval. If $[\nu,\mu)$ is good, then for all $t\in L^\prime\cap[\nu,\mu)$ we have
\begin{equation}\label{eq:concentration_of_dim}
    \dimH (k^{-1}(t))=D(t).
\end{equation}

\end{proposition}

\begin{proof}
Consider $t\in[x,y)$ as in \Cref{def:good_interval}. Fix $\eta\in(0,1)$ so that there is $\delta>0$ and a complete shift $\Sigma(A)$ such that $\Sigma(A)\subset\Sigma_{t-\delta}$ and $\dimH(K(A))>(1-\eta)\dimH(K_t)$. As before, since $t\in L^\prime$, for each $k\in\N^{*}$, there are regular Cantor sets $K(A_1^{(k)}),K(A_2^{(k)})$ defined by iterates of the Gauss map such that $$\diam\left(t-\big(K(A_1^{(k)})+K(A_2^{(k)})\big)\right)\leq2^{3-k},$$ and $K(A_1^{(k)})+K(A_2^{(k)})\subset L$. Fix an element $m(\gamma^{(k)})=\lambda_0(\gamma^{(k)})\in K(A_1^{(k)})+K(A_2^{(k)})$ or choose $\gamma^{(k)}\in\Sigma(B)$ with $\lambda_0(\gamma^{(k)})=m(\gamma^{(k)})=t$ if $t=x$. By repeating the same construction done in \eqref{eq:all_connected}, there is a finite word

\begin{multline*}
    \hat{\theta}^{(k)}=\widetilde{\alpha}_{-k}^{(k)}\dots\widetilde{\alpha}_0^{(k)}\dots\widetilde{\alpha}_k^{(k)}\underline{\tau}^{(n)}(\underline{ v}^{l+k}\ast\widetilde{ v}^{l+k})^{(k)}\widetilde{\tau}^{(k)}\gamma_{-k}^{(k)}\dots\gamma_0^{(k)}\dots \\
    \dots\gamma_{k}^{(k)}\tau^{(k)}( v^{l+k}\ast\underline{\widetilde{ v}}^{l+k})^{(k)}\underline{\widetilde{\tau}}\widetilde{\alpha}_k^{(k)}\dots\widetilde{\alpha}_0^{(k)}\dots\widetilde{\alpha}_k^{(k)},
\end{multline*}
such that for any choice $ \alpha_{|i|}\in A$ for $|i|\geq 1$ one has that
\begin{equation}\label{eq:final}
    |m(\dots  \alpha_{-2} \alpha_{-1}\hat{\theta}^{(k)} \alpha_1 \alpha_2 \dots)-t|<2^{4-k}.
\end{equation}
Denote by $r_k$ the length (in the alphabet $\{1,\dots,\lfloor t\rfloor\}$) of the finite word $\hat{\theta}^{(k)}$ and let $s_k=r_1+\dots+r_k$. Given $(\alpha_1,\alpha_2,\dots)\in \Sigma^{+}(A)$ where $\alpha_i\in A$, denote $\alpha(k)=\alpha_{s_k!+1}\dots \alpha_{s_{k+1}!}$. For each $z=[0;\alpha_1,\alpha_2,\dots]\in K(A)$, define the map $h:K(A)\to k^{-1}(t)$ by
\begin{equation*}
    h(z)=[0;\alpha_1,\dots,\alpha_{s_1!},\hat{\theta}^{(1)},\alpha(1),\hat{\theta}^{(2)},\alpha(2),\hat{\theta}^{(3)},\dots].
\end{equation*}
The fact that $k(h(z))=t$ comes from \eqref{eq:final}. This map is clearly injective and continuous. Moreover, for any small $\rho>0$, one has that $|z-z^\prime|=O(|h(z)-h(z^{\prime})|^{1-\rho})$ uniformly for all $z,z^\prime\in K(A)$. To see this, first note that by the bounded distortion property, given any finite word $a_1\dotsb a_n$ in the alphabet $\{1,\dots,T\}$ where $T=\lfloor t\rfloor$, there are positive constants $C_T$ and $\lambda_1(T)<\lambda_2(T)<1$ all depending only on $T$ such that one has
\begin{equation*}
    C_T^{-1}\lambda_1^{n}(T)<\sizes(a_1\dots a_n)<C_T\lambda_2^n(T),
\end{equation*}
where 
\begin{equation*}
    \sizes(a_1\dots a_n)=\left| \frac{p_n}{q_n}-\frac{p_n+p_{n-1}}{q_n+q_{n-1}} \right| = \frac1{q_n(q_n+q_{n-1})},
\end{equation*}
and $[0;a_1,a_2,\dots,a_n]=p_n/q_n$. Another useful property (\cite[Lemma A.2]{geometricproperties}) is the fact that given any finite words $\alpha$ and $\beta$ one has that $\frac{1}{2}\sizes(\alpha)\sizes(\beta)<\sizes(\alpha\beta)<2\sizes(\alpha)\sizes(\beta)$. Moreover, it follows from \cite[Lemma A.1]{geometricproperties} that if $x=[0;a_1,\dots,a_n,a_{n+1},\dots]$, $x^\prime=[0;a_1,\dots,a_n,b_{n+1},\dots]$ differ precisely at the $n+1$ term ($a_{n+1}\neq b_{n+1}$) and the coefficients are bounded by $T$, then $|x-x^\prime|\geq \widetilde{C}_T\sizes(a_1\dots a_n)$ for some positive constant $\widetilde{C}_T$ that depends only on $T$. Finally, since the denominators $q_k$ are increasing, it is clear that  $\sizes(a_1\dots a_n)>\sizes(a_1\dots a_m)$ for $m>n$.

Now, given $z,z\in K(A)$, let $\alpha_{s+1}$ be the first letter of $A$ in the continued fractions of $z$ and $z^\prime$ where $z$ and $z^\prime$ differ and consider $k$ maximal such that $s+1 \leq s_{k+1}!$. Note that the continued fractions of $z$ and $z^\prime$ will differ in some digit inside the letter $\alpha_{s+1}$. In particular, since $z,z^\prime \in I(\alpha_1\dots \alpha_s)$, we have that $|z-z^\prime|\leq \sizes(\alpha_1\dots \alpha_s)$. On the other hand 
\begin{align*}
    |h(z)-h(z^\prime)|^{1-\rho}&\geq \widetilde{C}_T^{1-\rho}\sizes(\alpha_1\dotsb \alpha_{s_1!}\hat{\theta}^{(1)}\alpha(1)\dotsb  \alpha_{s+1})^{1-\rho} \\
    &=\widetilde{C}_T^{1-\rho}\sizes(\alpha_1\dotsb \alpha_{s_1!}\hat{\theta}^{(1)}\alpha(1)\dotsb  \alpha_{s+1})\sizes(\alpha_1\dotsb \alpha_{s_1!}\hat{\theta}^{(1)}\alpha(1)\dots  \alpha_{s+1})^{-\rho} \\
    &\geq 2^{-3(k+1)}\widetilde{C}_T^{1-\rho}\sizes(\alpha_1\dotsb \alpha_{s})\sizes(\hat{\theta}^{(1)})\dotsb\sizes(\hat{\theta}^{(k)})\sizes(\alpha_{s+1})\sizes(\alpha_1\dotsb \alpha_{s_1!}\hat{\theta}^{(1)}\alpha(1)\dotsb  \alpha_{s+1})^{-\rho} \\
    &\geq 2^{-3(k+1)}\widetilde{C}_T^{1-\rho}\sizes(\alpha_1\dotsb \alpha_{s})C_T^{-(k+1)}\lambda_1(T)^{s_k+|\alpha_{s+1}|}C_T^{-\rho}\lambda_2(T)^{-\rho\cdot s_k!} \\
    &\geq \lambda_1(T)^{N_A}C_{T,\rho}\sizes(\alpha_1\dotsb \alpha_{s})
\end{align*}
for some constant $C_{T,\rho}$ that only depends on $T=\lfloor t\rfloor$ and $\rho$ and where $|\alpha_{s+1}|$ is the length of $\alpha_{s+1}$ in the alphabet $\{1,\dots,T\}$ and $N_A=\max_{\alpha\in A}|\alpha|$ (we also used that $e^{c_1\cdot n!-c_2\cdot n}$ is bounded by below for all $n\in\N$ by some constant that depends on the parameters $c_1,c_2>0$).

Therefore using this map $h$, for fixed $\eta\in(0,1)$ (and thus fixed $A$), one obtains that for any small $\rho>0$
\begin{align*}
    \dimH(k^{-1}(t))\geq(1-\rho)^{-1}\dimH(K(A))\geq(1-\rho)^{-1}(1-\eta)D(t).
\end{align*}
Letting $\rho\to 0$ and then $\eta\to 0$ shows that $\dimH(k^{-1}(t))\geq D(t)$. The other inequality follows from the fact that $D(t)=\dimH(k^{-1}(-\infty,t])\geq \dimH(k^{-1}(t))$.
\end{proof}

\begin{remark}
The reason why we restrict $t\in L^\prime$ in Proposition~\ref{prop:concentration} is because in general, when $t\in L$ is an isolated point in $L$, the preimage $k^{-1}(t)$ is countable: this follows from the same proof of \cite[Theorem 3]{geometricproperties} (in fact the same proof shows that there are only finitely many periodic orbits in $\ell^{-1}(t)$). In particular when $t\in L$ is isolated in $L$ and $t>3$, the equation \eqref{eq:concentration_of_dim} fails because $\dim_H(k^{-1}(t))=0$ but $D(t)>0$.
\end{remark}

We will use a slight variant of \cite[Lemma 2.5]{Phasetransition}. The following lemma states that if we remove one of the parts $J$ of a Markov partition $\mathcal{P}$ of a regular Cantor set, then the maximal invariant set of $\mathcal{P}\setminus\{J\}$ has Hausdorff dimension strictly smaller than $K$.

\begin{lemma}\label{lem:proper_cantor_set}
Let $(K,\mathcal{P},\psi)$ be a $C^k$-regular cantor set, $k>1$ and $J\in\mathcal{P}$. Then the maximal invariant set
\begin{equation*}
    C:=\bigcap_{n\geq 0}\psi^{-n}\left(\bigcup_{I\in\mathcal{P},I\neq J}I\right)
\end{equation*}
satisfies $\dimH(C)<\dimH(K)$.
\end{lemma}

\begin{proposition}
Let $[\nu,\mu)\subset\R$ be a good interval. If $s_1<t<s_2$ where $t\in L^\prime\cap[\nu,\mu)$, then
\begin{equation}\label{eq:Disstrictlyincreasing}
    D(s_1)<D(s_2).
\end{equation}

\end{proposition}

\begin{proof}
Consider $t\in[x,y)$ as in \Cref{def:good_interval}. First, we will show that for any $t\in[x,y)$ and $\varepsilon>0$, any sequence $\gamma\in\Sigma_{t-\varepsilon}$ is nonwandering in $\Sigma_{t}$. Assume $n\in\N$ is so big so that $t-\varepsilon+2^{-n}<\min\{y,t-\varepsilon/2\}$. Consider the finite set of words
\begin{equation*}
    F(n,t-\varepsilon)=\{a_1\dots a_{2n+1}\in\{1,\dots,\lfloor t\rfloor\}^{2n+1}:\lambda_0(a_1\dots a_{n}a_{n+1}^*a_{n+2}\dots a_{2n+1})>t-\varepsilon\}. 
\end{equation*}    

Clearly, one has that $\Sigma_{t-\varepsilon}\subset\Sigma(F)\subset\Sigma_{t-\varepsilon+2^{-n}}$ where $\Sigma(F)$ is the subshift of finite type that forbids all the words in $F(n,t-\varepsilon)$. As explained in \cite[Chapters 1 and 5]{Kitchens}, to a subshift of finite type we can associate a transition matrix $M$ that encodes the subshift. Let $\Sigma_{M_1},\dots,\Sigma_{M_r}\subset\Sigma_M$ be the irreducible components of $\Sigma_M=\Sigma(F)$.

Since each $\Sigma_{M_i}$ is a transitive subshift and the interval $[x,y)$ is good, we have that any $\gamma\in\Sigma_{M_i}$ connects positively and negatively with $B$ before $t-\varepsilon/2$, unless $\Sigma_{M_i}$ is trivial, in which case $\gamma$ is a periodic orbit, which is clearly nonwandering. By \cite[Observation 5.1.2]{Kitchens}, we know that for any $\gamma\in\Sigma(F)$ its positive and negative limit sets are each contained in an irreducible component of $\Sigma_M$. In particular, this implies that for any sufficiently large $N$, there is a finite word $\tau$ and a continuation $\upsilon\in\Sigma^+(B)$ such that $m(\gamma_{-\infty,N}\tau\upsilon)=m(\dots\gamma_{N-1}\gamma_{N}\tau\upsilon)< t-\varepsilon/4$. Similarly, there is a finite sequence $\tilde{\tau}$ and a continuation $\tilde{\upsilon}\in\Sigma^+(B)$ such that
\begin{equation}\label{eq:nonwanderingconnection}
    m(\tilde{\upsilon}^T\tilde{\tau}^T\gamma_{-N}\dots\gamma_N\tau\upsilon)<t.
\end{equation} 
Since $\Sigma(B)$ is transitive this implies that $\gamma$ is nonwandering in $\Sigma_{t}$. 

Now we consider the subshift of finite type $\Sigma(\tilde{F})$ that forbids all words in $F(n,t)$. In particular $\Sigma_{t}\subset\Sigma(\tilde{F})\subset\Sigma_{t+\varepsilon}$. By repeating the above construction, there is a transition matrix $\tilde{M}$ for this subshift $\Sigma(\tilde{F})=\Sigma_{\tilde{M}}$ and there are irreducible components $\Sigma_{\tilde{M}_1},\dots,\Sigma_{\tilde{M}_{\tilde{r}}}\subset\Sigma_{\tilde{M}}$. We claim that all nontrivial components $\Sigma_{M_j}$ are contained in a single irreducible component of $\Sigma_{\tilde{M}}$, say $\Sigma_{\tilde{M}_i}$. First note that since $\Sigma(B)$ is transitive, it must be contained in a single irreducible component, say $\Sigma_{\tilde{M}_i}$. Now, given any $\gamma\in\Sigma_{M_j}$ with $\Sigma_{M_j}$ nontrivial, by \eqref{eq:nonwanderingconnection} we see that $\gamma$ must be in the same transitive component as $\Sigma(B)$. In conclusion we have that $\bigcup_{j}\Sigma_{M_j}\subset\Sigma_{\tilde{M}_i}\subset\Sigma_{t+\varepsilon}$ where the union is over nontrivial components $\Sigma_{M_j}$. 

Now let $\gamma\in\Sigma_t$ be such that $m(\gamma)=\lambda_0(\gamma)=t$. Since $t\in L^\prime\cap[x,y)$, we can assume that $\gamma$ connects positively and negatively to $\Sigma(B)$, so in particular belongs to the same irreducible component as $\Sigma(B)$, so it must belong to $\Sigma_{\Tilde{M}_i}$. On the other hand, for some $k$ large we will have that any sequence containing $\gamma_{-k}\dots\gamma_k$ will have Markov value at least $t-2^{-k+1}>t-\varepsilon/2$. Hence, if we forbid this finite word from the subshift of finite type $\Sigma_{\tilde{M}_i}$, we will get a new subshift of finite type $\Sigma(\hat{F})$ that contains $\bigcup_{j}\Sigma_{M_j}$. By \Cref{lem:proper_cantor_set} we must have that $\dimH(K(\Sigma(\hat{F})))<\dimH(K(\Sigma_{\tilde{M}_i}))$.

Finally,
\begin{align*}
    D(t-\varepsilon)=\dimH(K_{t-\varepsilon})&\leq\sup_{j}\dimH(K(\Sigma_{\hat{M}_j})) \\
    &\leq\dimH(K(\Sigma(\hat{F}))<\dimH(K(\Sigma_{\Tilde{M}_i}))\leq \dimH(K_{t+\varepsilon})=D(t+\varepsilon).
\end{align*}

Now take $\varepsilon>0$ small enough so that $s_1<t-\varepsilon<t+\varepsilon<s_2$ and use the monotonicity of $D$.
\end{proof}

\section{Construction of good intervals}\label{sec:good_ints}

In this section, we will prove Theorem~\ref{thm:dloc} by demonstrating that the following is a list of good intervals (in the sense of \Cref{def:good_interval}). One important feature of the intervals found so far is the following: they contain almost all of the known elements of $M\setminus L$. The only exceptions are the finite set of elements $m(\gamma_k^1)$, $k\in\{1,2,3,4\}$ found in \cite{MminusLisnotclosed}. The only reason that these elements do not feature, is that they lie very close to 3 which is a region of the spectra that we have not yet investigated.

\begin{theorem}\label{thm:good_intervals_list}
The following are good intervals (in the sense of \Cref{def:good_interval}).
\begin{itemize}
\item $[3.05082, 3.122183)$
\item $[3.1299, 3.285441)$
\item $[3.28603, 3.28729)$
\item $[3.29296, 3.29335)$
\item $[3.33396, 3.33475)$
\item $[3.359, 3.423)$
\item $[\sqrt{12},3.8465)$
\item $[3.873,3.930691)$
\item $[3.93616, 3.943767)$
\item $[3.944054, 3.971606)$
\item $[3.97995, 3.9857)$
\item $[4.520781, 4.523103)$
\item $[4.5251, 4.5279)$.
\end{itemize}
\end{theorem}

\begin{remark}
    The precise endpoints of the good intervals we construct correspond to certain extremal Markov values, whose explicit symbolic representations can be found on each respective subsection. These values are explicit sums of quadratic irrationals, so we decided to just give the first decimal digits of their representation.
\end{remark}

Our method of proof will be to construct two subshifts $\Sigma(B)$ and $\Sigma(C)$ and calculate (or bound) the maximum and minimum Markov values for sequences in $\Sigma(B)$ and $\Sigma(C)$, respectively. Then, given two continuations $\ldots a_{-1}a_{0}a_{1}\ldots a_{N}v^{i}$ of a sequence $\underline{a} = (a_{i})_{i\in\Z}$ satisfying the hypotheses of the definition of a good interval, we will find a continuation $\ldots a_{-1}a_{0}a_{1}\ldots a_{N}v_{B}$ with $v_{B}\in\Sigma^{+}(B)$ and
\[[0;v^{1}]>[0;v_{B}]>[0;v^{2}].\]
We will then make use of Lemma~\ref{lem:tech} and the argument discussed in Remark~\ref{rem:MtoL} to bound the Markov value of the continuation $\ldots a_{-1}a_{0}a_{1}\ldots a_{N}v_{B}$ and, in doing so, prove that the continuation $v_{B}$ works as intended. Finally, one can check explicitly that the subshifts $\Sigma(B)$ and $\Sigma(C)$ constructed are transitive through \Cref{prop:mixing_criteria} by connecting all possible orbits with the periodic orbit $\overline{1}$, which for the intervals constructed here, will always belong to both subshifts. This last part (verifying transitivity) can be verified computationally and is omitted from the paper.

We will argue the first few intervals in full detail before dropping some details once the method has become clear to the reader.

\bigskip

\noindent{\Large\textbf{Intervals before $\boldsymbol{t_1=3.334384\dots}$}}
\vspace*{5pt}

Note here that, for $t\in L'\cap[x,y)$ for a good interval $[x,y)$, we will have $\dloc(t)=d(t) < 1$. Moreover, all sequences here will be from $\{1,2\}^{\Z}$.

\subsection{Interval [3.05082, 3.122183)}\hfill

Let $\Sigma(C)=\{\underline{a}\in\{1,2\}^\Z\colon\text{121, 122212 and 212221 are not substrings of}~\underline{a}\}$. The minimum word containing 121 is $m(\overline{12^*1})=3.162\dots$ (indeed, if 1212 is forbidden then a continued fraction in $\{1,2\}$ that begins with $[2;1]$ is minimized when $[2;\overline{1,1,2}]$). Note that
\begin{equation*}
    \lambda_0(1222^*122)\geq 3.1216
\end{equation*}
which is currently below our claimed right endpoint of the interval. We claim that the right endpoint of the interval can be pushed to 
\begin{equation*}
    m(\overline{1222^*12212221})=3.122183\dots
\end{equation*}
which corresponds to the minimum Markov value of a word containing 1222122.

\begin{remark}
The minimum is a palindrome. 
\end{remark}

Indeed, since $\lambda_0(1222^*1222)>3.1257$, $\lambda_0(1222^*12211)>3.1228$, $\lambda_0(212^*2212)>3.1248$,
\linebreak \hbox{$\lambda_0(111222^*12)>3.1225$} are greater than the above candidate and since 121 is forbidden, the subword 1222122 must extend to 2211222122122. Since $\lambda_0(2211222^*1221221)>3.1222$ and \linebreak $\lambda_0(2211222^*12212222)>3.122187$, this subword must extend to 221122212212221 and using the previous forbidden words it extends to 221122212212221122. Hence it has the form $w_{1}^T\beta^t2^*\theta2\beta w_{2}$ where $\theta=1221$ is an even palindrome, $\beta=221122$ and $w_{1},w_{2}\in\{1,2\}^\N$. Hence, by \Cref{lem:minimizing}, this is minimized when $w_{1}=w_{2}$ and $[0;\beta,w_{1}]$ is minimal. Using the fact that 121, 12221222, 122212211, 22112221221221, 221122212212222 are forbidden, a continued fraction with coefficients in $\{1,2\}$ that begins with $[0;2,2,1,1,2,2]$ is minimized when $[0;\overline{2,2,1,1,2,2,2,1,2,2,1,2}]$.

\begin{remark}\label{rem:forbid121}
In order to explain the choice of 122212 as a forbidden word, observe that the word $\overline{1222}$ is the word with the maximum Markov value over words without 121  ($m(\overline{12^*22})>3.1298$). So, when building $\Sigma(C)$, we need to avoid this word somehow. We do this by forbidding arbitrarily large subsequences of $\overline{1222}$. We begin with the word 122212 because it is the first one with dangerous positions. Since $121$ is forbidden, it extends to 1222122.
\end{remark}

Let $\Sigma(B)=\{\underline{a}\in\{1,2\}^\Z\colon\text{121, 212 are not substrings of}~\underline{a}\}$. We have that if $\underline{b}\in\Sigma(B)$ then $m(\underline{b})\leq m(\overline{1112^*22})=3.050816\ldots$ (in fact according to Cusick-Flahive~\cite[Table 1, Chapter 5]{Cusick-Flahive}, we have a gap $(m(\overline{1112^*22}),m(\overline{12^*2}))=(3.050816\ldots, 3.073\ldots)$).

So we set $[x,y) = [m(\overline{1112^*22}),m(\overline{1222^*12212221}))=[3.050816\ldots,3.122183\ldots)$.

We will prove that the interval $[x,y)$ is good by using $\Sigma(B)$ and $\Sigma(C)$ for every $m := m(\underline{a})<y$. That is, we will not require local intervals of transitivity. Note that $\Sigma(C)$ is transitive because it contains the periodic orbit $\overline{1}$ and so by \Cref{prop:mixing_criteria} it rests to verify that the connections $122\overline{1}$, $1222\overline{1}$, $212222\overline{1}$ do not contain any of the forbidden strings defining $\Sigma(C)$. Similarly for $\Sigma(B)$ it rests to verify that $122\overline{1}$ and $2\overline{1}$ do not contain 121 nor 212.

Let $m := m(\underline{a})=\lambda_{0}(\underline{a})<y$, then $\underline{a}\in\Sigma(C)$. Suppose that, as in the definition of a good interval, we have two distinct continuations $\ldots a_{-1}a_{0}a_{1}\ldots a_{N}v^{i}$ of $\underline{a}$ with $m(\ldots a_{-1}a_{0}a_{1}\ldots a_{N}\upsilon^{i}) < y$. Let $v^{i} = v^{i}_{1}v^{i}_{2}\ldots$. We may assume that $v^{1}_{1}\neq v^{2}_{1}$. Since $\Sigma(C)\subset\{1,2\}^{\Z}$, we may further assume that $v^{1}_{1} = 1$ and $v^{2}_{1} = 2$. We observe that, since 121 and 122212 are both forbidden,
\[[0;v^{1}] = [0;1,v^{1}_{2},v^{1}_{3},\ldots] \geq [0;1,1,\overline{2,2,2,1,1,1}] =:[0;v_{B}]\in K(B).\]
So either $v^{1} = v_{B} \in\Sigma^{+}(B)$ and we had nothing to prove, or since $v^{2}_{1} = 2$ we have
\[[0;v^{1}]>[0;v_{B}]>[0;v^{2}].\]
Hence, by Lemma~\ref{lem:tech}, we have
\[\lambda_{j}(\ldots a_{-1}a_{0}a_{1}\ldots a_{N}v_{B}) \leq \max_{i=1,2} m(\ldots a_{-1}a_{0}a_{1}\ldots a_{N}v^{i}),\,\,\,\,\,\,j\leq N+1,\]
and since for positions $j>N+1$ the value is small because $\lambda_{0}(112^{*}2)\leq x$, it follows that
\[m(\dots a_{-1}\alpha_0\dots a_Nv_B)\leq \max\{m(\dots a_{-1}a_0\dots a_Nv^1),m(\dots a_{-1}a_0\dots a_Nv^2),x\},\]
as required.

Hence the interval $[x,y) = [m(\overline{1112^*22}),m(\overline{1222^*12212221}))=[3.050816\ldots,3.122183\dots)$ is good.

\subsection{Interval [3.1299, 3.2854441)}\hfill

We will prove that this interval is good by showing that the intervals [3.1299, 3.2811) and [3.2659, 3.2854441) are both good. 
We will again give full details.

\subsubsection{Interval $[3.12984\dots,3.2811\dots)$}\hfill

Let $\Sigma(C)=\{\underline{a}\in\{1,2\}^\Z\colon\text{1212 and 2121 are not substrings of}~\underline{a}\}$. The minimum Markov value of a bi-infinite sequence containing 2121 (or its transpose) corresponds to $m(\overline{1212212^*1})=3.2811\dots$. In fact, this point is an isolated point between $m(\overline{2^*111})=3.2659\dots$ and \linebreak $m(\overline{2122}12112^*1\overline{2212})=3.2812\dots$. This was proved by Matheus-Moreira-Vytnova~\cite[Subsection 3.1.5]{MMV}. Hence we take
$y=m(\overline{1212212^*1})=3.2811\dots$.

Let $\Sigma(B)=\{\underline{a}\in\{1,2\}^\Z\colon\text{121 is not a substring of}~\underline{a}\}$. If $\underline{b}\in\Sigma(B)$ then, as discussed in Remark~\ref{rem:forbid121}, $m(\underline{b})\leq m(\overline{12^*22})=3.12984\ldots$.

So we set $[x,y)=[m(\overline{12^*22}),m(\overline{1212212^*1}))=[3.12984\dots,3.2811\dots)$.

Let $m := m(\underline{a})=\lambda_{0}(\underline{a})<y$, then $\underline{a}\in\Sigma(C)$. Suppose that, as in the definition of a good interval, we have two distinct continuations $\ldots a_{-1}a_{0}a_{1}\ldots a_{N}v^{i}$ of $\underline{a}$ with $m(\ldots a_{-1}a_{0}a_{1}\ldots a_{N}v^{i}) < y$. Let $v^{i} = v^{i}_{1}v^{i}_{2}\ldots$. We may assume that $v^{1}_{1}\neq v^{2}_{1}$. Since $\Sigma(C)\subset\{1,2\}^{\Z}$, we may further assume that $v^{1}_{1} = 1$ and $v^{2}_{1} = 2$. We observe that, since 1212 and 2121 are both forbidden,
\[[0;v^{2}] = [0;2,v^{2}_{2},v^{2}_{3},\ldots] \leq [0;2,2,\overline{1,2,2,2}] =:[0;v_{B}]\in K(B).\]
So either $v^{2} = v_{B} \in\Sigma^{+}(B)$ and we had nothing to prove, or since $v^{1}_{1} = 1$ we have
\[[0;v^{1}]>[0;v_{B}]>[0;v^{2}].\]
Hence, by Lemma~\ref{lem:tech}, we have
\[\lambda_{j}(\ldots a_{-1}a_{0}a_{1}\ldots a_{N}v_{B}) < \max_{i=1,2} m(\ldots a_{-1}a_{0}a_{1}\ldots a_{N}v^{i}),\,\,\,\,\,\,j\leq N+1,\]
so, by the argument in Remark~\ref{rem:MtoL}, it follows from that fact that $\lambda_{0}(22^{*}\overline{1222}),\lambda_{0}(2212^{*}22\overline{1222}) \leq m(\overline{12^*22}) = x$ that
\[m(\ldots a_{-1}a_{0}a_{1}\ldots a_{N}v_{B}) < \max\{m(\dots a_{-1}a_0\dots a_Nv^1),m(\dots a_{-1}a_0\dots a_Nv^2),x\},\]
as required.

Hence the interval $[x,y)=[m(\overline{12^*22}),m(\overline{1212212^*1}))=[3.12984\dots,3.2811\dots)$ is good.

\subsubsection{Interval $[3.2659\dots,3.28544419\dots)$}\hfill

Let $F_{C_{0}} = \{21212,212111,12121,2121122,12221211211\}$.

Now consider 
\begin{equation*}
    \Sigma(C_0)=\{\underline{a}\in\{1,2\}^{\Z}:\forall w\in F_{C_{0}},\,w\text{ and its transpose are not substrings of}~\underline{a}\}.
\end{equation*}

Observe that
\begin{align*}
    \lambda_0(12^*121)&>3.297 \\
    \lambda_0(212^*12)&>3.4 \\
    \lambda_0(212^*111)&>3.314 \\
    \lambda_0(212^*1122)&>3.2884 \\
    \lambda_0(122212^*11211)&>3.2872
\end{align*}

The right extreme of this region is determined by $m(\overline{1211}2^*1222\overline{1121})=3.2872978\dots$, which corresponds to the minimum Markov value of a string containing $12221211211$. Indeed, one only has to use that 212111 and 2121122 are forbidden.

Let
\begin{equation*}
    \Sigma(B_0)=\{\underline{a}\in\{1,2\}^{\Z}:1212,2121\,\text{are not substrings of}~\underline{a}\}
\end{equation*}
Thus we have that $m(\underline{b})\leq m(\overline{2^*111})=3.2659\dots$ for all $\underline{b}\in\Sigma(B_0)$. Indeed, one just need to use that if 2121, 1212 are forbidden then a continued fraction in $\{1,2\}$ that begins with $[2;1,1]$ is maximized with $[2;\overline{1,1,1,2}]$ and the fact that $\lambda_0(12^*2)<3.16$ and $\lambda_0(22^*2)<3$.

So $\Sigma(B_{0})$ and $\Sigma(C_{0})$ can be used to analyse the interval [3.2659\ldots, 3.28729\ldots). We will start our analysis here before specialising since we will also use this preliminary analysis for the interval in the next subsection.

So let $m = m(\underline{a}) = \lambda_{0}(\underline{a})<3.28729\ldots$ and suppose that we have two continuations $v^{1} = 1v^{1}_{2}\ldots$ and $v^{2} = 2v^{2}_{2}\ldots$ satisfying the hypothesis of the definition of a good interval. If $a_{N} = 1$, then since 111212 is forbidden
\[[0;v^{1}] = [0;1,v^{1}_{2},\ldots] \geq [0;1,\overline{1,2,1,1}]\in K(B_{0}).\]
So either $v^{1}\in\Sigma^{+}(B_{0})$ or, since $\lambda_{0}(112^{*}\overline{1112})\leq m(\overline{2^*111})=3.2659\ldots$, we can use the continuation $1\overline{1211}\in\Sigma^{+}(B_{0})$. So from here on we only need to consider the situation with $a_{N} = 2$.

If $a_{N} = 2$ and $a_{N-1} = 2$, then since 2211212 is forbidden we can again show that
\[[0;v^{1}] \geq [0;1,\overline{1,2,1,1}]\in K(B_{0}),\]
so that the continuation $1\overline{1211}$ can be used.

So, going forward, we need only consider the case $a_{N-1}a_{N} = 12$. At this point, however, if one tries to lower bound $[0;v^{1}]$ or upper bound $[0;v^{2}]$ then one is led to continuations that are not contained in $\Sigma^{+}(B_{0})$, and so similar arguments to the above do not follow immediately. Instead, we must further restrict the Cantor sets.

Now consider the $\Sigma(C) = \{\underline{a}\in\Sigma(C_{0}):\text{1222121 and 1212221 is not a substring of}~\underline{a}\}$. Note that if a sequence in $\{1,2\}^{\Z}$ contains the subword 1222121, then this subword must extend to 1222121121 since 21212, 212111 and 2121122 are forbidden. Now observe that
\begin{equation*}
    \lambda_0(122212^*112122)>3.2854.
\end{equation*}
In fact, the minimum Markov value of a sequence containing 122212112122 is
\begin{equation*}
    m(\overline{1211}22212112^*1222\overline{1121})=3.28544419\dots
\end{equation*}
Indeed, since $\lambda_0(122212^*1121221)>3.2855$ and $\lambda_0(122212^*11212222)>3.285447$ are greater than the above candidate, if 122212112122 is a subword then it must extend as 12221211212221. Using \Cref{lem:minimizing} and the fact that 2121122 and 212111 are forbidden, one confirms that the above is the minimum.

We set $\Sigma(B) = \Sigma(B_{0})$.

So, we set $[x,y) = [m(\overline{2^*111}),m(\overline{1211}22212112^*1222\overline{1121}))=[3.2659\dots,3.28544419\dots)$.

Let $m = m(\underline{a}) = \lambda_{0}(\underline{a}) <y$. By the above discussion, we need only consider the continuations $v^{1} = 1\ldots$ and $v^{2} = 2\ldots$ with $a_{N-1}a_{N} = 12$. In such a case, since 1222121 is forbidden we have
\[[0;v^{2}] = [0;2,v^{2}_{2},\ldots] \leq [0;\overline{2,2,1,2}]\in K(B).\]
So, since $\lambda_{0}(222^{*}122) < 3.14$, we can use the continuation $\overline{2212}\in\Sigma^{+}(B)$.

Hence, the interval $[3.2659\dots,3.28544419\dots)$ is good. Therefore, the larger interval $$[m(\overline{12^*22}),m(\overline{1211}22212112^*1222\overline{1121})) = [3.12984\dots,3.28544419\dots)$$ is good.

\subsection{Interval [3.28603, 3.28729)}\hfill

Consider again the Cantor set $\Sigma(C_{0})$ from the previous subsection.

We will set 
\begin{equation*}
    \Sigma(C)=\{\underline{a}\in\Sigma(C_0):\text{22212112112 and 21121121222 are not substrings of}~\underline{a}\}
\end{equation*}

Observe that
\begin{equation*}
    \lambda_0(2112112^*1222)>3.2879
\end{equation*}

The right extreme of this  interval is determined by $y=m(\overline{1211}2^*1222\overline{1121})=3.2872978\dots$, which corresponds to the minimum Markov value of a word containing $11211212221$.

We will set

\begin{gather*}
    \Sigma(B)=\{\underline{a}\in\Sigma(C_0):1121121222, 22111211221211211212, \\ 2111211221211211212211, 11121112112212112112122122 ~\text{and their transposes}\\\text{ are not substrings of}~\underline{a}\}.
\end{gather*}

Suppose that $\underline{b}\in\Sigma(B)$ satisfies $m(\underline{b})=\lambda_0(\underline{b})$. Since $\lambda_0(12^*2)<3.16$, we must have $b_{-1}b_0^*b_1=12^*1$. Since $\lambda_0(112^*11)<3.268$ we can assume that $b_{-1}b_0^*b_1b_2=12^*12$. Since 21212, 212111, 2121122, 12121 are forbidden it is forced to $12112^*122$. Since $\lambda_0(212112^*122)<3.28589$ we must continue as $112112^*122$. Since 1121121222 is forbidden we are forced to $112112^*1221$. Since $\lambda_0(112112^*12212)<3.2846$, $\lambda_0(1112112^*1221)<3.28518$, $\lambda_0(112112^*122111)<3.2856298$,\linebreak $\lambda_0(2112112^*1221122)<3.28598$, we should continue as $2112112^*1221121$. Since 2121122 is forbidden we are forced to $2112112^*12211211$. Since $\lambda_0(22112112^*1221)<3.28587$ and\linebreak $\lambda_0(112112112^*12211211)<3.28599$, we should continue as $212112112^*12211211$. Since 12121 and 1121121222 are forbidden we are forced to $12212112112^*12211211$. Since we have \linebreak $\lambda_0(12212112112^*122112112)<3.286015$, $\lambda_0(12212112112^*1221121111)<3.286026$, we should continue as $12212112112^*1221121112$. Since 22111211221211211212 and 212111 are forbidden in $\Sigma(B)$ we are forced to $12212112112^*122112111211$. Since 2111211221211211212211 is forbidden and $\lambda_0(212212112112^*1221121112112)<3.28602838$, we should continue as $212212112112^*1221121112111$. Since 11121112112212112112122122 is forbidden we must continue as $1212212112112^*1221121112111$. 

Since 12121, 21212, 212111, 2121122, 2221211211 are forbidden, a continued fraction in $\{1,2\}$ is maximized when
\begin{equation*}
    [0;1,2,1,2,2,2,1,2,1,1,2,1,2,2,1,2,1,1,2,1,1,2,1,2,2,1,1,\overline{2,1,1,1}]
\end{equation*}
On the other hand, since 212111 is forbidden, a continued fraction in $\{1,2\}$ that begins with $[0;1,1,1]$ is maximized when $[0;\overline{1,1,1,2}]$.

Therefore if $\underline{b}=\dots 1212212112112^*1221121112111\dots$, then $\lambda_0(\underline{b})$ is maximized when
\begin{equation*}
    m(\overline{1112}1122121121121221211212221211212212112112^*12211\overline{2111})=3.2860284\dots
\end{equation*}

\begin{remark}
Note that the middle word is a palindrome and the period 2111 is semi-symmetric.
\end{remark}

Hence we set $$x = m(\overline{1112}1122121121121221211212221211212212112112^*12211\overline{2111})=3.2860284\dots.$$

Let $m=m(\underline{a}) = \lambda_{0}(\underline{a})\in [x,y)$. From the arguments of the previous subsection, we need only consider the situation where we have two continuations $v^{1} = 1\ldots$ and $v^{2} = 2\ldots$ and $a_{N-1}a_{N} = 12$. In such a case, by considering all of the forbidden words giving $\Sigma(C)$, we get
\[[0;v^{2}] = [0;2,v^{2}_{2},\ldots]\leq[0;2,2,1,2,1,1,2,1,2,2,1,2,1,1,2,1,1,2,1,2,2,1,1,2,\overline{1,1,1,2}]\in K(B).\]
So, since $\lambda_{0}(122212^{*}112122) < 3.28589$ and $\lambda_{0}(122212112^{*}122) < 3.2858$, by arguments similar to the previous subsections, we can use the continuation $221211212212112112122112\overline{1112}\in\Sigma^{+}(B)$.

\vspace*{5pt}
\noindent{\Large\textbf{An interval containing Freiman's second example}}
\vspace*{5pt}

\subsection{Interval [3.29296, 3.29335)}\hfill

In \cite[Page 44]{Cusick-Flahive} it is mentioned that the list of Berstein intervals (see \Cref{app:ber}) does not cover Freiman's second example $\alpha_\infty=m(\overline{2122211}2^*12221121\overline{2})=3.293044265\dots\in M\setminus L$. Here we construct an interval that does.

We set
\begin{gather*}
    \Sigma(C)=\{\underline{a}\in\{1,2\}^{\Z}:\text{21212, 212111, 12121, 2221211221, 21222121122 and their }\\\text{ transposes are not substrings of}~\underline{a}\}
\end{gather*}

Observe that
\begin{align*}
    \lambda_0(22212^*11221)&>3.2935 \\
    \lambda_0(12122212^*1122)&>3.2934
\end{align*}

The right extreme of the interval can be pushed to 3.29335. Indeed, we claim that the minimum value of a word containing 121222121122 is
\begin{equation*}
    m(\overline{21121222}112122212112122212^*11\overline{22212112})=3.293351\dots
\end{equation*}
where we used the forbidden words 111212, 1221121222, 21212.

\begin{remark}
    The middle word is a palindrome and the period is non-semisymmetric of even length.
\end{remark}

Now let
\begin{equation*}
    \Sigma(B)=\{\underline{a}\in\Sigma(C):2211212221~\text{and 1222121122 are not substrings of}~\underline{a}\}
\end{equation*}

Suppose that $\underline{b}\in\Sigma(B_1)$ satisfies $m(\underline{b})=\lambda_0(\underline{b})$. Since $\lambda_0(12^*2)<3.16$, we must have $b_{-1}b_0^*b_1=12^*1$. Since $\lambda_0(112^*11)<3.268$ we can assume that $b_{-1}b_0^*b_1b_2=12^*12$. Since 21212, 212111, 12121 are forbidden it is forced to $2112^*122$. Since $\lambda_0(2112^*1221)<3.2921$, $\lambda_0(12112^*122)<3.2891$ and 2211212221, 1221121222 are forbidden in $\Sigma(B)$ we can assume that its equal to $222112^*12222$. Now we use the fact if 12121, 212111, 2221211221 are forbidden a continued fraction that begins with $[0;1,2,2,2,2]$ is maximized with $[0;1,2,\overline{2,2,2,2,1,2,1,1}]$ and similarly a continued fraction that begins with $[0;1,1,2,2,2]$ is maximized with $[0;\overline{1,1,2,2,2,2,1,2}]$. This shows that 

\begin{equation*}
    m(\underline{b})\leq m(\overline{21222211}2^*12222212\overline{11222212})=3.292954\dots
\end{equation*}

In fact this Markov value corresponds to the left endpoint of a gap. Recall the numbers
\[b_\infty = [2; \overline{1,1,2,2,2,1,2}] + [0; \overline{1, 2,2,2, 1,1, 2}] = 3.2930442439\ldots\]
and
\begin{gather*}
B_\infty = [2; 1, \overline{1, 2,2,2, 1, 2, 1,1, 2, 1,1, 2}] + [0; 1, 2,2,2, 1,1, 2, 1, 2,2,2,\\ 1,1, 2, 1, 2,2, \overline{1, 2,2,2, 1, 2, 1,1, 2, 1,1, 2}] = 3.2930444814\ldots.
\end{gather*}
discussed in Section~\ref{sec:goodprops}. It can be demonstrated that the minimum Markov value of a word containing 2211212221 is
\begin{equation*}
    m(\overline{2^*112221})=3.2930442439\dots=b_{\infty}.
\end{equation*}
Indeed, it follows from \cite[Lemma 3.7]{MminusL0353} and \cite[Lemma 3.8]{MminusL0353} that any bi-infinite sequence $\underline{c}\in\{1,2\}^\Z$ containing 2211212221 with $m(\underline{c})<B_\infty=3.2930444814\dots$ must contain  $221\overline{2112221}$. Now since
\begin{align*}
    \lambda_0(122121122212^*1122212112)&>b_\infty \\
    \lambda_0(2222121122212^*112221211222)&>b_\infty
\end{align*}
this subsequence must extend to $12221\overline{2112221}$. Applying \cite[Lemma 3.7]{MminusL0353} inductively yields that $\underline{c}=\overline{2112221}=b_\infty$ (in particular $b_\infty$ is isolated). So there is a gap between $3.292954\dots$ and $b_{\infty}$.

So we set 
\begin{align*}
[x,y) &= [m(\overline{21222211}2^*12222212\overline{11222212}),m(\overline{121221}22212^*112221211222\overline{2111})) \\
&=[3.292954\dots,3.293315\dots).
\end{align*}

Let $m = m(\underline{a}) = \lambda_{0}(\underline{a})<y$, and suppose that we have two continuations $v^{1} = 1\ldots$ and $v^{2} = 2\ldots$. If $a_{N} = 1$, since 111212 is forbidden, we have
\[[0;v^{1}] \geq [0;\overline{1,1,2,1}]\in K(B)\]
so since $\lambda_{0}(112^{*}11)<3.269$ we have the continuation $\overline{1121}\in\Sigma^{+}(B)$. So we need only consider the case $a_{N} = 2$.

If $a_{N-1}=1$, then since 12121, 21212, 212111 and 121222121122 are forbidden, we can choose the continuation $\overline{112122212112}$ which is allowed because $\lambda_0(12112^*122)<3.28907$.

So we need only consider the case where $a_{N-1}a_{N} = 22$. In such a case, since 21212, 212111, 2221211221 are forbidden, we have
\[[0;v^{2}] = [0;2,v^{2}_{2},\ldots] \leq [0;\overline{2,2,1,2,1,1,2,2}]\in K(B),\]
and since 
\begin{align*}
    \lambda_0(\ldots a_{N-1}a_N\overline{2212^*1122})&=[2;\overline{1,1,2,2,2,2,1,2}]+[0;1,2,2,2,2,\dots] \\
    &\leq[2;\overline{1,1,2,2,2,2,1,2}]+[0;1,2,\overline{2,2,2,2,1,2,1,1}] \\
    &=m(\overline{21222211}2^*12222212\overline{11222212})=x < m
\end{align*}
we can use the continuation $\overline{22121122}\in\Sigma^{+}(B)$.

\bigskip
\noindent{\Large{\textbf{An interval containing $\boldsymbol{t_1}$}}}
\vspace*{5pt}

Since $t_{1} = \min\{t\in\R:d(t) = 1\}$, we see that $d_{loc}$ starts below 1 in the interval below, but reaches 1 by the end of the interval.

\begin{remark}
We want to highlight that we found this good interval by constructing two transitive subshifts that work for one of Berstein's intervals.
\end{remark}

\subsection{Interval [3.33396, 3.33475)}\hfill

We set 
\begin{multline*}
    \Sigma(C)=\{\underline{a}\in\{1,2\}^{\Z}:21212, 21112121, 111112121, 211121222, 12111212, \\ 111212111, 21122121112  \text{ and their transposes are not substrings of } \underline{a}\}.
\end{multline*}

Observe that

\begin{align*}
    \lambda_0(212^*12)&>3.4 \\
    \lambda_0(21112^*121)&>3.35 \\
    \lambda_0(111112^*121)&>3.337 \\
    \lambda_0(21112^*1222)&>3.337 \\
    \lambda_0(121112^*12)&>3.335 \\
    \lambda_0(1112^*12111)&>3.335 \\
    \lambda_0(2112212^*1112)&>3.3347
\end{align*}

The minimum Markov value of a bi-infinite sequence containing 21122121112 is 
\begin{equation*}
    m(\overline{121111121222}2211121221112122112212^*111221211122\overline{222121111121})=3.3347525256\dots
\end{equation*}
where we used that 21212, 111112121, 21211121, 222121112, 1122121112221 are forbidden, as well as \Cref{lem:minimizing}. So we set $y = m(\overline{121111121222}2211121221112122112212^*111221211122\overline{222121111121})=3.3347525256\dots$.

\begin{remark}
The minimum is a palindrome and the period is semi-symmetric.
\end{remark}

Set
\begin{equation*}
    \Sigma(B)=\{\underline{a}\in\Sigma(C):121111212112, 2111212211~\text{and its transposes are not substrings of}~\underline{a}\}.
\end{equation*}

Suppose that $\underline{b}\in\Sigma(B)$ satisfies $m(\underline{b})=\lambda_0(\underline{b})$. Since $\lambda_0(12^*2)<3.16$, we must have $b_{-1}b_0^*b_1=12^*1$. Since $\lambda_0(112^*11)<3.268$ and 21212 is forbidden we can assume that $b_{-2}b_{-1}b_0^*b_1b_2=112^*12$. Since $\lambda_0(2112^*12)<3.32$ we must continue as $1112^*12$. If it continues as $1112^*121$, then since 21112121, 111112121, 21212, 111212111 are forbidden and $\lambda_0(1211112^*12112)<3.33383$ we must continue as $1211112^*12112$, which is forbidden in $\Sigma(B)$. Hence we must continue as $1112^*122$. Since $\lambda_0(11112^*122)<3.329$ and 211121222 is forbidden we must continue as $21112^*1221$. Since 12111212, 2111212211 are forbidden we must continue as $221112^*12212$. Hence we have that
\begin{equation*}
    \lambda_0(221112^*12212)\leq m(\overline{211111212221}112^*1221221211\overline{122212111112})=3.333958\dots
\end{equation*}
We set $x$ to be this value.

Let $m = m(\underline{a}) = \lambda_{0}(\underline{a})<y$ and suppose that we have two continuations $v^{1}=1\ldots$ and $v^{2} = 2\ldots$. If $a_{N} = 2$, then since 21212, 222121112 and 111112121 are forbidden, we have
\[[0;v^{2}] = [0;2,v^{2}_{2},\ldots] \leq [0;\overline{2,2,1,2,1,1,1,1,1,2,1,2}]\in K(B),\]
so that, since $\lambda_{0}(2212^{*}1111) < 3.3282$, we can choose the continuation $\overline{221211111212}\in\Sigma^{+}(B)$.

So we need only consider the case where $a_{N} = 1$. If $a_{N-1} = 2$, then, since 21212, 21211121, 22212112 and 111112121 are forbidden, we have
\[[0;v^{2}] = [0;2,v^{2}_{2},\ldots] \leq [0;2,2,1,2,1,1,\overline{1,2,2,2,1,2,1,1,1,1,1,2}]\in K(B)\]
so that we can take the continuation $221211\overline{122212111112}\in\Sigma^{+}(B)$, which is allowed because $\lambda_0(2212^*1111)<3.3282$ and 
\begin{equation*}
    \lambda_0(a_{N-1}a_{N}2212^*11\overline{122212111112})\leq m(\overline{211111212221}11212212212^*11\overline{122212111112})=3.333958\dots.
\end{equation*}

So, we are left considering the case of $a_{N-1}a_{N} = 11$. If $a_{N-2} = 2$, then since 21212, 21122121112 and 111112121 are forbidden, we have
\[[0;v^{2}] = [0;2,v^{2}_{2},\ldots] \leq [0;\overline{2,2,1,2,1,1,1,1,1,2,1,2}]\in K(B)\]
and so, since $\lambda_0(2212^*1111)<3.3282$, we have the continuation $\overline{221211111212}\in\Sigma^{+}(B)$. If $a_{N-2} = 1$, then since 111112121, 21212 and 222121112 are forbidden, we have
\[[0;v^{1}] = [0;1,v^{1}_{2}] \geq [0;\overline{1,1,2,1,2,2,2,1,2,1,1,1}]\in K(B)\]
and, since $\lambda_{0}(2212^{*}1111) < 3.3282$, we can use the continuation $\overline{112122212111}\in\Sigma^{+}(B)$.

\bigskip
\noindent{\Large{\textbf{An interval between $\boldsymbol{t_1}$ and $\boldsymbol{\sqrt{12}=3.4641\dots}$}}}
\vspace*{5pt}

Recall that $t_{1} = \min\{t\in\R : d(t) = 1\}$, so now that we are above $t_{1}$ we have that $d_{loc}(t) = 1$ for all $t\in L'\cap[\nu,\mu)$ where $[\nu,\mu)$ is a good interval. We are still in the situation where all of our Cantor sets will correspond to subshifts of $\{1,2\}^{\Z}$ since we remain below $\sqrt{12}$.

\subsection{Interval [3.359, 3.423)}\hfill

Consider $\Sigma(C)=\{\underline{a}\in\{1,2\}^{\Z}:\text{121212, 212121 are not substrings of}~\underline{a}\}$.

The minimum Markov value of a bi-infinite sequence containing 121212 is
\begin{equation*}
    m(\overline{121212212^*121})=3.42339101\dots
\end{equation*}
This Markov value was calculated using the algorithm discussed in Appendix~\ref{app:alg}.

\begin{remark}
The period is a palindrome.
\end{remark}

Indeed, since $\lambda_0(1212^*121)>3.44$, $\lambda_0(1212^*1222)>3.427$, $\lambda_0(1212^*12211)>3.424$,\linebreak $\lambda_0(1212^*122122)>3.4236$, $\lambda_0(1212^*1221211)>3.4234$, $\lambda_0(1212^*1221212)>3.423397$ are all greater than the above candidate the subword must extend to the word 121212212121 which has the form $\beta^T2\theta 2\beta$ with $\beta=121$ and the even palindrome $\theta=1221$. Hence by \Cref{lem:minimizing} and the forbidden words 1212121, 12121222, 121212211, 1212122122, 12121221211, 121212212122, the minimum is attained at the above candidate. So we can set $y = m(\overline{121212212^*121})=3.42339101\dots$.

Let $\Sigma(B)=\{\underline{a}\in\{1,2\}^{\Z}:\text{21212 is not a substring of}~\underline{a}\}$. From Cusick-Flahive~\cite[Chapter 5, Table 1]{Cusick-Flahive}, we know that $m(\underline{b})\leq m(\overline{1112^*12})=3.35871\ldots$ for all $\underline{b}\in\Sigma(B)$. So we set $x = m(\overline{1112^*12})=3.35871\ldots$.

Suppose there are two continuations $v^{1}=1\ldots$ and $v^{2} = 2\ldots$. Since, 121212 is forbidden, we have
\[[0;v^{1}] = [0;1,v^{1}_{2},\ldots]\geq [0;\overline{1,1,2,1,2,1}]\in K(B),\]
so that, since $\lambda_{0}(112^{*}121) \leq x$, we can use the continuation $\overline{112121}\in\Sigma^{+}(B)$.

\bigskip

\noindent{\Large{\textbf{Intervals after $\boldsymbol{\sqrt{12}=3.464\dots}$}}}
\vspace*{5pt}

Since we are now above $\sqrt{12}$, our shifts will be subshifts of $\{1,2,3\}^{\Z}$ not contained in $\{1,2\}^{\Z}$; that is, our sequences must now contain 3's.

\subsection{Interval [$\boldsymbol{\sqrt{12}}$, 3.8465)}\hfill 

This interval was done in the paper \cite{MMPV} but we will give the details here for completeness. The subshifts used are $\Sigma(C)=\{\underline{a}\in\{1,2,3\}^{\Z}:\text{13, 31 are not substrings of}~\underline{a}\}$ and $\Sigma(B)=\{1,2\}^{\Z}$. The right extreme of the interval corresponds to the minimum Markov value of a bi-infinite sequence containing 13, namely $m(\overline{3^*113})=3.846546\dots$. In fact this Markov value is isolated: note that the words 131, 132, 312, 313, 1332, 1333, 3111, 1133112 all produce Markov values greater than 3.8488 so if a bi-infinite sequence contains 13 then it must be equal to $\overline{1331}$ and since 32323, 323222, 323221, 223232233, 2232322321, 22323223223 all give Markov values greater than 3.84656, so if a bi-infinite sequence does not contain 13 and neither of the previous forbidden words, the subword 32 must extend (up to transposition) to either 332, 22322, 1232, 3323223, 12323223, 32232322323 (all which have values at most 3.846511), 322323223222 which has value at least $m(\overline{12}2322232232322323^*22322232\overline{21})=3.846553\dots$ or 322323223221 which has value at most $m(\overline{21}1223223^*23223221\overline{12})=3.8465437\dots$.

The minimum Markov value of a bi-infinite sequence containing 1133112 is 
\begin{equation*}
    m(\overline{21}13^*31\overline{12})=3.856886\dots
\end{equation*}
where we used that 1113 and 312 are forbidden, as well as \Cref{lem:minimizing}.

The minimum Markov value of a bi-infinite sequence containing 1332 is 
\begin{equation*}
    m(\overline{12}23233113^*3232\overline{21})=3.856958\dots
\end{equation*}
Indeed, since the words 312, 313, 1113, 21132, 131, 231, 33311, 13311332 all produce larger Markov values than the above candidate, the subword 1332 must extend to 23311332. Now use \Cref{lem:minimizing} and the forbidden words 132, 312, 32323.

\begin{remark}
In both cases the minimum has middle word palindrome and period semi-symmetric.
\end{remark}

If we have two continuations $v^{1} = 1\ldots$ and $v^{2} = 3\ldots$, then we can choose $v_{B} = \overline{21}$ since we will have
\[[0;v^{1}]>[0;v_{B}]>[0;v^{2}],\]
so that, since the Markov value must occur at a position of the sequence with a 3 and $v_{B}\in\{1,2\}^{\N}$, $v_{B}$ is allowed.

For this reason, for the remainder of the section, we need only consider continuations of the form $v^{1} = 1\ldots$ and $v^{2} = 2\ldots$, or continuations $v^{1} = 2\ldots$ and $v^{2} = 3\ldots$.

In the former case, since 13 and 31 are forbidden, we will have
\[[0;v^{1}] \geq [0;1,\overline{1,2}]\in K(B)\]
so the continuation $v_{B} = 1\overline{12}\in\Sigma^{+}(B)$ can be used. In the latter case, we will have
\[[0;v^{1}] \geq [0;\overline{2,1}]\in K(B)\]
so the continuation $v_{B} = \overline{21}\in\Sigma^{+}(B)$ can be used.

\subsection{Interval [3.873, 3.930691)}\label{subsec:39306}\hfill 

Consider $\Sigma(C)=\{\underline{a}\in\{1,2,3\}^{\Z}:\text{132, 231, 312, 213, 313, 131 are not substrings of}~\underline{a}\}$. Note that
\begin{equation*}
    \lambda_0(3^*12)\geq [3;1,2,\overline{3,1}]+[0;\overline{3,1}]>3.95
\end{equation*}
\begin{equation*}
    \lambda_0(3^*13)\geq[3;1,3,\overline{3,1}]+[0;\overline{3,1}]>4.02
\end{equation*}
\begin{equation*}
    \lambda_0(13^*1)\geq[3;1,\overline{1,3}]+[0;1,\overline{1,3}]>4.11
\end{equation*}
Assuming that 312, 313 and 131 are forbidden, we have that the least Markov value of a word containing 231 is $m(\overline{12}3^*113\overline{21})=3.930691\dots$.
\begin{remark}
Note that the middle word 3113 is palindromic.
\end{remark}

Indeed, since
\begin{align*}
    \lambda_0(13^*23)&\geq [3;2,3,\overline{3,1}]+[0;1,\overline{1,3}]>3.99 \\
    \lambda_0(13^*22)&\geq [3;2,2,\overline{3,1}]+[0;1,\overline{1,3}]>3.967 \\
\end{align*}
we have that the subword $23^*1$ must continue as $123^*11$. Since $312$ is forbidden, the existence of the above word implies that it must continue as $123^*11321$ so it is $w_1^T123^*11321w_2$ for some $w_1,w_2\in\{1,2,3\}^\N$. By the \Cref{lem:minimizing} we have that  this minimum value is attained when $[0;2,1,w_1]=[0;2,1,w_2]$ is minimum, which implies $w_1=w_2=\overline{21}$ since $312$ is forbidden.

We can set $y=m(\overline{12}3^*113\overline{21})=3.930691\dots$.

On the other hand we can use $$\Sigma(B)=\{\underline{a}\in\{1,2,3\}^{\Z}:\text{13, 31 are not substrings of}~\underline{a}\}.$$ We have that $m(\underline{b})\leq x=m(\overline{32})=3.8729\dots$.

If we have two continuations $v^{1} = 1\ldots$ and $v^{2} = 2\ldots$, then, since 231 is forbidden, we have
\[[0;v^{2}]\leq[0;\overline{2,3}]\in K(B)\]
so the continuation $v_{B} = \overline{23}\in\Sigma^{+}(B)$ can be used since $\lambda_{0}(23^{*}2)\leq x$. In the case $v^{1} = 2\ldots$ and $v^{2} = 3\ldots$, since 213 is forbidden, we have
\[[0;v^{1}]\geq [0;\overline{2,1}]\in K(B)\]
so we can use the continuation $v_{B} = \overline{21}\in\Sigma^{+}(B)$.

\subsection{Interval [3.93616, 3.943767)}\hfill

We will build this interval as the union of the intervals [3.93616, 3.9373), [3.93727, 3.94) and [3.93931, 3.943767).

First, consider the subshift
\begin{gather*}\label{eq:sigmac39656}
    \Sigma(C_0)=\{\underline{a}\in\{1,2,3\}^{\Z}:\text{313, 131, 1322, 1323, 312, 11132 and}\\ \text{ their transposes are not substrings of}~\underline{a}\}
\end{gather*}

As above, if $v^{1} = 2\ldots$ and $v^{2} = 3\ldots$, we can use the continuation $v^{B} = \overline{21}$. So we need only consider the case $v^{1} = 1\ldots$ and $v^{2} = 2\ldots$. If $a_{N} \in \{2,3\}$ then, since 3231 and 2231 are forbidden, we can show that
\[[0;v^{2}] \leq [0;\overline{2,3}]\]
so the continuation $v_{B} = \overline{23}$ can again be used.

So we are left to consider the case $a_{N} = 1$. Since 312 is forbidden and $v^{2}$ begins with 2, we must have $a_{N-1}a_{N} \in\{21,11\}$.

At this point, we have
\[[0;v^{1}] \geq [0;1,1,\overline{3,3,3,1,1,1}]\]
and
\[[0;v^{2}] \leq [0;\overline{2,3,1,1}]\]
so we will construct several $\Sigma(C)$ and $\Sigma(B)$ by forbidding longer subwords of these periods.

\subsubsection{Interval [3.93616, 3.9373)}\hfill

Consider the subshift
\begin{gather*}
    \Sigma(C)=\{\underline{a}\in\Sigma(C_0):23112, 11231, 3111333, 3211133311121 ~\text{and their}\\\text{ transposes are not substrings of}~\underline{a}\}.
\end{gather*}

Observe that

\begin{align*}
    \lambda_0(1123^*112)&>3.954 \\
    \lambda_0(1123^*1133)&>3.944 \\
    \lambda_0(1123^*1132)&>3.941 \\
    \lambda_0(2123^*112)&>3.941 \\
    \lambda_0(31113^*33)&>3.94 \\
    \lambda_0(32111333^*11121)&>3.9375
\end{align*}

Note that if 23112 is a subword, then it must extend to either 2123112 or 1123112 because 1323, 1322 and 312 are forbidden in $\Sigma(C_0)$. Similarly 11231 must extend to 11231132, 11231133 or 1123112 because 313, 312 and 131 are forbidden in $\Sigma(C_0)$.

We set
\begin{gather*}
    \Sigma(B)=\{\alpha\in\Sigma(C): 33311121, 321231133, 2212311331, \\ 1211133311122, 2211133311122~\text{and their transposes are not substrings of}~\alpha\}
\end{gather*}

We claim that for all $\underline{b}\in\Sigma(B)$
\begin{equation*}
    m(\underline{b})\leq m(\overline{12}33111332111333^*1112\overline{23})=3.936154\dots
\end{equation*}

Indeed, let $\underline{b}\in\Sigma(B)$ such that $m(\underline{b})=\lambda_0(\underline{b})$ assumes the maximum of $\Sigma(B)$, so we have either $b_{-1}b_0^*b_1$ is equal to $33^*1$ or $23^*1$. If it is equal to $33^*1$, then since 313 and 312 are forbidden it is forced to $33^*11$. Since $\lambda_0(33^*11b_3)<3.9$ for $b_2\in\{2,3\}$ we can assume that it continues as $33^*111$. Since $\lambda_0(33^*1111)<3.928$, $\lambda_0(233^*1112)<3.933$, $\lambda_0(233^*11133)<3.9365$ and $\lambda_0(133^*11)<3.923$ we see that we can assume that it continues as $333^*111$. Since $\lambda_0(333^*1111)<3.928$, $\lambda_0(333^*11123)<3.9357$ and 3331113, 33311121 are forbidden, we can assume it continues as $333^*11122$. Since $\lambda_0(3333^*11122)<3.934$, $\lambda_0(2333^*11122)<3.9349$, $\lambda_0(311333^*11122)<3.9358$, $\lambda_0(211333^*11122)<3.936$ and 213, 313, 3111333 are forbidden, we must continue with either $2111333^*11122$ or $\lambda_0(1111333^*11122)\leq\lambda_0(\overline{31}1111333^*1112\overline{23})<3.9361$ (where we used that 2231 is forbidden). Since 1211133311122 and 2211133311122 are forbidden in $\Sigma(B)$ we can assume the maximum has the form $32111333^*11122$. Finally, we use that if 11132, 131, 213, 313, 3111333 are forbidden then a continued fraction in $\{1,2,3\}$ that begins with $[3;3,3,1,1,1,2,3]$ is maximized with $[3;3,3,1,1,1,2,3,3,1,1,1,3,3,\overline{2,1}]$ and similarly if 2231 is forbidden a continued fraction in $\{1,2,3\}$ that begins with $[0;1,1,1,2,2]$ is maximized with $[0;1,1,1,2,\overline{2,3}]$.

Now assume that $b_{-1}b_0^*b_1=23^*1$. Since 2231, 3231, 312, 313, 23112, 23111, 11231 are forbidden, we have that $23^*1$ is forced to $2123^*113$. Since $\lambda_0(2123^*1132)<3.93605$ and 131 is forbidden we must continue as $2123^*1133$. Since $\lambda_0(12123^*1133)<3.9343$ and 321231133, 2212311331 are forbidden, and moreover $\lambda_0(22123^*11333)<3.9359$, $\lambda_0(22123^*11332)<3.9361$, we see that this continuations leads to smaller Markov values.

Recall that we have $v^{1} = 1\ldots$, $v^{2} = 2\ldots$, and $a_{N-1}a_{N}\in\{21,11\}$. If $a_{N-1} = 1$, since 11231 and 3231 are forbidden, we have
\[[0;v^{2}]\leq [0;\overline{2,3}]\]
and the continuation $v_{B} = \overline{23}$ can again be used.

So we have $a_{N-1}a_{N} = 21$. Since 131, 313, 312, 23111, 23112 and 3111333 are forbidden, we have
\[[0;v^{2}]\leq[0;2,3,1,1,3,3,1,1,1,3,3,\overline{2,1}]\]
and, since $\lambda_{0}(a_{N-2}2123^{*}113311133\overline{21})<3.3936$ when $a_{N-2}\in\{1,2\}$, we can use the continuation $v_{B} = 23113311133\overline{21}$ if $a_{N-2}\neq 3$. Otherwise, since 11132, 313, 312, 3331113, 3211133311121 and 2231 are forbidden, we have
\[[0;v^{1}]\geq [0;1,1,3,3,3,1,1,1,2,\overline{2,3}]\]
so that, since $\lambda_0(2111333^*11122)<3.9362$, the continuation $v_{B} = 113331112\overline{23}$ can be used.

\subsubsection{Interval [3.93727, 3.94)}\hfill

Now, consider the subshift
\begin{equation*}
    \Sigma(C)=\{\underline{a}\in\Sigma(C_0):23112, 11231, 3111333 ~\text{and their transposes are not substrings of}~\underline{a}\}.
\end{equation*}
Observe that
\begin{equation*}
    \lambda_0(31113^*33)>3.94
\end{equation*}

Now we set $\Sigma(B)$ so that the above argument works. We choose
\begin{equation*}
    \Sigma(B)=\{\underline{a}\in\Sigma(C): 33311121, 13331112  ~\text{are not substrings of}~\underline{a}\}
\end{equation*}

We claim that for all $\underline{b}\in\Sigma(B)$
\begin{equation}\label{eq:maxsigmab394}
    m(\underline{b})\leq m(\overline{12}33111331132123113311133\overline{21}) =: x = 3.93726\dots
\end{equation}

Indeed, let $\underline{b}\in\Sigma(B)$ such that $m(\underline{b})=\lambda_0(\underline{b})$ assumes the maximum of $\Sigma(B)$, so we have either $b_{-1}b_0^*b_1$ is equal to $33^*1$ or $23^*1$. If it is equal to $33^*1$, then since 313 and 312 are forbidden it is forced to $33^*11$. Since $\lambda_0(33^*11b_3)<3.9$ for $b_2\in\{2,3\}$ we can assume that it continues as $33^*111$. Since $\lambda_0(33^*1111)<3.928$, $\lambda_0(233^*1112)<3.933$, $\lambda_0(233^*11133)<3.9365$ and $\lambda_0(133^*11)<3.923$ we see that we can assume that it continues as $333^*111$. Since 3331113 is forbidden, we can assume it continues as $333^*1112$. Since $\lambda_0(b_{-3}333^*1112b_5)<3.935$ for $b_{-3},b_5\in\{2,3\}$ and 13331112 and 33311121 are forbidden, we can see that this continuations are smaller than the above candidate.

Now assume that $b_{-1}b_0^*b_1=23^*1$. Since 2231, 3231, 312, 313, 23112, 23111, 11231 are forbidden, we have that $23^*1$ is forced to $2123^*113$. If 313, 312, 11132, 3111333 are forbidden a continued fraction in $\{1,2,3\}$ that begins in $[0;1,1,3]$ is maximized with $[0,1,1,3,3,1,1,1,3,3,\overline{2,1}]$ and similarly since 313, 312, 23111, 23112, 3111333 are forbidden, a continued fraction in $\{1,2,3\}$ that begins with $[0;2,1,2]$ is maximized with $[0;2,1,2,3,1,1,3,3,1,1,1,3,3,\overline{2,1}]$.

The continuations can be argued as in the case of the previous interval. However, in the case of $a_{N-1}a_{N} = 21$ we can now use the continuation $v_{B} = 23113311133\overline{21}$ for any $a_{N-2}$ as it can be checked that
\[\lambda_{0}(a_{N-2}2123^{*}113311133\overline{21})\leq x\]
for any $a_{N-2}\in\{1,2,3\}$.

\subsubsection{Interval [3.93931, 3.943767)}\hfill

Now let
\begin{equation*}
    \Sigma(C)=\{\underline{a}\in\Sigma(C_0):1113331113, 3111333111 ~\text{are not substrings of}~\underline{a}\}.
\end{equation*}
Observe that
\begin{equation*}
    \lambda_0(111333^*1113)>3.9435
\end{equation*}

The minimum Markov value of a word containing 21113331113 is the same as that of a word containing  the palindrome 2111333111331113331112 which is
\begin{equation*}
    m(\overline{12}31132111333^*1113311133311123113\overline{21})=3.943767\dots
\end{equation*}
where we used the forbidden words 11132, 1123112, 11231133, 213, 313. So we can choose $y = 3.943767\dots$.

\begin{remark}
The middle word is palindrome an the period is semi-symmetric.
\end{remark}

Now we define $\Sigma(B)$ by
\begin{equation*}
    \Sigma(B)=\{\underline{a}\in\Sigma(C): 23112, 11231, 3331113 ~\text{and their transposes are no substring of}~\underline{a}\}.
\end{equation*}

We claim that $m(\underline{b})\leq m(\overline{21}11333^*11\overline{12})=3.939301\dots$ for all $\underline{b}\in\Sigma(B)$. Indeed, let $\underline{b}\in\Sigma(B)$ such that $m(\underline{b})=\lambda_0(\underline{b})$ assumes the maximum of $\Sigma(B)$, so we have either $b_{-1}b_0^*b_1$ is equal to $33^*1$ or $23^*1$. Since 2231, 3231, 312, 313, 23112, 23111, 11231 are forbidden in $\Sigma(B)$, we have that $23^*1$ is forced to $\lambda_0(2123^*113)<3.937672$. Hence we can assume that $b_{-1}b_0^*b_1=33^*1$. Since 312 and 313 are forbidden we are forced to $33^*11$. Since $\lambda_0(b_{-2}33^*11)<3.9396$ for $b_{-2}\in\{1,2\}$ and $\lambda_0(33^*11b_{3})<3.9$ for $b_{3}\in\{2,3\}$, we can assume that it continues as $333^*111$. Since $\lambda_0(333^*1111)<3.928$ and 3331113 is forbidden, we can assume that it continues as $333^*1112$. Since $\lambda_0(b_{-2}333^*1112)<3.9385$ for $b_{-2}\in\{2,3\}$ and since 312 and 313 are forbidden, we can assume that it continues as $11333^*1112$. Since 312 is forbidden, a continued fraction in $\{1,2,3\}$ that begins with $[0;1,1,1,2]$ is maximized with $[0;1,1,\overline{1,2}]$. Since $\lambda_0(b_{-5}11333^*11\overline{12})<3.939$ for $b_{-5}\in\{2,3\}$, we can assume that it continues as $111333^*1112$. Finally since 3331113 and 312 are forbidden, a continued fraction in $\{1,2,3\}$ that begins with $[3;3,3,1,1,1]$ is maximized with $[3;3,3,1,1,\overline{1,2}]$.

Recall that we need only consider continuations of the form $v^{1} = 1\ldots, v^{2} = 2\ldots$, with $a_{N-1}a_{N} \in\{21,11\}$. Therefore, since 131, 11132, 313, 312 and 1113331113 are forbidden, we have
\[[0;v^{1}]\geq [0;1,1,3,3,3,1,1,\overline{1,2}]\]
and, since $\lambda_0(11333^*111212)<3.9394$ and $\lambda_0(21113^*3311\overline{12})<3.93987$, we can use the continuation $v_{B} = 1133311\overline{12}$.

\subsection{Interval [3.944054, 3.971606)}\label{eq:Sigma_C39857}\hfill

Before specialising to this interval, we consider subshifts for the region [3.94405, 3.9857).

Consider the subshift $\Sigma(C_0)$ defined by the forbidden subwords
\begin{gather}
    131, 313, 2132, 111322, 1323, 1213, 33312,  23312, 211132, 311132, 1123111, 223112. \nonumber
\end{gather}

Observe that
\begin{align*}
    \lambda_0(213^*2)&> 4.05 & \lambda_0(1123^*111)&>3.987 \\ 
    \lambda_0(1113^*22)&>4.01 & \lambda_0(21113^*21)&>3.987 \\ 
    \lambda_0(223^*1123)&>3.998 & \lambda_0(223^*112)&>3.9855 \\ 
    \lambda_0(333^*12)&>3.996 & \lambda_0(1213^*3)&> 3.982 \\
    \lambda_0(31113^*2)&>3.996  & \lambda_0(233^*12)&>3.984
\end{align*}
To determine precisely the right extreme of this region we have to compute explicitly the minimum Markov values of some of the forbidden words above. The word that (allegedly) is determining the right extreme of the region is 12133, but we will see that in fact it is 23312.

We choose
\begin{equation}\label{eq:Sigma_B94405}
    \Sigma(B_0)=\{\alpha\in\Sigma(C_0):  312, 23112, 11231  ~\text{and their transposes are not substrings of}~\alpha\}
\end{equation}

We claim that $m(\underline{b})\leq m(\overline{111333^*})=3.944053\dots$ for all $\underline{b}\in\Sigma(B_{0})$. Indeed, let $\underline{b}\in\Sigma(B_{0})$ such that $m(\underline{b})=\lambda_0(\underline{b})$ assumes the maximum of $\Sigma(B_{0})$, so we have either $b_{-1}b_0^*b_1$ is equal to $33^*1$ or $23^*1$. If $b_{-1}b_0^*b_1=23^*1$, then since 312, 313, 2231, 3231 are forbidden it is forced to $123^*11$. Since 23111, 23112, 11231, 131, 312 are forbidden in $\Sigma(B_{0})$, it is forced to $\lambda_0(2123^*113)<3.9377$.

If $b_{-1}b_0^*b_1=33^*1$ then since 312 and 313 are forbidden in $\Sigma(B_{0})$ it must extend to $33^*11$. If 131, 311132, 312 are forbidden, a continued fraction in $\{1,2,3\}$ that begins with $[3;1,1]$ is maximized with $[3;\overline{1,1,1,3,3,3}]$ and similarly a continued fraction in $\{1,2,3\}$ that begins with $[3;3]$ is maximized with $[3;\overline{3,3,1,1,1,3}]$. This proves the claim.

Now that we have chosen $\Sigma(B_0)$ and $\Sigma(C_0)$, we shall consider possible continuations $v^{1}$ and $v^{2}$ for strings in $\Sigma(C_{0})$.

Firstly, suppose that $v^{1} = 2\ldots$ and $v^{2} = 3\ldots$. If $a_{N} = 1$, then since $1213$ is forbidden, we have
\[[0;v^{1}] = [0;2,\ldots] \geq [0;\overline{2,1}]\in K(B_{0}).\]
If $a_{N} \neq 1$, then since 313, 23312, 33312, 131 and 311132 are forbidden, we have
\[[0;v^{2}] = [0;3,\ldots] \leq [0;3,3,\overline{1,1,1,3,3,3}]\in K(B_{0}),\]
which is allowed since $\lambda_{0}(33^{*}11) \leq m(\overline{333^{*}111}) = 3.94405\ldots$.

So we now consider $v^{1} = 1\ldots$ and $v^{2} = 2\ldots$. If $v^{2} = 22\ldots$ or $v^{2} = 21\ldots$, then
\[[0;v^{2}] < [0;\overline{2,3}]\in K(B_{0}),\]
and so we can assume that $v^{2} = 23\ldots$. If $v^{1} = 13\ldots$, $v^{1} = 12\ldots$, $v^{1} = 111\ldots$ or $v^{1} = 112\ldots$, then
\[[0;v^{1}] \geq [0;1,1,\overline{2,1}]\in K(B_{0}),\]
and so we can assume that $v^{1} = 113\ldots$. If $v^{1}$ continues as $v^{1} = 1133\ldots$ then, since 313, 33312, 131 and 311132 are forbidden, we have that
\[[0;v^{1}] = [0;1,1,3,3,\ldots] \geq [0;1,1,\overline{3,3,3,1,1,1}]\in K(B_{0})\]
which is allowed by arguments from above. Moreover, if $v^{2} = 232\ldots$ or $v^{2} = 233\ldots$, then we have 
\[[0;v^{2}] \leq [0;2,\overline{3,2}]\in K(B_{0}),\]
which is again allowed.

So in summary, and since 131, 313 and 2312 are forbidden, we must have
\begin{equation}\label{eq:reg3982}
    v^{1} = 1132\ldots\,\,\,\,\,\,\text{and}\,\,\,\,\,\,v^{2} = 2311\ldots.
\end{equation}

Now specialising to the interval [3.944054, 3.971606). Consider the subshift $\Sigma(C)$ defined by the forbidden subwords
\begin{gather*}
    131, 313, 2132, 1323, 1213, 33312,  23312, 23111, 2231.
\end{gather*}

Observe that
\begin{align*}
    \lambda_0(213^*2)&> 4.05 & \lambda_0(1213^*3)&> 3.982  \\ 
    \lambda_0(333^*12)&>3.996 & \lambda_0(223^*1)&>3.967 \\
    \lambda_0(233^*12)&>3.984 & \lambda_0(23^*111)&>3.967
\end{align*}

For this interval we keep $\Sigma(B)=\Sigma(B_0)$ defined in \eqref{eq:Sigma_B94405}. 

The minimum Markov value of a word containing 23111 is
\begin{equation*}
    m(\overline{12}22111133123^*111132133111122\overline{21})=3.9716067\dots
\end{equation*}

\begin{remark}
The middle word is palindrome and the period is semi-symmetric.    
\end{remark}

Indeed, since 1323, 231112, 231113 are forbidden and $\lambda_0(223^*111)>4$, we must extend the subword to 1231111. Since  $\lambda_0(1123^*111)>3.987$, $\lambda_0(2123^*111)>3.974$ are greater than the above candidate and 131, 2312, 33312, 23312 are forbidden, we must extend as 1331231111. Since $\lambda_0(23^*11111)>3.974$, $\lambda_0(33123^*11112)>3.973$, $\lambda_0(133123^*111133)>3.9717$ and 131 is forbidden, we must extend as 133123111132. By the previous arguments we further extend to 1331231111321331. Since 313 is forbidden and $\lambda_0(13312311113^*213312)>3.9717$ we must extend as 113312311113213311. Since $\lambda_0(2133^*113)>3.973$, $\lambda_0(2133^*112)>3.973$ we must extend as 11133123111132133111. Thus the word has the form $w_1^T\beta^t3\theta3\beta w_2$ where $\theta=1111$ is an even palindrome, $\beta=2133111$ and $w_1,w_2\in\{1,2\}^\N$. Hence by \Cref{lem:minimizing} this is minimized when $w_1=w_2$ and $[0;\beta,w_2]$ is minimal. Using the fact that $\lambda_0(11113213^*311113)>3.9717$, $\lambda_0(311113213^*3111121)>3.97164$, $\lambda_0(12311113213^*31111223)>3.9716069$ and the fact that 22213 is forbidden (use 131, 2132, 21333, 21332, 313 and $\lambda_0(2213^*311)>3.976$, $\lambda_0(22213312)>3.974$), we see that $w_1=w_2=122\overline{21}$, which confirms the candidate above.

The minimum value of a word containing $2231$ is
\begin{equation*}
    m(\overline{23^*1132})=3.97402\dots
\end{equation*}

\begin{remark}
The minimum is a palindrome.
\end{remark}

Indeed, since 313, 2312 are forbidden, $\lambda_0(223^*111)>4$ and  $\lambda_0(223^*112)>3.985$, we must extend the subword to 223113. Since $\lambda_0(223^*1133)>3.975$ and since 131, 1323 is forbidden, we must extend the subword to 2231132. Note that if 313, 2312, 23111, 223112 are forbidden, a continued fraction that begins with $[0;2,2]$ is minimized with $[0;\overline{2,2,3,1,1,3}]$. In particular
\begin{equation*}
    \lambda_0(223^*11321)\geq [3;\overline{2,2,3,1,1,3}]+[0;1,1,3,2,1,\dots]>[3;\overline{2,2,3,1,1,3}]+[0;1,1,3,2,2,\dots] 
\end{equation*}
so we must extend the subword to 22311322 since 1323 is forbidden. Thus the word has the form $w_1^T\beta^t3\theta3\beta w_2$ where $\theta=11$ is an even palindrome, $\beta=22$ and $w_1,w_2\in\{1,2\}^\N$. Hence by \Cref{lem:minimizing} this is minimized when $w_1=w_2$ and $[0;\beta,w_2]$ is minimal, which proves that indeed the above candidate is the minimum. 

Recall that the only continuations that need to be considered at this stage are $v^{1} = 1132\ldots$ and $v^{2} = 2311\ldots$.

Since 3231 and 22311 are forbidden and $a_N v^{2} = a_N 2311\ldots$ we must have $a_N=1$. However, by looking at $a_{N-1}a_N v^{1}=a_{N-1}11132\ldots$ we see that we cannot have these continuations at all since 11132 is forbidden.

\subsection{Interval [3.97995, 3.9857)}\hfill

For this interval we keep $\Sigma(C)=\Sigma(C_0)$ defined in Subsection~\ref{eq:Sigma_C39857}.

We choose
\begin{multline*}
    \Sigma(B)=\{\underline{a}\in\Sigma(C):  23312, 33121, 2231, 23111, 213312, 1133122 \\ ~\text{and their transposes are not substrings of}~\underline{a}\}
\end{multline*}

From the inequalities of Subsection~\ref{eq:Sigma_C39857} it seems that the right border of this interval would be 3.984. To determine the exact right border, we will compute explicitly the minimum Markov values of some of the forbidden words in $\Sigma(C)$ after showing that there exists indeed a good interval with these choices of $\Sigma(B)$ and $\Sigma(C)$.

We claim that $m(\underline{b})\leq m(\overline{113^*2})=3.97994\dots$ for all $\underline{b}\in\Sigma(B)$. Indeed, let $\underline{b}\in\Sigma(B)$ such that $m(\underline{b})=\lambda_0(\underline{b})$ assumes the maximum of $\Sigma(B)$, so we have either $b_{-1}b_0^*b_1$ is equal to $33^*1$ or $23^*1$. As before if $b_{-1}b_0^*b_1=33^*1$ extends to $33^*11$, then since 131, 311132, 33312 are forbidden we will have that $m(\underline{b})\leq m(\overline{111333^*})=3.944\dots$. If $33^*1$ extends to $33^*12$, then since 33312, 23312, 213312, 313 are forbidden it must extend to $1133^*12$. Then since 33121 and 1133122 are forbidden in $\Sigma(B)$ it is forced to $\lambda_0(1133^*123)<3.9786$.

Now assume that $b_{-1}b_0^*b_1=23^*1$. Since 313, 2312, 3231, 22311 are forbidden in $\Sigma(B)$, it must continue as $123^*11$. Since $\lambda_0(123^*113)<3.958$ and since 123111 is forbidden, we can assume it continues as $123^*112$. If 313, 2312, 31123111 are forbidden, a continued fraction that begins with $[3;1,1,2]$ is maximized with $[3;\overline{1,1,2,3}]$ and similarly if 131, 1323, 3211322 are forbidden, if it begins with $[3;2,1]$ then it is maximized with $[3;\overline{2,1,1,3}]$. This finishes the claim.

Now that we have chosen $\Sigma(B)$ and $\Sigma(C)$, by \eqref{eq:reg3982} we have  $v^1=1132\dots$ and $v^2=2311\dots$. Since 3231 is forbidden we have $a_N\neq 3$. If $a_N=1$, then since 311132, 211132, we have $a_{N-1}=1$. Since 1123111, 313, 2312 are forbidden, we see that $v^2$ is connecting to $\overline{2311}$ which is allowed because if 131, 1323, 211322 are forbidden a continued fraction that begins with $[3;2,1,1]$ is maximized with $[3;\overline{2,1,1,3}]$, so $\lambda_0(1123^*11\overline{2311})\leq m(\overline{1123^*})=3.97994\dots$ 

Now assume $a_N=2$. Since 1323 and 211322 are forbidden we must continue as $v^1=11321\dots$. Since 131, 1323, 211322 are forbidden we see that $v^1$ is connecting to $\overline{1132}$ which is allowed because if 313, 2312, 1123111 are forbidden, a continued fraction that begins with $[3;1,1,2]$ is maximized with $[3;\overline{1,1,2,3}]$, so $\lambda_0(2\overline{113^*2})\leq m(\overline{113^*2})=3.97994\dots$.

The right extreme of the interval can be pushed to 3.9857. We claim that the minimum value of a word containing $23312$ is
\begin{equation*}
    m(\overline{23}3^*123111113213312\overline{23})=3.98575\dots
\end{equation*}

\begin{remark}
This minimum is very asymmetric in the sense that the middle is not palindrome.
\end{remark}

Note that all the words defining $\Sigma(C)$ but 12133 and 223112 automatically produce larger values than the above candidate. The minimum of 223112 is easy to compute because 1213, 2132, 111322 are forbidden and corresponds to $m(\overline{21}13^*2231\overline{12})=3.9874\dots$. The minimum of the word 12133 necessarily has to extend to 112331211 (otherwise $\lambda_0$ is bigger than some candidate containing 12133). Note that $\lambda_0(11213^*31211)>3.987$ is bigger than the minimum of 23312. Similarly 223112 has to extend to 22311213, which extends to 223112133 by the forbidden words 131 and 2132, however $\lambda_0(223^*112133)>3.986$ is bigger than the minimum of 23312, so the right endpoint of the interval is determined by 23312.

First note that $\lambda_0(23^*312)<3.72$. Clearly to minimize $233^*12$ to the left, the smallest possible continuation is $\overline{23}3^*12\dots$ which is contained in the above candidate, so we only have to minimize to the right. Since 233121 and 233122 are forbidden it must extend as $233^*123$. Because of the candidate word given above we must continue as $233^*1231$. Since 2312 and 131 are forbidden it must continue as $233^*123111$. Since 1231113 and 231112 are forbidden it must continue as $233^*12311111$. Since 131 is forbidden it should continue as $233^*1231111132$. Since 111322, 1323, 1123111, 3111113212 are forbidden this is forced to $233^*123111113213$. Hence it should continue as $233^*12311111321331$ because of the candidate above. Since 313 and 2133121 are forbidden it should continue as $233^*12311111321331223$. If there is 1 to the right, since 313, 2132, 1223112, 223111, 131, 312231132, 312231133 are forbidden then the word would have no extensions (in other words 312231 is forbidden). Thus it must continue as $233^*123111113213312232$, so it is minimized with $\overline{32}$ because 1323 is forbidden. 

\begin{remark}
The minimum of 12133 is the same as the palindrome 12133121 which is \[m(\overline{333111}\beta^T3^*3\beta\overline{111333})=3.98827021\dots\] where \[\beta=121132233123112221331132212212113223113221331122311331222.\]

The middle of this word is palindromic and the period is semi--symmetric. 
\end{remark}

\bigskip

\noindent{\Large{\textbf{Intervals near to and containing Freiman's constant $\boldsymbol{c_{F}}$}}}
\vspace*{5pt}

In this setting, our subshifts will be subshifts of $\{1,2,3,4\}^{\Z}$. Recall that Freiman's constant $c_{F} = \lambda_{0}(\overline{121313}22344^{*}3211\overline{313121}) = 4.52782956616\ldots$ is the beginning of Hall's ray $[c_{F},\infty)\subset L\subset M$.

\subsection{Interval [4.520781, 4.523103)}\hfill

Consider the subshift $\Sigma(C)$ defined by the forbidden subwords
\begin{gather*}
    41, 42, 334, 343, 434, 234, 31313, 131312131312, 231312131312, 331312131312
\end{gather*}
Observe that
\begin{align*}
    \lambda_0(4^*1)&>4.75 & \lambda_0(3234^*4323)&>4.5235 \\
    \lambda_0(4^*2)&>4.56 & \lambda_0(33131213^*1312)&>4.5228 \\
    \lambda_0(344^*34)&>4.54 \\
    \lambda_0(334^*43)&>4.533 & \lambda_0(34^*3)&>4.5224 \\
    \lambda_0(313^*13)&>4.524
\end{align*}

Observe that if 334 is a subword, since $\lambda_0(334^*3)>4.56$ and 41, 42 are forbidden, then it must extend to 33443 or 33444 which are both forbidden. Similarly, if 434 is a subword then it must extend to 4434, so $\lambda_0(4^*34)$ is at least $\lambda_0(344^*34)>4.54$. 
Is easy to see that if 234 is a subword, then it must extend to $3234^*4323$ (otherwise the value is very large) and we can bound by below the minimum of this subword by $\lambda_0(3234^*4323)>4.5235$.

We choose
\begin{equation*}
    \Sigma(B)=\{\alpha\in\Sigma(C): 12131312  ~\text{is not a substring of}~\alpha\}
\end{equation*}

We claim that $m(\underline{b})\leq m(\overline{313111313^*121})=4.520780004\dots$ for all $\underline{b}\in\Sigma(B)$. Indeed, let $\underline{b}\in\Sigma(B)$ such that $m(\underline{b})=\lambda_0(\underline{b})$ assumes the maximum of $\Sigma(B)$, so we have either $b_0^*$ is equal to $3^*$ or $4^*$. First suppose $b_0^*=4^*$. Since 41, 42, 343 are forbidden, we have to extend to $34^*4$. Since 434, 334, 234 are forbidden we must extend to $\lambda_0(134^*4)<4.52$. Now suppose $b_0^*=3^*$. Since $\lambda_0(3^*4)<4.07$, $\lambda_0(3^*3)<4.15$, $\lambda_0(3^*2)<4.29$, we must continue as $13^*1$. Since $\lambda_0(13^*11)<4.48$, $\lambda_0(213^*12)<4.48$ and 31313, 14 are forbidden we must continue as $313^*12$. Finally since 41, 31313, 12131312 are forbidden we have
\begin{equation*}
    \lambda_0(313^*12)\leq m(\overline{313111313^*121})=4.520780004\dots
\end{equation*}

As in the cases considered above, if we have two continuations $v^{1} = v^{1}_{1}v^{1}_{2}$ and $v^{2} = v^{2}_{1}v^{2}_{2}$ with $|v^{1}_{1}-v^{2}_{1}| > 1$ (e.g., $v^{1} = 1\ldots$ and $v^{2} = 3\ldots$) then we can always find a continuation using a $v_{B} = v_{1}v_{2}\ldots$ with, say, $v^{1}_{1} > v_{1} > v^{2}_{1}$.

Hence, we need only consider the cases $v^{1} = 1\ldots, v^{2} = 2\ldots$, or $v^{1} = 2\ldots, v^{2} = 3\ldots$, or $v^{1} = 3\ldots, v^{2} = 4\ldots$.

In the first case, since 24, 14, 31313 and 231312131312 are forbidden, we have
\[[0;v^{2}] = [0;2,\ldots] \leq [0;2,3,\overline{1,3,1,2,1,3,1,3,1,1,1,3}]\]
which is allowed since $\lambda_0(113^*1)<4.48$ and
\begin{equation*}
    \lambda_0(11313^*1213)\leq m(\overline{313111313^*121})=4.520780004\dots.
\end{equation*}
In the second case, 42 being forbidden forces $a_{N}\neq 4$. If $a_{N}\in\{2,3\}$, then since 234, 334, 41, 31313 and 331312131312 are forbidden, we have
\[[0;v^{2}] = [0;3,\ldots] \leq [0;3,3,\overline{1,3,1,2,1,3,1,3,1,1,1,3}],\]
which is again allowed. Whereas if $a_{N} = 1$, since 41, 42, 343, 434, 433, 432, 31313, and 131312131312 are forbidden, we have
\[[0;v^{2}] = [0;3,\ldots] \leq [0;3,4,4,4,3,1,\overline{1,3,1,3,1,2,1,3,1,3,1,1}\]
which is also allowed. In the final case, since 14 and 24 are forbidden, we have $a_{N} \in \{3,4\}$. Hence, since 41, 42, 434, 433, 432, 31313 and 131312131312 are forbidden, we have
\[[0;v^{2}] = [0;4,\ldots] \leq [0;4,4,3,1,\overline{1,3,1,3,1,2,1,3,1,3,1,1}]\]
which is again allowed.

The minimum Markov value of a bi-infinite sequence containing 331312131312 is
\begin{equation*}
    m(\overline{111313121313}213131213133131213^*1312\overline{313121313111})=4.5231035\dots
\end{equation*}
\begin{remark}
The middle word is palindrome and the period is semi-symmetric.
\end{remark}

The minimum Markov value of a bi-infinite sequence containing 231312131312 is
\begin{equation*}
    m(\overline{3131213^*1312})=4.523119130\dots
\end{equation*}

\begin{remark}
This period is \textbf{not} semi-symmetric. 
\end{remark}

The minimum Markov value of a bi-infinite sequence containing 131312131312 is
\begin{equation*}
    m(\overline{21313121313213}13121313113131213^*1\overline{31231312131312})=4.52314985\dots
\end{equation*}

\begin{remark}
The middle word is palindrome but the period is \textbf{not} semi-symmetric.
\end{remark}

The minimum Markov value of a bi-infinite sequence containing 343 is
\begin{equation*}
    m(\overline{3134^*31})=\sqrt{82}/2=4.52769\dots
\end{equation*}
\begin{remark}
The period is semi-symmetric.
\end{remark}
Indeed, if 41, 313131, 313132, 4313133, 43131344 are forbidden, then a continued fraction in $\{1,2,3,4\}$ that begins with $[4;3]$ is minimized with $[4;\overline{3,1,3,1,3,4}]$.

The minimum Markov value of a bi-infinite sequence containing 31313 is the same as the minimum of 343 which is
\begin{equation*}
    m(\overline{4313^*13})=\sqrt{82}/2=4.52769\dots
\end{equation*}
Indeed, since 313131, 313132, 3313133, 3313134 are forbidden the subword 31313 must extend to 4313134. Now since 41, 42 and 43131344 are forbidden it must extend to 343131343. Moreover, since 3434, 3433, 3432, 34311, 34312, 14, 343134, 343133, 343132, 3431311 431313431312 are forbidden this subword extends to 31313431313431313. In other words, we have self-replication so in fact 31313 extends to $\overline{431313}$.

Similarly, if 343 is a subword, then it must extend to 343131343131343, where we used the forbidden words of the previous paragraph plus 21313431312, 31313431312.

\subsection{Interval [4.5251, 4.5279)}\hfill

We realise the interval [4.5251, 4.5279) as a union of two good intervals; namely, [4.5251, 4.52769) and [4.52753, 4.5279).

\subsubsection{Interval [4.5251, 4.52769)}\hfill

Consider the subshift
\begin{gather*}
    \Sigma(C) = \{\underline{a} \in\{1,2,3,4\}^{\Z}: 41, 42, 334, 31313~\text{and their transposes are not substrings of}~\underline{a}\}.
\end{gather*}

We choose
\begin{equation*}
    \Sigma(B)=\{\underline{a}\in\Sigma(C): 343, 234, 434  ~\text{are not substrings of}~\underline{a}\}
\end{equation*}

We claim that $m(\underline{b})\leq m(\overline{13^*1312})=4.52509\dots$ for all $\underline{b}\in\Sigma(B)$. Indeed, let $\underline{b}\in\Sigma(B)$ such that $m(\underline{b})=\lambda_0(\underline{b})$ assumes the maximum of $\Sigma(B)$, so we have either $b_0^*$ is equal to $3^*$ or $4^*$. First suppose $b_0^*=4^*$. Since 41, 42, 343 are forbidden, we have to extend to $34^*4$. Since 434, 334, 234 are forbidden we must extend to $\lambda_0(134^*4)<4.52$. Now suppose $b_0^*=3^*$. Since $\lambda_0(3^*4)<4.07$, $\lambda_0(3^*3)<4.15$, $\lambda_0(3^*2)<4.29$, we must continue as $13^*1$. Now since 31313 and 14 are forbidden, we see that $\lambda_0(13^*1)\leq m(\overline{13^*1312})=4.52509\dots$.

As above, we need only consider the cases $v^{1} = 1\ldots, v^{2} = 2\ldots$, or $v^{1} = 2\ldots, v^{2} = 3\ldots$, or $v^{1} = 3\ldots, v^{2} = 4\ldots$.

In the first case, since 41, 14 and 31313 are forbidden, we have
\[[0;v^{1}] \geq [0;1,1,\overline{3,1,3,1,2,1}]\]
so, since $\lambda_0(13^*1)\leq m(\overline{13^*1312})=4.52509\dots$, we can use the continuation $v_{B} = 11\overline{313121}$.

In the second case, as similar argument shows that
\[[0;v^{1}]\geq [0;\overline{2,1,3,1,3,1}]\]
so that $v_{B} = \overline{213131}$ can be used.

Finally, when $v^{1} = 3\ldots$ and $v^{2} = 4\ldots$, since 14 and 24 are forbidden, we must have $a_{N}\in\{3,4\}$ and so we can show that
\[[0;v^{1}]\geq [0;\overline{3,1,3,1,2,1}]\]
so that the continuation $v_{B} = \overline{313121}$ can be used.

\subsubsection{Interval [4.52753, 4.5279)}\hfill

Now consider the subshift $\Sigma(C)$ defined by the forbidden subwords
\begin{gather*}
    41, 42, 334, 313131, 313132, 313133, 4434, 2343, 444321, 444322   
\end{gather*}
and their transposes.

Observe that
\begin{align*}
    \lambda_0(234^*3)&>4.55 \\
    \lambda_0(313^*131)&>4.542 \\
    \lambda_0(344^*34)&>4.54 \\
    \lambda_0(313^*132)&>4.532 \\
    \lambda_0(444^*321)&>4.53 \\
    \lambda_0(34313^*133)&>4.529 \\
    \lambda_0(444^*322)&>4.5279
\end{align*}

Note that if 313133 is a subword, then since 41 and 42 and forbidden we have that $\lambda_0(313^*133)$ is at least $\lambda_0(34313^*133)>4.529$. Similarly if 4434 is a subword, then $\lambda_0(44^*34)$ is at least $\lambda_0(344^*34)>4.54$. 

Now set
\begin{equation*}
    \Sigma(B)=\{\underline{a}\in\Sigma(C): 343, 23443, 34432, 31313 ~\text{are not substrings of}~\underline{a}\}
\end{equation*}

We claim that $m(\underline{b})\leq m(\overline{23444^*3})=4.52752\dots$ for all $\underline{b}\in\Sigma(B)$. Indeed, let $\underline{b}\in\Sigma(B)$ such that $m(\underline{b})=\lambda_0(\underline{b})$ assumes the maximum of $\Sigma(B)$, so we have either $b_0^*$ is equal to $3^*$ or $4^*$. First suppose $b_0^*=3^*$. Since $\lambda_0(3^*4)<4.07$, $\lambda_0(3^*3)<4.15$, $\lambda_0(3^*2)<4.29$, we must continue as $13^*1$. Now since 31313 and 14 are forbidden, we see that $\lambda_0(13^*1)\leq m(\overline{13^*1312})=4.52509\dots$. Now suppose $b_0^*=4^*$. Since $\lambda_0(44^*4)<4.48$ and 41, 42, 343 are forbidden, we have to extend to $34^*4$. Since 4434, 334 are forbidden and $\lambda_0(134^*4)<4.52$, we must extend to $234^*4$. Since 23443 is forbidden it extends to $234^*44$. Finally using that 444321, 444322, 24, 41, 42, 2343, 4434, 433 are forbidden this is maximized with $\lambda_0(234^*44)\leq m(\overline{234^*443})=4.5275206\dots$

The arguments of the previous interval will allow us to only have to consider the continuations $v^{1} = 3\ldots, \,v^{2} = 4\ldots$, with $a_{N}\in\{3,4\}$. If $a_{N} = 3$, then, since 331313 is now forbidden, we will have
\[[0;v^{1}]\geq [0;\overline{3,1,3,1,2,1}]\]
so the continuation $v_{B} = \overline{313121}$ can be used. Otherwise, since 41, 42, 4434, 433, 444321, 444322 and 2343 are forbidden, we have
\[[0;v^{2}]\leq [0;\overline{4,4,3,2,3,4}]\]
so that, since
\[\lambda_0(444^*\overline{323444})\leq m(\overline{444^*323})=4.5275206\dots\]
we can use the continuation $v_{B} = \overline{443234}$.

\begin{remark}
Even though this fact is not used here, we note that Freiman's gap $(\nu_{F},c_F)$, satisfies
\begin{equation*}
    \nu_{F} = \lambda_{0}(\overline{323444}313134^{*}313121133\overline{313121})=\lim_{n\to\infty} \lambda_0(\overline{(323444)^n313134^{*}313121133(313121)^n1})\in L^\prime.
\end{equation*}

Indeed, one only needs to use that $m(\overline{313121}1\overline{323444})=m(\overline{323444^*})=4.5275206\dots<\nu_{F}$.
\end{remark}

\begin{remark}
As mentioned at the beginning of this section, we wish to point out again that all but finitely many of the currently known regions of $M\setminus L$ occur inside good intervals. It is an open question as to whether or not there exist elements of $M\setminus L$ close to $\nu$. The fact that one can place $\nu$ within a good interval is an interesting observation.
\end{remark}

%%%%%%%%%%%%%%%%%%%%%%%%%%%%%%%%%%%%%%%%%%%%

\section{$M\setminus L$ near 3.942}

\subsection{Brief description of the algorithm}\label{subsec:M-Lalg}\hfill

Here we give a description of the algorithm used in the computer investigations used to find new regions of $M\setminus L$. See also~\cite{JMM}.

Let $w$ be a (finite) odd not semi--symmetric word  and $j_0:=m(\overline{w})$. By rewriting the period of this Markov value we can write $w^*=\eta_13^*\eta_2$ where both $\eta_1,\eta_2$ have the same length and where $m(\overline{w})=\lambda_0(\overline{w^*})$. The first step is to establish local uniqueness: we want to find an $\varepsilon>0$ such that if $\lambda_0(b)\in(j_0-\varepsilon,j_0+\varepsilon)$ for some $b\in\{1,2,3\}^\Z$, then we must have $b=\eta_2w^*\eta_1$. The next step is to prove self-replication, that is: there is $\varepsilon>\delta>0$  such that if $\lambda_0(b)\in(j_0-\delta,j_0+\delta)$ for some $b\in\{1,2,3\}^\Z$, then $b=\eta_2^\prime ww^*w\eta_1^\prime$ where either $\eta_1^\prime=\eta_1$ and $\eta_2^\prime$ is a suffix of $\eta_2$ or $\eta_2^\prime=\eta_2$ and $\eta_1^\prime$ is a prefix of $\eta_1$. In particular we must have that either $b=\overline{w}w^*w\eta_1^\prime$ or $b=\eta_2^\prime ww^*w\overline{w}$. To prove this it is \textbf{essential} to just use forbidden words: that is any words of the form $\tau_2\eta_2w\eta_1\tau_1$, where either $\tau_2\eta_2$ is a suffix of $ww$ and $\eta_1\tau_1$ minus the last digit is a prefix of $ww$ but $\eta_1\tau_1$ is not, or $\eta_1\tau_1$ is a prefix of $ww$ and $\tau_2\eta_2$ minus the first digit is a suffix of $ww$ but $\tau_2\eta_2$ is not. In particular $\tau_2\eta_2w\eta_1\tau_1$ is not a subword of $wwwww$. This will give us a finite sequence of forbidden words $f_1,f_2,\dots$ such that $\lambda_0(f_i)>j_0+\varepsilon_i$ with $\varepsilon_i<\varepsilon$. It turns out that there should exist a minimal forbidden word $f_{\textrm{min}}$ that has the lowest lower bound $\varepsilon_{\textrm{min}}=\delta$. Next we consider $t_{\textrm{min}}$, which is the minimum Markov value of a bi-infinite sequence that contains $f_{\textrm{min}}$. If we have luck, we will have that $t_{\textrm{min}}\in L$. Otherwise we minimize the next minimal forbidden word until this minimum belong to $L$. We call $j_1=j_1(w)\in L$ to be the right border of the maximal gap $(j_0,j_1)\cap L=\varnothing$.

\subsection{Local uniqueness}\label{sec:local_uniq}\hfill

We denote throughout this section $w=12111233311133232$ and $w^*=121112333^*11133232$. Note that we then have $j_0=m(\overline{w^*}) = 3.942001159911341469213548\dots$.

\begin{lemma}[Small forbidden words]\label{lem:forwords1}
Let $b\in\{1,2,3\}^\Z$. 
\begin{enumerate}

\item If $b=13^*1$, $3^*13$, $3^*12$, $323^*1$, $223^*1$, $23^*111$, $1123^*112$, $2123^*1123$, $2123^*1122$, $2123^*11211$, $2123^*11212$, $31123^*1$, $21123^*1$, then $\lambda_0(b)>j_0+10^{-4}$.
\item $b=1333^*1113$, $32333^*11133$, $22333^*11133$, $3331113^*33211$ then $\lambda_0(b)>j_0+10^{-5}$.
\item If $b=1123^*1133$, $1123^*11321$ then $\lambda_0(b)>j_0+10^{-4}$.
\end{enumerate}
\end{lemma}

\begin{corollary}\label{cor:simplified_forbidden1}
If $m(b)<j_0+10^{-4}$, then $b$ does not contain  23112, 11231.    
\end{corollary}

\begin{lemma}
Let $b\in\{1,2,3\}^\Z$ be such that $\lambda_0(b)<j_0+10^{-6}$. Then $b$ or $b^T$ must be of the form:
\begin{enumerate}
\item $b=1^*,2^*$ and $\lambda_0(b)<j_0-10^{-1}$.
\item\label{localuniq_item:2} $b=33^*3,33^*2,23^*2$ and $\lambda_0(b)<j_0-10^{-2}$.
\item $b=2123^*113$ and $\lambda_0(b)<j_0-10^{-3}$.
\item\label{localuniq_item:4} $b=133^*11,233^*11,333^*112,333^*113,333^*1112,333^*1111$ and $\lambda_0(b)<j_0-10^{-3}$. 
\item $b=3333^*11133$ and $\lambda_0(b)<j_0-10^{-4}$.
\item $b=12333^*11133$.
\end{enumerate}
\end{lemma}

\begin{proof}
Let us assume that $b=3^*$. If $b$ is not of the forms in item \ref{localuniq_item:2}, without loss of generality we can assume that $b=3^*1$. Since 313 and 312 are forbidden we must extend this word to $b=3^*11$. Since 131 is forbidden we have that $b$ is either $b=23^*11$ or $b=33^*11$. 

If $b=23^*11$, then since 23111, 23112, 3231, 2231, 312 and 11231 are forbidden by \Cref{cor:simplified_forbidden1} and \Cref{lem:forwords1}, it is forced to $b=2123^*113$ and $\lambda_0(b)<j_0-10^{-3}$.

Assume $b=33^*11$. If $b$ is not of the forms in item \ref{localuniq_item:4}, then $b=333^*1113$. Since 131 and 11132 are forbidden we have $b=333^*11133$. Since 13331113 is forbidden we can assume that $b=2333^*11133$ or $3333^*11133$. Since 323331113 and 223331113 are forbidden, we must extend in the former case to $b=12333^*11133$.
\end{proof}

By the previous lemma and since 312 is forbidden, it suffices to analyse the extensions to the left of $12333^*11133$, i.e. $212333^*11133$, $112333^*11133$.

\subsubsection{Extensions of the word $212333^*11133$}

\begin{lemma}\hfill
\begin{enumerate}
\item If $b=212333^*111331,212333^*111332$ then $\lambda_0(b)<j_0-10^{-5}$.
\item If $b=3212333^*111333$ then $\lambda_0(b)>j_0+10^{-6}$.
\item If $b=1212333^*1113332, 1212333^*1113333$ then $\lambda_0(b)<j_0-10^{-6}$.
\item If $b=2212333^*111333212$ then $\lambda_0(b)>j_0+10^{-7}$.
\end{enumerate}
\end{lemma}

\begin{corollary}\label{cor:fw1}
If $m(b)<j_0+10^{-7}$, then $b$ does not contain 22123331113332. 
\end{corollary}
\begin{proof}
Use the forbidden words 311133322, 311133323, 333111333211, 312 and the previous lemma.
\end{proof}

Since 31113331 is forbidden,  by the previous lemma it suffices to analyse the extensions of $2212333^*1113333$.

\begin{lemma}
If $b=22212333^*1113333$ then $\lambda_0(b)<j_0-10^{-8}$.
\end{lemma}
\begin{proof}
\begin{equation*}
    \lambda_0(22212333^*1113333)\leq\lambda_0(\overline{13}22212333^*1113333311\overline{13})<j_0-10^{-8}.
\end{equation*}
\end{proof}

\begin{lemma}\hfill
\begin{enumerate}
\item If $b=32212333^*1113333$ then $\lambda_0(b)<j_0-10^{-6}$.
\item If $b=12212333^*11133333,12212333^*11133332$ then $\lambda_0(b)>j_0+10^{-7}$.
\item If $b=212212333^*111333$ then $\lambda_0(b)>j_0+10^{-6}$.
\item If $b=1112212333^*111333311, 2112212333^*111333311$ then $\lambda_0(b)>j_0+10^{-7}$.
\end{enumerate}    
\end{lemma}

\begin{corollary}
If $m(b)<j_0+10^{-8}$, then $b$ does not contain 12212333111333.    
\end{corollary}
\begin{proof}
If $b$ contains the subword 12212333111333, then using the previous lemma, \Cref{cor:fw1} and the forbidden words 31113331, 312, 313, 23112, 131 this subword must extend to the word 33112212333111333311. Now using that 131 and 11132 are forbidden, we have the inequality
\begin{equation*}
    \lambda_0(33112212333^*111333311)\geq \lambda_0(\overline{13}33112212333^*111333311133\overline{31})>j_0+10^{-8}.
\end{equation*}
\end{proof}

\begin{corollary}
If $b\in\{1,2,3\}^\Z$ is such that $b=\ldots 212333^*11133\ldots$, then $\abs{m(b)-j_0}>10^{-8}$.
\end{corollary}

\subsubsection{Local uniqueness up to $w^*$}

Now we continue analyzing the extensions of the candidate $112333^*11133$, which is the one that
will converge to the non semi-symmetric word we are looking for (the other branches are already discarded).

\begin{lemma}\hfill
\begin{enumerate}
\item If $b=112333^*111333$ then $\lambda_0(b)>j_0+10^{-5}$.
\item If $b=2112333^*1113311, 1112333^*1113311$ then $\lambda_0(b)<j_0-10^{-6}$.
\item If $b=3112333^*11133113$ then $\lambda_0(b)>j_0+10^{-6}$.
\item If $b=33112333^*11133112, 33112333^*11133111$ then $\lambda_0(b)<j_0-10^{-7}$.
\end{enumerate}    
\end{lemma}

By the previous lemma and since 312, 313, 23112 are forbidden, it suffices to analyse the extensions of $112333^*111332$.

\begin{lemma}\hfill
\begin{enumerate}
\item If $b=3112333^*111332,2112333^*111332$ then $\lambda_0(b)>j_0+10^{-5}$.
\item If $b=112333^*1113321$ then $\lambda_0(b)>j_0+10^{-5}$.
\item If $b=1112333^*1113322$ then $m(b)>j_0+10^{-7}$.
\end{enumerate}
\end{lemma}
\begin{proof}
If $b=1112333^*1113322$ then $\lambda_0(b)\geq\lambda_0(\overline{13}331112333^*1113322\overline{13})>j_0+10^{-7}$ where we used that 131, 11132 are forbidden.
\end{proof}

By the previous lemma, it suffices to analyse the extensions of $1112333^*1113323$.

\begin{lemma}\hfill
\begin{enumerate}
\item If $b=11112333^*111332$ then $\lambda_0(b)>j_0+10^{-5}$.
\item If $b=31112333^*1113323$ then $\lambda_0(b)<j_0-10^{-6}$.
\item If $b=321112333^*111332, 221112333^*111332$ then $m(b)>j_0+10^{-6}$.
\item If $b=1121112333^*11133233,2121112333^*11133233$ then $\lambda_0(b)>j_0+10^{-7}$.
\end{enumerate}
\end{lemma}
\begin{proof}
If $b=221112333^*111332$ then $\lambda(b)\geq\lambda_0(\overline{13}221112333^*11133232\overline{31})>j_0+10^{-6}$ where we used that 3231 is forbidden.   
\end{proof}

\begin{corollary}
If $m(b)<j_0+10^{-7}$ then $b$ does not contain 2111233311133233.    
\end{corollary}

Since 3231 is forbidden we obtain
\begin{corollary}\label{cor:localuniq_w}
Let $b\in\{1,2,3\}^\Z$. If $\abs{\lambda_0(b)-j_0}<10^{-8}$ then $b=\ldots w^*\ldots$. 
\end{corollary}

Since 312 is forbidden, it suffices to analyze the extensions of $1w^*$
and $2w^*$.

\begin{lemma}\hfill
\begin{enumerate}
\item If $b=1w^*1,1w^*2,31w^*3,21w^*3$ then $\lambda_0(b)>j_0+10^{-8}$.
\item If $b=111w^*3,211w^*3$ then $m(b)>j_0+10^{-8}$.
\item If $b=311w^*32$ then $\lambda_0(b)<j_0-10^{-9}$.
\item If $b=311w^*33$ then $m(b)>j_0+10^{-9}$.
\end{enumerate}    
\end{lemma}

\begin{proof}
If $b=211w^*3$, then since 1323 is forbidden $\lambda_0(b)\geq\lambda_0(211w^*32)>j_0+10^{-8}$. 

If $b=3311w^*33$, then since 23111 and 131 are forbidden, we will have $\lambda_0(b)\geq\lambda_0(3311w^*33)>j_0+10^{-9}$.
\end{proof}

\begin{corollary}
If $m(b)<j_0+10^{-9}$ then $b$ does not contain $1w33$.
\end{corollary}

\begin{corollary}
Let $b\in\{1,2,3\}^\Z$. If $b=\ldots 1w^*\ldots$ then $\abs{\lambda_0(b)-j_0}>10^{-9}$.
\end{corollary}

\begin{lemma}\hfill
\begin{enumerate}
\item If $b=2w^*2,2w^*3,2w^*11,12w^*1,22w^*12$ then $\lambda_0(b)<j_0-10^{-8}$.
\item If $b=32w^*123,32w^*122$ then $\lambda_0(b)<j_0-10^{-9}$.
\item If $b=132w^*121$ then $\lambda_0(b)>j_0+10^{-9}$.
\item If $b=332w^*1211,332w^*1212$ then $\lambda_0(b)<j_0-10^{-9}$.
\item If $b=232w^*1212$ then $\lambda_0(b)>j_0+10^{-10}$.
\item If $b=232w^*12112,232w^*12113,1232w^*12111$ then $\lambda_0(b)<j_0-10^{-9}$.
\item If $b=2232w^*12111$ then $\lambda_0(b)<j_0-10^{-11}$.
\item If $b=3232w^*121111$ then $\lambda_0(b)<j_0-10^{-11}$.
\item If $b=3232w^*121113,3232w^*1211121,3232w^*1211122,23232w^*121112$ then $\lambda_0(b)>j_0+10^{-11}$.
\end{enumerate}    
\end{lemma}
\begin{proof}
If $b=32w^*122$ then $\lambda_0(b)\leq\lambda_0(1132w^*122)<j_0-10^{-9}$ where we used that 312 and 313 are forbidden.   

If $b=232w^*1212$ then $\lambda_0(b)\geq\lambda_0(21232w^*1212)>j_0+10^{-10}$ where we used that 312 is forbidden.

If $b=2232w^*12111$ then $\lambda_0(b)\leq\lambda_0(\overline{12}232w^*1211133\overline{31})<j_0-10^{-11}$ where we used that 312, 131, 11132 are forbidden.
\end{proof}

By the previous lemma, it suffices to analyse the extensions of $33232w1211123$.

\begin{lemma}\hfill
\begin{enumerate}
\item If $b=333232w^*1211123$ then $m(b)>j_0+10^{-12}$.
\item If $b=233232w^*12111232,233232w^*12111233$ then $\lambda_0(b)>j_0+10^{-12}$.
\end{enumerate}
\end{lemma}

\begin{proof}
If $b=333232w^*1211123$ then $\lambda_0(b)\geq\lambda_0(333232w^*121112311)>j_0+10^{-12}$ where we used that 313, 312 are forbidden. 
\end{proof}

\begin{corollary}
If $m(b)<j_0+10^{-12}$, then $b$ does not contain 233232w1211123. 
\end{corollary}
\begin{proof}
Use the forbidden word 11231 and the previous lemma.
\end{proof}

\begin{corollary}\label{cor:localuniq_partial}
Let $b\in\{1,2,3\}^\Z$ such that $m(b)=\lambda_0(b)$. If $\abs{m(b)-j_0}<10^{-12}$ then either $b=\ldots 1133232w^*12111232\ldots$ or $b=\ldots 1133232w^*12111233\ldots$.    
\end{corollary}
\begin{proof}
Use the previous lemma and that 313, 312, 11231 are forbidden.    
\end{proof}

\subsubsection{Extensions of the word $\tau=1133232w12111232$}
Denote $\tau^*=1133232w^*12111232$.

\begin{lemma}\label{lem:first_branch_tau}\hfill
\begin{enumerate}
\item If $b=1\tau^*,\tau^*2,\tau^*3$ then $\lambda_0(b)<j_0-10^{-12}$.
\item If $b=3\tau^*12$ then $\lambda_0(b)>j_0+10^{-12}$.
\item If $b=23\tau^*112,23\tau^*113,33\tau^*112,33\tau^*113,33\tau^*1111$ then $\lambda_0(b)<j_0-10^{-14}$.
\item If $b=23\tau^*1112,23\tau^*1113,33\tau^*1113,33\tau^*1112$ then $\lambda_0(b)>j_0+10^{-15}$.
\item If $b=2123\tau^*11112,23\tau^*11113,$ then $\lambda_0(b)<j_0-10^{-14}$.
\item If $b=2123\tau^*11111$ then $\lambda_0(b)>j_0+10^{-14}$.

\end{enumerate}
\end{lemma}
\begin{proof}
If $b=33\tau^*1112$ then $\lambda(b)\geq\lambda_0(\overline{13}3\tau^*111\overline{23})>j_0+10^{-15}$ where we used that 11231 and 3231 are forbidden.    
\end{proof}

Using that 131 is forbidden we get
\begin{corollary}
If $m(b)<j_0+10^{-15}$ then $b$ does not contain $3\tau 1112$ neither $2\tau 1113$. 
\end{corollary}

\begin{corollary}
If $b=3\tau^*$ then $\abs{m(b)-j_0}>10^{-15}$.
\end{corollary}
\begin{proof}
Use the forbidden words 3231, 2231, 11231, 312 we see that $b=23\tau^*$ must extend to $2123\tau^*$ and then we use the previous lemma. The case $b=33\tau^*$ follows from items 1 to 4 of the previous lemma.
\end{proof}

By the previous lemma, it suffices to analyse the extensions of $2\tau 1$.

\begin{lemma}\hfill
\begin{enumerate}
\item If $b=2\tau^*11,22\tau^*12,32\tau^*12,12\tau^*123,212\tau^*122,112\tau^*122,112\tau^*1211$ then $\lambda_0(b)<j_0-10^{-14}$.
\item If $b=1212\tau^*1211,2212\tau^*1211,3212\tau^*12112,3212\tau^*12111,23212\tau^*121133$, $33212\tau^*121133$, $2113212\tau^*121133$, $3113212\tau^*121133$ then $\lambda_0(b)>j_0+10^{-16}$.
\end{enumerate}    
\end{lemma}

\begin{corollary}
If $m(b)<j_0+10^{-16}$ then $b$ does not contain $212\tau^*1211$. 
\end{corollary}

\begin{corollary}
If $b=2\tau^*1211$ then $\abs{m(b)-j_0}>j_0+10^{-16}$.
\end{corollary}

By the previous lemma and since 312 is forbidden, it suffices to analyse the extensions of $12\tau 1212$.

\begin{lemma}\hfill
\begin{enumerate}
\item If $b=212\tau^*1212$ then $\lambda_0(b)>j_0+10^{-13}$.
\item If $b=2112\tau^*1212,3112\tau^*1212,1112\tau^*12123,1112\tau^*12122,11112\tau^*12121$ then \newline $\lambda_0(b)<j_0-10^{-15}$.
\item If $b=31112\tau^*12121$ then $\lambda_0(b)>j_0+10^{-15}$.
\item If $b=221112\tau^*121211,221112\tau^*121212$ then $\lambda_0(b)<j_0-10^{-15}$.
\end{enumerate}    
\end{lemma}

\begin{proof}
If $b=1112\tau^*12122$ then $\lambda_0(b)\leq\lambda_0(\overline{13}331112\tau^*12122\overline{31})<j_0-10^{-15}$ where we used the forbidden words 131, 23111.   
\end{proof}

By the previous lemma, it suffices to analyse the extensions of $121112\tau 12121$. 

\begin{lemma}\label{lem:last_branch_tau}\hfill
\begin{enumerate}
\item If $b=1121112\tau^* 121211$ then $\lambda_0(b)<j_0-10^{-15}$.
\item If $b=2121112\tau^* 121212,1121112\tau^* 1212121,1121112\tau^* 1212122$ then $m(b)>j_0+10^{-17}$.
\item If $b=31121112\tau^* 1212123$ then $\lambda_0(b)<j_0-10^{-16}$.
\item If $b=11121112\tau^* 121212,321121112\tau^* 121212,221121112\tau^* 121212$ then $\lambda_0(b)>j_0+10^{-16}$.
\item If $b=1121121112\tau^* 121212,22121121112\tau^* 121212,32121121112\tau^* 121212$,\newline $112121121112\tau^* 121212$ then $m(b)>j_0+10^{-19}$.
\item If $b=21121112\tau^* 12121232,21121112\tau^* 12121233$ then $\lambda_0(b)>j_0+10^{-17}$.
\item If $b=212121121112\tau^* 12121231133,212121121112\tau^* 12121231132$ then $\lambda_0(b)<j_0-10^{-19}$.
\item If $b=32121112\tau^*121211,22121112\tau^*121211,12121112\tau^*1212112,12121112\tau^*1212113,\newline 12121112\tau^*12121111,12121112\tau^*12121112$ then $\lambda_0(b)<j_0-10^{-17}$.
\item If $b=212121112\tau^*12121113$ then $\lambda_0(b)>j_0+10^{-17}$.
\item If $b=112121112\tau^*121211133$ then $\lambda_0(b)<j_0-10^{-18}$.
\end{enumerate}    
\end{lemma}
\begin{proof}
If $b=1121112\tau^* 1212122$ then $\lambda(b)\geq\lambda_0(\overline{13}331121112\tau^* 1212122\overline{13})>j_0+10^{-16}$ where we used that 131, 23112 are forbidden.
If $b=1121121112\tau^* 121212$ then $$\lambda(b)\geq\lambda_0(\overline{31}1121121112\tau^* 121212311\overline{13})>j_0+10^{-17}$$ where we used that 312 and 313 are forbidden.

If $b=22121121112\tau^* 121212,32121121112\tau^* 121212,112121121112\tau^* 121212$ then 
\begin{equation*}
    \lambda_0(b)\geq\lambda_0(\overline{31}112121121112\tau^* 1212123113\overline{31})>j_0+10^{-19}
\end{equation*}
where we used that 313, 312, 23111, 23112 are forbidden.
\end{proof}

\begin{corollary}\label{cor:localuniq_tau}
If $b=\tau^*$ then $\abs{m(b)-j_0}>10^{-19}$.
\end{corollary}
\begin{proof}
This is consequence of Lemmas and Corollaries \ref{lem:first_branch_tau} through \ref{lem:last_branch_tau}.
\end{proof}

Now we collect all of the forbidden words and the self-replicating word $w_r:=1133232w12111233$ from \Cref{cor:localuniq_partial}. Define $F$ as the following set of words and their transposes:
\begin{itemize}
    \item 131, 313, 312, 3231, 2231, 23111, 23112, 11231.
    \item 13331113, 323331113, 223331113, 112333111333.  

    \item 3212333111333, 22123331113332, 12212333111333, 112333111333, 311233311133113,\linebreak 3112333111332, 2112333111332, 1123331113321, 1123331113322, 11112333111332,\linebreak   321112333111332, 221112333111332, 2111233311133233.
    
    \item $1w1$, $1w2$, $31w$, $21w$, $111w, 211w, 1w33, 132w121, 232w1212, 3232w121113, 3232w1211121$, $3232w1211122$, $23232w121112$, $333232w1211123$, $233232w1211123$.
    
    \item $3\tau12, 3\tau1112, 3\tau1113, 2\tau1113,  2123\tau11111, 212\tau121, 31112\tau12121, 2121112\tau121212,\linebreak 1121112\tau 1212121,
    1121112\tau1212122, 11121112\tau121212,321121112\tau121212,\linebreak 221121112\tau121212,1121121112\tau121212,22121121112\tau121212, 32121121112\tau121212,\linebreak 112121121112\tau121212,21121112\tau12121232,  21121112\tau12121233,212121112\tau12121113$.

    \item $w_r=1133232w12111233$.
\end{itemize}

All the results of this section imply the following:

\begin{corollary}\label{cor:not_forbid_words}
If $m(b)<j_0+10^{-19}$ then $b$ does not contain words from $F\setminus\{w_r,w_r^T\}$.
\end{corollary}

\begin{corollary}[Local uniqueness]\label{cor:fulllocal_uniq}
Let $b\in\{1,2,3\}^\Z$ be such that $j_0-10^{-19}<\lambda_0(b)$. If $b$ does not contain words from $F\setminus\{w_r,w_r^T\}$ then up to transposition we have
\begin{equation*}
    b_{-15}\dots b_0^*\dots b_{16}=1133232w^*12111233.
\end{equation*}
\end{corollary}

\subsection{Self-replication}
Denote $w_r=1133232w12111233$ and $w_r^*=1133232w^*12111233$. By the previous section we know that if $m(b)<j_0+10^{-19}$ and $j_0\leq\lambda_0(b)$ then $b=\ldots w_r^*\ldots$.

\begin{lemma}[Self replication forbidden words]\label{lem:forwords3}\hfill
\begin{enumerate}
\item If $b=2w_r^*,3w_r^*,1w_r^*1, 1w_r^*2, 11w_r^*3, 21w_r^*3, 331w_r^*33, 331w_r^*32$ then $\lambda_0(b)>j_0+10^{-13}$.
\item If $b=331w_r^*3113, 331w_r^*3112, 1331w_r^*3111, 2331w_r^*3111, 3331w_r^*31111, 3331w_r^*31112$ then\linebreak $\lambda_0(b)>j_0+10^{-16}$.
\item If $b=33331w_r^*311133$ then $\lambda_0(b)>j_0+10^{-16}$.
\item If $b=2123331w_r^*3111331, 11123331w_r^*3111331, 2123331w_r^*3111332, 311123331w_r^*31113323$ then $\lambda_0(b)>j_0+10^{-18}$.
\item If $b=2ww^*w3, 2ww^*w2, 12ww^*w1, 22ww^*w1, 32ww^*w11$ then $m(b)>j_0+5\cdot10^{-20}$.
\item If $b=32ww^*w123, 32ww^*w122, 332ww^*w1211, 232ww^*w12113, 232ww^*w12112$ then\linebreak $\lambda_0(b)>j_0+2\cdot 10^{-21}$. 
\item If $b=1232ww^*w12111$ then $\lambda_0(b)>j_0+2\cdot 10^{-21}$.
\item If $b=2232ww^*w121111$ then $\lambda_0(b)>j_0+10^{-21}$.
\item\label{item:min_word} If $b=2232ww^*w121112$ then $\lambda_0(b)>j_0+2\cdot 10^{-22}$.
\end{enumerate}
\end{lemma}
\begin{proof}
If $b=22ww^*w1$ then $\lambda(b)\geq\lambda_0(\overline{21}22ww^*w\overline{12})>j_0+5\cdot 10^{-20}$ where we used that 213 is forbidden.     

If $b=32ww^*w122$ then $\lambda_0(b)\geq\lambda_0(\overline{13}31132ww^*w122\overline{31})>j_0+7\cdot 10^{-21}$ where we used that 213, 313, 11132, 21132 are forbidden. 
\end{proof}

\begin{remark}
The lower bound is decreasing on each item. It is essential for self-replication that there exists a minimal forbidden word. In this case, this is given by the word at item \ref{item:min_word}. 
\end{remark}

We need some auxiliary forbidden words that will help us in the characterisation of $M\setminus L$ but are not needed in the description of the Cantor set.
\begin{lemma}\label{lem:aux_forbidwords}\hfill
\begin{enumerate}
\item If $b=232ww^*w1211113, 33232ww^*w1211112, 33232ww^*w12111111$ then $\lambda_0(b)>j_0+10^{-22}$.
\item If $b=133232ww^*w121111$ then $\lambda_0(b)>j_0+10^{-22}$.
\item If $b=133232ww^*w12111231$ then $\lambda_0(b)>j_0+10^{-23}$.
\item If $b=11133232ww^*w12111232$ then $\lambda_0(b)>j_0+10^{-24}$.
\end{enumerate}    
\end{lemma}

\begin{lemma}[Self-replication]\label{lem:selfreplication}
Let $b\in\{1,2,3\}^\Z$ be such that $m(b)<j_0+10^{-19}$ and suppose that for some index $i\in\Z$ we have $b_{i-15}\dots b_i^*\dots b_{i+16}=1133232w^*12111233=w_r^*$. Then one has
\begin{enumerate}
\item If $\lambda_0(b)<j_0+2\cdot 10^{-21}$ then  
\begin{equation*}
    b_{i-32}\dots b_i^*\dots b_{i+30}=232ww^*w12111.
\end{equation*}
\item If $\lambda_0(b)<j_0+10^{-21}$ and does not contain $2232www121112$ then  
\begin{equation*}
    b_{i-32}\dots b_i^*\dots b_{i+30}=1133232ww^*w12111.
\end{equation*}
In particular
\begin{equation*}
    \dots b_i^*\dots b_{i+30}=\overline{w}w^*w12111.
\end{equation*}
\item If $\lambda_0(b)<j_0+10^{-23}$ and does not contain $11133232www12111232$ then  
\begin{equation*}
    b_{i-32}\dots b_i^*\dots b_{i+33}=1133232ww^*w12111233.
\end{equation*}
In this particular case $b=\overline{w}$.
\end{enumerate}
\end{lemma}

\begin{proof}
By the \Cref{lem:forwords3}, this subword must extend to $b=31w_r^*3$. Since 23111 and 131 are forbidden, we are forced to $b=331w_r^*3$. By \Cref{lem:forwords3} and since 313, 312 are forbidden, it must extend to $b=3331w_r^*31113$. Using that 11132 and 131 are forbidden it is forced to $b=3331w_r^*311133$. By \Cref{lem:forwords3} and since 13331113, 323331113, 223331113 are forbidden, it must extend to $b=123331w_r^*311133$. By \Cref{lem:forwords3} and since 112333111333, 3112333111332, 2112333111332 are forbidden, it must extend to $b=11123331w_r^*3111332$. Using the forbidden words 1123331113321, 11123331113322, 12111233311133233, 3231, 11112333111332, 321112333111332, 221112333111332 and \Cref{lem:forwords3} it must extend to $b=ww^*w$. Now using that $1w1$, 312 are forbidden and \Cref{lem:forwords3} it must extend to $b=32ww^*w121$. Now using that $232w1212$, 213, $132w121$, are forbidden and \Cref{lem:forwords3} it must extend to $b=232ww^*w12111$. 

By using the last three forbidden words of \Cref{lem:forwords3} and since $3232w121113$ is forbidden, we must continued as $b=3232ww^*w12111$. Finally using the forbidden words $23232w121112$, 1323, $333232w1211123$, $233232w1211123$, 313, 213 it is forced to $b=1133232ww^*w12111$.

Now to prove the last claim, just use the inequalities of \Cref{lem:aux_forbidwords} and the forbidden words $3232w121113$, $2w_r$ and $3w_r$. 
\end{proof}

\begin{corollary}
The point $j_0=m(\overline{w})$ is an isolated point in $M$.    
\end{corollary}

The reason we stop extending the word in item 2 of \Cref{lem:selfreplication} is because we want to apply an inductive argument, namely we know that we have local uniqueness at $\eta_1w\eta_2$, so to self replicate we need to extend this subword until we reach a word of the form $\eta_1www\eta_2^\prime$ or $\eta_1^\prime www\eta_2$ where $\eta_i^\prime$ is a \textbf{strict} suffix or prefix of $\eta_i$.

\subsection{Description of $M\setminus L$}

Define the constant $j_1$ to be the minimum Markov value of a bi-infinite sequence $b\in\{1,2,3\}^\Z$ containing the word $2232www121112$ and $m_0$ to be the minimum Markov value of a bi-infinite sequence $b\in\{1,2,3\}^\Z$ containing the word $11133232www12111232$.

The reason to define $j_1$ is the following. First note that this word has the smallest lower bound in \Cref{lem:forwords3}. In fact since we have that $\lambda_0(2232ww^*w121112)<j_0+1.8\cdot 10^{-21}$, there should exist a Markov value smaller than all other forbidden words. This minimal Markov value should be the one that is determining the right extreme of the gap of $L$, provided it belongs to $L$.

\begin{lemma}\hfill
\begin{equation*}
    j_1=m(\overline{12}232wwww1211111\overline{23})=\lambda_0(\overline{12}232ww^*ww1211111\overline{23})\in M.
\end{equation*}
\begin{equation*}
    m_0=m(\overline{w}1211123\overline{21})=\lambda_0(\overline{w}w^*w1211123\overline{21})\in M.
\end{equation*}
Numerically we have $j_1\approx j_0+8.326\cdot 10^{-22}$ and $m_0\approx j_0+5.9\cdot 10^{-24}$.
\end{lemma}
\begin{proof}
Let $b\in\{1,2,3\}^\Z$ be such that $b=2232ww^*w121112$. Since $3232w1211121, 3232w1211122$ are forbidden we must extend to $b=2232ww^*w1211123$. Since 312 is forbidden, we have $\lambda_0(b)\geq\lambda_0(\overline{12}232ww^*w1211123)$. Hence at each step we have to minimize the right with respect to this position avoiding forbidden words. So we have to continue as $b=2232ww^*w12111233$. In particular, because of \Cref{lem:selfreplication} item (i) (note we can assume $m(b)<j_0+10^{-21}$ because of the above candidate) we have self-replication so the word extends to $b=2232ww^*ww12111$. Now we have to be careful because of the dangerous positions inside of each copy of $w$. Since we are minimizing we should continue as $b=2232ww^*ww121111$. Note that $\lambda_0(3232w^*121111)<j_0-10^{-11}$ and $\lambda_0(2232w^*12111)<j_0-10^{-11}$ so the only dangerous position remaining is inside the third copy of $w$. Now since $33232www1211112, 232www1211113$ are forbidden by \Cref{lem:aux_forbidwords}, we must extend to $b=2232ww^*ww1211111$. Since $33232www12111111$ is forbidden by \Cref{lem:aux_forbidwords}, we should extend as $b=2232ww^*ww12111112$, which is good because now
\begin{equation*}
    \lambda_0(www^*w12111112)<j_0+8.31\cdot 10^{-22}<\lambda_0(21212232ww^*ww12111112),
\end{equation*}
so we have no longer dangerous positions. Now using that 11231 and 3231 are forbidden, we see that $\lambda_0(2232ww^*ww12111112)$ is minimized with the above candidate.

Now let $b\in\{1,2,3\}^\Z$ be such that $b=11133232ww^*w12111232$. By \Cref{lem:selfreplication} we must have that $b=\overline{w}w^*w12111232$. Now we only have to minimize $\lambda_0(b)$, which is achieved with the above candidate because 213 is forbidden.
\end{proof}

In particular for all $m(b)\in(j_0,j_1)$ we have that $b$ is connecting in the past with $\overline{w}$.

\begin{corollary}\label{cor:self-replication}
Let $b\in\{1,2,3\}^\Z$ be a sequence such that $m(b)<j_0+10^{-19}$ and $\lambda_0(b)\in(j_0,j_1)$. Then one has
\begin{equation*}
    \dots b_0^*\dots b_{30}=\overline{w}w^*w12111.
\end{equation*}

\end{corollary}
\begin{proof}
Just use \Cref{cor:fulllocal_uniq} and \Cref{lem:selfreplication}.
\end{proof}

\begin{lemma}\label{lem:max_element}\hfill

$m_1:=\max(M\setminus L)\cap(j_0,j_1)=\lambda_0(\overline{w}w^*w1211111\overline{23})=3.94200115991134146921437465\dots\approx j_0+8.26\cdot 10^{-22}$.    
\end{lemma}
\begin{proof}
By the previous corollary we have $b=\overline{w}w^*w12111$. Since we want to maximize $\lambda_0(b)$, we should continue as $b=\overline{w}w^*w121111$. Now using the inequalities of \Cref{lem:aux_forbidwords} this word must be extended to $b=\overline{w}w^*w1211111$. Since we want to maximize and since $33232www12111111$ is forbidden by that same lemma we should continue with $b=\overline{w}w^*w12111112$. Now using that 11231 and 3231 are forbidden, we have that $\lambda_0(b)$ is maximized with the given candidate. 
\end{proof}

\begin{corollary}
The intervals $(j_0,m_0)$ and $(m_1,j_1)$ are maximal gaps of $M$. 
\end{corollary}

\begin{corollary}
The interval $(j_0,j_1)$ is a maximal gap of $L$ and $j_1\in L^\prime$.
\end{corollary}
\begin{proof}
We will use the fact that Lagrange values of periodic words are dense in $L$. Suppose that $\ell(\overline{a})\in(j_0,j_1)$ where $a\in\{1,2,3\}^N$ is a finite word in $\{1,2,3\}$ that attains its Markov value at $\ell(\overline{a})=m(\overline{a})=\lambda_0(\overline{a_1^*\dots a_N})$.  By \Cref{cor:fulllocal_uniq} and \Cref{cor:self-replication} we must have that
\begin{equation*}
    \overline{a}a_1^*\dots =\overline{w}w^*w12111\dots
\end{equation*}
Hence we must have $\overline{a}=\overline{w}$ and so $\ell(\overline{a})=\ell(\overline{w})=j_0$, a contradiction.

Now we must prove that $j_1\in L^\prime$. For this we use the same characterisation of $L$, by proving that it is the limit of a strictly decreasing sequence of Markov values of periodic words:
\begin{equation*}
    j_1=\lim_{n\to\infty}\lambda_0\left(\overline{2(12)^n232ww^*ww1211111(23)^n}\right)=\lim_{n\to\infty}m\left(\overline{2(12)^n232ww^*ww1211111(23)^n}\right).
\end{equation*}
We only need to check that the Markov value of these periodic words is really attained at the position $*$. We only have three potential dangerous positions (each one inside each copy of $w$)
\begin{itemize}
\item $\lambda_0(1232w^*1211)<j_0-10^{-9}$
\item $\lambda_0(232w^*121111)<j_0-10^{-11}$
\item $\lambda_0(12232www^*w121111123)<j_0+8.31\cdot 10^{-22}<j_1$.
\end{itemize}
Although we checked directly that $j_1\in L^\prime$, this also follows from \cite[Theorem 2]{Flahive77} and the symbolic expression for $j_1$.
\end{proof}

Finally we can prove the complete characterization of $M\setminus L$ in this region.
\begin{theorem}\label{thm:characterization}
We have that
\begin{equation*}
    (M\setminus L)\cap(j_0,j_1)=C_1\cup C_2\cup C_3\cup D_1\cup D_2\cup D_3\cup X,
\end{equation*}
where 
\begin{equation*}
    C_1=\{\lambda_0(\overline{w}w^*w12111112\gamma):2111112\gamma\in\{1,2,3\}^{\N}~\text{not contain any word from $F$}\},
\end{equation*}
\begin{equation*}
    C_2=\{\lambda_0(\overline{w}w^*w12111113\gamma):2111113\gamma\in\{1,2,3\}^{\N}~\text{not contain any word from $F$}\},
\end{equation*}
\begin{equation*}
    C_3=\{\lambda_0(\overline{w}w^*w12111232\gamma):2111232\gamma\in\{1,2,3\}^{\N}~\text{not contain any word from $F$}\},
\end{equation*}
are Cantor sets and
\begin{multline*}
    D_1=\{\lambda_0(\overline{w}w^*w1211111\theta1111121\overline{w^T}):\theta~\text{finite word in 1,2,3},\theta_1\neq 1,(\theta^T)_1\neq 1, [0;\theta]\geq[0;\theta^T], \\ \text{and $211111\theta111112$ not contain any word from $F$}\},
\end{multline*}
\begin{multline*}
    D_2=\{\lambda_0(\overline{w}w^*w1211111\theta23211121\overline{w^T}):\theta~\text{finite word in 1,2,3},\theta_1\neq 1 \\ \text{and $211111\theta232$ not contain any word from $F$}\},
\end{multline*}
\begin{multline*}
    D_3=\{\lambda_0(\overline{w}w^*w12111232\theta23211121\overline{w^T}):\theta~\text{finite word in 1,2,3,}~ [0;\theta]\leq[0;\theta^T], \\ \text{and $2111232\theta2321112$ not contain any word from $F$}\},
\end{multline*}
\begin{multline*}
X=\bigg\{m\left(\overline{w}121111121111121\overline{w^T}\right), m\left(\overline{w}121111121\overline{w^T}\right), m\left(\overline{w}121111123211121\overline{w^T}\right), \\ 
m\left(\overline{w}121112321111121\overline{w^T}\right), m\left(\overline{w}121112323211121\overline{w^T}\right), m\left(\overline{w}1211123211121\overline{w^T}\right)\bigg\},
\end{multline*}
are sets of isolated points in $M$.
\end{theorem}

\begin{proof}
Let $b\in\{1,2,3\}^\Z$ be such that $m(b)=\lambda_0(b)\in(j_0,j_1)$. By \Cref{cor:fulllocal_uniq} and \Cref{cor:self-replication} we must have that $b=\overline{w}w^*w12111\widetilde{\gamma}$ with $\widetilde{\gamma}\in\{1,2,3\}^{\N_{>0}}$. Since $3232w121113$ is forbidden we have $\widetilde{\gamma}_1\neq 3$. In case that $\widetilde{\gamma}_1=2$, then since $3232w1211121, 3232w1211122$ are forbidden we must continue as $b=\overline{w}w^*w1211123$. Since 11231 is forbidden we have that either $b=\overline{w}w^*w12111232$ or $b=\overline{w}w^*w12111233$. In the latter case using \Cref{lem:selfreplication} this implies $b=\overline{w}w^*ww12111$, but since $j_0=\lambda_0(\overline{w}w^*ww12111)<j_0+10^{-35}<m_0$ we must have $\lambda_0(b)=j_0$. In case that $\widetilde{\gamma}_1=1$, since $33232www1211112$, $33232www1211113$, $33232www12111111$ are forbidden we are forced to $b=\overline{w}w^*w1211111\widetilde{\gamma}_3\dots$ with $\widetilde{\gamma}_3\neq 1$. 

Note that in any case $\widetilde{\gamma}$ can not contain $w_r$, because by self-replication (\Cref{cor:self-replication}) we would get that $\overline{w}$ contains 11111, a contradiction. Thus $\gamma:=\widetilde{\gamma}_3\widetilde{\gamma}_4\dots$ does not contain $w_r$. Observe that all the words in $F$ do not contain 1211111, 1111121 nor $1\tau$ so if there are forbidden words they can only appear in $211111\gamma$ or in $2111232\gamma$. In particular if $\gamma$ does not contain $w_r^T$ then $\lambda_0(b)\in C_{1}\cup C_{2}\cup C_{3}$.

If $\gamma$ contains $w_r^T$, pick $N$ minimal such that $\gamma_N\gamma_{N+1}\dots\gamma_{N+31}=w_r^T$, then since $m(b)=m(b^T)<j_1$, by applying \Cref{lem:selfreplication} to $b^T$ we get that either $\gamma_{N-14}\gamma_{N-13}\dots=11121\overline{w^T}$ and so $b=\overline{w}w^*w12111\tau_L\widetilde{\theta}11121\overline{w^T}$ where $\tau_L\in\{112,113,232\}$ with $\widetilde{\theta}=\gamma_1\dots\gamma_{N-15}$ or that $\lambda_0(b)\in X$. 

In the former case, using the transpose of the forbidden words above we know that $\widetilde{\theta}$ must end with 1 or 2. If it ends with 2 then it must end with 32, and since 13211 is forbidden and by the minimality of $N$ we must have that $\widetilde{\theta}=\theta232$. In case that $\widetilde{\theta}$ ends with 1, then by using the transpose of the above forbidden words we must have $\widetilde{\theta}=\theta11$ with the last digit of $\theta$ different from 1. Since 1211111, 1111121 and $23211121w^T23233111=\tau^T1$ are not subwords of any word in $F$, we know that forbidden words can only appear in $211111\theta 111112$, $211111\theta232$ or $2111232\theta2321112$ respectively. So $\lambda_0(b)\in D_1\cup D_2\cup D_3$.

The above argument shows that $(M\setminus L)\cap(j_0,j_1)\subset C_1\cup C_2\cup C_3\cup D_1\cup D_2\cup D_3\cup X$. Conversely given $\lambda_0(b)\in C_1\cup C_2\cup C_3\cup D_1\cup D_2\cup D_3$ (by direct computation we can check that $X\subset (j_0,j_1)$), first note that we have $j_0=\lambda_0(\overline{w}w^*\overline{w})<\lambda_0(\overline{w}w^*w1211111\dots)=\lambda_0(b)\leq m_1<j_1$ (because of the forbidden words 11231 and 3231). Now we need to show that $\lambda_0(b)=m(b)\in M$. It suffices to guarantee that the supremum is attained at the given position. Suppose that for some $N\in\N_{>0}$ we have that $\lambda_0(b)<\lambda_N(b)$. Since $j_0<\lambda_0(b)$ and because of \Cref{cor:fulllocal_uniq}, we will have that $b_{N-15}\dots b_N\dots b_{N+16}$ is equal to $w_r$ or $w_r^T$. But by the definition of the sets $C_i,D_i,X$ and since $w_r\in F$, we will have that is equal to $w_r^T$, that $\lambda_0(b)\not\in C_i$ and moreover the position $N$ must occur after $\theta$, more specifically inside the period $\overline{w^T}$. Note that we have the inequalities
\begin{itemize}
\item $\lambda_0(www^*ww)<j_0+10^{-32}<m_0$,
\item $\lambda_0(ww^*121111)<j_0$,
\item $\lambda_0(ww^*12111232)<j_0$.
\end{itemize}
Hence we must have 
$b_{N-30}\dots b_N^*b_{N+1}\dots=11121w^T(w^*)^T\overline{w^T}$ and either $b_{N-33}b_{N-32}b_{N-31}\in \{232,211,311\}$ (since we assume $\lambda_0(b)\not\in X$). In the case $\lambda_0(b)\in D_2$ we use the inequality $\lambda_0(\overline{w}w^*w121112)<\lambda_0(\overline{w}w^*w121111)$ and in the case that $\lambda_0(b)\in D_1$ the condition $[0;\theta]\geq[0;\theta^T]$ guarantees that $\lambda_N(b)\leq\lambda_0(b)$. 

Finally we prove that $D_1\cup D_2\cup D_3$ is an isolated set in $M$. Given $\lambda_0(b)\in D_1\cup D_2\cup D_3$, suppose we have a sequence $\lambda_0(a^{(n)})=m(a^{(n)})$ that converges to $\lambda_0(b)$. By \Cref{cor:self-replication} we have $a^{(n)}=\overline{w}w^*w12111\dots$. Hence the non-positive part of those sequences all coincide with the non-positive part of $b$, i.e. $\dots a_{-1}^{(n)}a_0^{(n)}=\dots b_{-1}b_0$. In particular using a basic inequality in continued fractions (see \cite[Lemma A.1.]{geometricproperties}) we must have for $n$ sufficiently large 
\begin{itemize}
\item $a^{(n)}=\overline{w}w^*w1211111\theta1111121w^Tw^Tw^T\dots$
\item $a^{(n)}=\overline{w}w^*w1211111\theta23211121w^Tw^Tw^T\dots$
\end{itemize}
In all cases the word $w_r^T$ is appearing, so by using \Cref{lem:selfreplication} we will obtain that $a^{(n)}$ is connecting to $\overline{w^T}$, that is $a^{(n)}=b$ for all large $n$. The same proof applies to all elements of $X$.
\end{proof}

\begin{remark}
By using \Cref{prop:mixing_criteria}, one can show that the set 
\begin{multline*}
    K=\{[0;1,1,1,3,3,2,3,2,w,1,2,1,1,1,1,1,\gamma]:11111\gamma\in\{1,2,3\}^{\N} \\
    ~\text{not contain any word from $F$ and}~\gamma_1\neq 1\}
\end{multline*}
is a topologically mixing dynamically defined Cantor set.
\end{remark}

\subsection{A lower bound on $\dim_H(M\setminus L)$}\label{subsec:lower_bound_dimM-L}

In this subsection we show that $D(3.94) = 0.5945611\dots$, so in particular by \eqref{eq:precise_dimension} we have 
\begin{equation*}
    \dim_H((M\setminus L)\cap(3.94,3.943)) > D(3.94) =  0.594561\dots
\end{equation*}
since $(M\setminus L)\cap(j_0,j_1)$ belongs to the good interval $[3.9362, 3.943767)$ (see \Cref{thm:good_intervals_list}). 

In fact, we claim that there is a gap $(\nu,\mu)$ of $M$ with extremes $\nu=m(\overline{21}113^*3311\overline{12})=3.939301\dots$ and $\mu=m(\overline{111331113333^*})=3.941015\dots$. Moreover, the subshift $\Sigma_{\nu}$ is of finite type given by
\begin{multline*}
    \Sigma_{\nu}=\Sigma(A)=\{\alpha\in\{1,2,3\}^\Z: 131, 312, 313, 1323, 1322, 21132, 13211, 11132, 3111333 \\ ~\text{and their transposes are not substrings of}~\alpha\}.
\end{multline*}

As a consequence, using the algorithms developed in \cite{MMPV} (the code is available at \url{https://github.com/Polevita/Gauss_Cantor_sets}), one can confirm that
\begin{equation*}
    D(3.94)=D(\nu)=\dim_H(K_\nu)=\dim_H(K(A))=0.5945611\dots.
\end{equation*}

Suppose $b\in\{1,2,3\}^\Z$ is such that  $m(b)=\lambda_0(b)\in(\nu,\mu)$, so in particular $b=3^*$. Moreover since $\lambda_0(33^*3)<3.7$, $\lambda_0(23^*3)<3.8$, $\lambda_0(23^*2)<3.9$, and 131 is forbidden we can assume that either $b=13^*2$ or $b=13^*3$. If $b=13^*2$, then using the forbidden words one is forced to continue as $\lambda_0(3113^*212)<3.9377$. Hence assume $b=113^*3$. Since $\lambda_0(113^*31)<3.923$ we continue as either $b=113^*32$ or $b=113^*33$. Since $\lambda_0(113^*32)\leq\lambda_0(\overline{31}113^*32\overline{12})<3.939$ (because 213 is forbidden), we can assume that $b=113^*33$. Now since $\lambda_0(3113^*33)<3.88$, $\lambda_0(2113^*33)<3.9$, $\lambda_0(11113^*33)<3.928$, we are forced to either $b=21113^*33$ or $b=31113^*33$. Since 313 and 312 are forbidden, we see that $\lambda_0(21113^*33)\leq\nu$. On the other hand, the minimum Markov value of a bi-infinite sequence containing 3111333 corresponds to $\mu$. Indeed, any continuation of that subword different from 33111333311133 gives values at least $3.941016>\mu$. Finally, since this continuation has the form $\beta^T3\theta3\beta$ with $\beta=1113311$ and $\theta=33$, using \Cref{lem:minimizing} and the fact that if 131, 11132, 31113331, 31113332, 311133333, 311133332, 11331113333113, 11331113333112, 113311133331111, 113311133331112 are forbidden, then a continued fraction in $\{1,2,3\}$ beginning with $[0;1,1,1,3,3]$ is minimized with $[0;\overline{1,1,1,3,3,1,1,1,3,3,3,3}]$, we conclude that $\lambda_0(31113^*33)\geq\mu$, a contradiction.

\begin{remark}
The point $\mu=m(\overline{111331113333^*})=3.941015605\dots$ is isolated in $M$. The right border of the other gap is given by $m(\overline{12}33111331113^*3331113311133\overline{21})=3.941015626\dots$
\end{remark}

\subsection{Characterisation of the region associated with $\boldsymbol{w = 21133311121}$}

As discussed in the introduction, the second and third authors in joint work with Matheus~\cite{JMM} recently studied a region of $M\setminus L$ associated with the odd non-semi-symmetric word $w=21133311121$. We take this opportunity to state the full description of this region.

Let $j_0=m(\overline{w})=3.938776241981\dots\in L$ and $j_1\in L^\prime$ given by 
\begin{equation*}
    j_1=m(\overline{21}2331113311321231133311121ww^*w22\overline{32})\approx 9.52145\cdot 10^{-12}.
\end{equation*}

Let $F$ be the following set of words and their transposes:
\begin{itemize}
\item 131, 312, 313, 2231, 3231, 23111, 23112, 11231, 3331113, 2111333111212,
\item 1113331112121, 1113331112122, 11111333111212, 21111333111212, 11133311121232,
\item 11133311121233,  133311121212, 21231133311121212, 
\item  $11w21$, $2w21$, $3w21$, $21w212$, $221w21$, $321w21$, $w2111$, $w2112$, $1121w21$.
\end{itemize}

\begin{theorem}
We have that
\begin{equation*}
    (M\setminus L)\cap(j_0,j_1)=C_1\cup C_2\cup C_3 \cup \bigcup_{i,j}D_{i,j}\cup X
\end{equation*}
where 
\begin{equation*}
    C_1=\{\lambda_0(\gamma 2121ww^*\overline{w}):\gamma21212 \in\{1,2,3\}^{\N}~\text{not contain any word from $F$}\},
\end{equation*}
\begin{equation*}
    C_2=\{\lambda_0(\gamma 21231133311121w^*\overline{w}):\gamma2123 \in\{1,2,3\}^{\N}~\text{not contain any word from $F$}\},
\end{equation*}
\begin{equation*}
    C_3=\{\lambda_0(\gamma 331133311121w^*\overline{w}):\gamma331133 \in\{1,2,3\}^{\N}~\text{not contain any word from $F$}\},
\end{equation*}
are Cantor sets, where
\begin{multline*}
    D_{i,j}=\{\lambda_0(\overline{w^T}w^T\tau_L\theta\tau_Rw^*\overline{w}):\theta~\text{finite word in 1,2,3},~ [0;\tau_R^T,\theta^T]\geq[0;\tau_L,\theta], \\ \text{and $\tau_L\theta\tau_R$ not contain any word from $F$}\},
\end{multline*}
are discrete sets of isolated points with $(\tau_1,\tau_2,\tau_3)=(2121w,21231133311121,331133311121)$ and $(\tau_R,\tau_L^T)=(\tau_i,\tau_j)$ and finally where
\begin{multline*}
X=\bigg\{m\left(\overline{w^T}12111333113331133311121\overline{w}\right), m\left(\overline{w^T}1211133311331133311121\overline{w}\right), \\ 
m\left(\overline{w^T}121113331133311121\overline{w}\right), m\left(\overline{w^T}1211133311321231133311121\overline{w}\right), m\left(\overline{w^T}121113331132121\overline{w}\right),  \\ m\left(\overline{w^T}121231133311121\overline{w}\right),
m\left(\overline{w^T}1212121\overline{w}\right), m\left(\overline{w^T}12121\overline{w}\right), m\left(\overline{w^T}121\overline{w}\right)\bigg\}
\end{multline*}
is a finite set of isolated points.
\end{theorem}

\begin{remark}
By using \Cref{prop:mixing_criteria}, one can show that each of the sets
\begin{equation*}
    K_i=\{[0;3,3,1,1,2,\tau_i^T,\gamma]:\gamma^T\tau_i\in\{1,2,3\}^{\N}
    ~\text{does not contain any word from $F$}\},
\end{equation*}
is a topologically mixing dynamically defined Cantor set.
\end{remark}

\section{New maximal gaps in $M$ near 3.943}\label{sec:gaps} 

We now demonstrate the existence of two new maximal gaps in the Markov spectrum $M$. These are in the vicinity of the last visible gaps in the computer approximations produced by Delecroix, Matheus and the third author~\cite[Figure 5]{DMM}.

\subsection{First gap}

There is gap in the Markov spectrum with left endpoint given by

\begin{equation*}
    m(\overline{23}331113^*3\overline{32})=3.94254\dots
\end{equation*}
and right endpoint given by

\begin{equation*}
    m(\overline{23}3311133113212311333^*111331113331132123113311133\overline{32})=3.943304\dots
\end{equation*}

\begin{lemma}\hfill

\begin{enumerate}
\item If $b=13^*1$, $3^*13$, $3^*12$, $323^*1$, $223^*1$, $23^*111$, $13^*2112$, $13^*2113$, $33113^*2111$, $23113^*21111$, $123113^*21112$ then $\lambda_0(b)>3.9438$
\item If $b=331113^*331$ then  $\lambda_0(b)>3.94317$.
\end{enumerate}
\end{lemma}

Suppose $m(\underline{b})\leq 3.94317$. Since 131 is forbidden, we have that either $b_{-1}b_0^*b_1$ is equal to $13^*2$ or $13^*3$. Suppose first it is equal to $13^*3$. Since 313, 213 are forbidden we are forced to $113^*3$. Since $\lambda_0(3113^*3)<3.9$, $\lambda_0(2113^*3)<3.9$, $\lambda_0(113^*31)<3.93$, $\lambda_0(113^*32)<3.94$, $\lambda_0(11113^*33)<3.93$, $\lambda_0(21113^*33)<3.941$ we should continue as $31113^*33$. Since 131, 23111, 331113331 are forbidden we just have the options  $\lambda_0(331113^*333)<3.942$ and $\lambda_0(331113^*332)<3.9427$. Actually one has
\begin{equation*}
    \lambda_0(331113^*332)\leq m(\overline{23}331113^*3\overline{32})=3.94254\dots
\end{equation*}
where we used that 3231 and 133311133 are forbidden.

Now suppose $b_{-1}b_0^*b_1=13^*2$. Since 313, 213, 1322, 1323 are forbidden we are forced to $113^*21$. Since 11132 is forbidden, we have two cases:
\begin{itemize}
\item If it continues as $2113^*21$, then since 312 is forbidden we have 
\begin{equation*}
    \lambda_0(2113^*21)\geq m(\overline{21}13^*\overline{21})=3.94337\dots
\end{equation*}

\item If it continues as $3113^*21$ then since $\lambda_0(3113^*212)<3.938$ and 213 is forbidden we have to continue as $3113^*211$. Now using that 132113, 132112, 331132111, 131, 2231, 3231, 2311321111, 2311321112, 131, 11132 are forbidden, we are forced to  

\begin{equation*}
    \lambda_0(123113^*211133)>3.943227.
\end{equation*}

\end{itemize}

The minimum Markov value of a bi-infinite sequence that contains 123113211133 is the same as the palindrome 3311123113211133 which is
\begin{equation*}
    m(\overline{21}1133311123^*113211133311\overline{12})=3.94342\dots
\end{equation*}

where we used that if the words 313, 312, 1113331113 are forbidden, then a continued fraction in $\{1,2,3\}$ that begins with $[0;2,1,1,1,3,3]$ is minimized when $[0;2,1,1,1,3,3,3,1,1,\overline{1,2}]$. 

\begin{remark}
The middle word is a palindrome and the period is semi-symmetric.
\end{remark}

The minimum Markov value of a bi-infinite sequence that contains 331113331 is the same as the palindrome 113331113311133311 which is
\begin{equation*}
    m(\overline{23}3311133113212311333111331113^*331132123113311133\overline{32})=3.943304\dots
\end{equation*}

where we used that if the words 131, 1322, 1323, 33113211, 23111, 32123112, 311331113331 are forbidden, then a continued fraction in $\{1,2,3\}$ that begins with $[0;3,3,1,1]$ is minimized when $[0;3,3,1,1,3,2,1,2,3,1,1,3,3,1,1,1,3,3,\overline{3,2}]$.

\begin{remark}
The middle word is a palindrome and the period is semi-symmetric.
\end{remark}

\subsection{Second gap}

The above argument shows that we have a local uniqueness on $(3.9433, 3.94337)$. Indeed, we have that  $b_{-5}\dots b_0^*\dots b_4=331113^*3311$.  

\begin{lemma}
If $b=2331113^*3311$, $3331113^*3311$, $331113^*33111$, $331113^*33112$, $11331113^*331133$, $311331113^*331132$, $211331113^*331132$, or $1111331113^*3311321$ then $\lambda_0(b)>3.94331$.
\end{lemma}

Recall that we have that $b_{-5}\dots b_0^*\dots b_4=331113^*3311$. Using the above forbidden words and 1323, 1322, 131, 313, 312, 33113211 we get that it must extend to $111331113^*33113212$. Since $\lambda_0(3111331113^*33113212)<3.9433055$ and  11113311133311321 is forbidden we should continue as $2111331113^*33113212$ which satisfies $\lambda_0(2111331113^*33113212)>3.9433068$. In fact we have
\begin{equation*}
    \lambda_0(3111331113^*33113212)\leq m(\overline{12}31133311133111331113^*33113\overline{21})=3.94330534\dots
\end{equation*}

This is the left endpoint of the gap.

The right endpoint of the gap is given by the minimum Markov value of a bi-infinite sequence that contains 211133111333113212, which is
\begin{equation*}
    m(\overline{21}11331113^*331132123113311133\overline{32})=3.94330716\dots
\end{equation*}
where we used that if 313, 312, 23111, 32123112, 311331113331 are forbidden then a continued fraction in $\{1,2,3\}$ that begins with $[0;3,2,1,2]$ is minimized with $[0;3,2,1,2,3,1,1,3,3,1,1,1,3,3,\overline{3,2}]$. 

\appendix

\section{Markov value algorithm}\label{app:alg}

We want to compute the Markov value of a doubly periodic sequence of the form $\omega=\overline{p_1}\tau\overline{p_2}$. 

\begin{enumerate}
\item The first step is to write $\omega$ in the simplest terms.
\begin{enumerate}
\item  First, one writes each $p_i$ with its minimal period.
\item If $\tau$ and $p_1$ start with the same digit, then we apply a cyclic shift to $p_1$ and erase the first digit of $\tau$. Do this until it is no longer possible.
\item If $\tau$ is non-empty, repeat the same process with $p_2$.
\item If $\tau$ is non-empty, then we are done. Otherwise, we check if $p_1=p_2$. If $p_1\neq p_2$ and they begin with the same first digit, then we do a cyclic shift in both until we arrive at $p_1,p_2$ with different first term.
\end{enumerate}
\item We claim that the Markov value of $\omega=\overline{p_1}\tau\overline{p_2}$ can be found on the maximum between the positions $p_1p_1\tau p_2p_2$ and the periodic Markov values $\overline{p_1},\overline{p_2}$. Indeed, observe that for any position inside a $p_1$, say $p_1^*$ the left will always be $\overline{p_1}$, while for the right we have two choices: $p_1$ or $\tau\overline{p_2}$. By the simplification above at least one of them \textbf{must win}. If $p_1$ is even, then $p_1$ always wins and the value at this position is maximized with $\overline{p_1}$ to the right, that is 
\begin{equation*}
    \lambda_i(\overline{p_1}p_1^*\dotsb)\leq\lambda_i(\overline{p_1}p_1^*\overline{p_1})\leq m(\overline{p_1}).
\end{equation*}
In case $p_1$ is odd then we maximize with $p_1$, but in the next question $p_1$ \textbf{must lose}, so the position is maximized with 
\begin{equation*}
    \lambda_i(\overline{p_1}p_1^*\dotsb)\leq\lambda_i(\overline{p_1}p_1^*p_1\tau\overline{p_2}).
\end{equation*}
The other side is analogous.
\end{enumerate}

\begin{remark}
The step (i).(a) is not always necessary, unless one is in the special situation described in step (i).(d), namely, when the word has the form $\omega=\overline{p_1}~\overline{p_2}$.
\end{remark}

\section{Berstein's intervals}\label{app:ber}\hfill

Here, for completeness, we list the intervals that Berstein was investigating in~\cite{Berstein}. They can be found in Table 3, pp. 117-120, of the book containing that article.

They are
\begin{itemize}
\item (3.001494,  3.001889)
\item (3.011905,  3.012844)
\item (3.026437,  3.032218)
\item (3.043396,  3.049177)
\item (3.049177,  3.049296)
\item (3.072944,  3.091490)
\item (3.116089,  3.121575)
\item (3.122174,  3.122890)
\item (3.123669,  3.123758)
\item (3.128940,  3.129192)
\item (3.166577,  3.265590)
\item (3.280963,  3.284892)
\item (3.284892,  3.285169)
\item (3.289043,  3.292022)
\item (3.325167,  3.326240)
\item (3.332493,  3.333233)
\item (3.333971,  3.334734)
\item (3.341564,  3.351944)
\item (3.405792,  3.415638)
\item (3.432835,  3.437413)
\item (3.440922,  3.449666)
\item (3.459882,  3.460606)
\item (3.462426,  3.463064).
\end{itemize}

We will not check whether Berstein's intervals are indeed good intervals. If we plot these intervals in green on our Figure~\ref{fig:goodints} (we restrict now to the region between 3 and 4), we obtain the following:

\begin{figure}[h]
\centering
\begin{tikzpicture}[scale=14]

% Interval 1
\fill[opacity = 0.2, blue] (3.05082,-.03) -- (3.122183, -.03) -- (3.122183, .03) -- (3.05082,.03) -- cycle;

% Interval 2
\fill[opacity = 0.2, blue] (3.1299,-.03) -- (3.285441, -.03) -- (3.285441, .03) -- (3.1299,.03) -- cycle;

% Interval 3
\fill[opacity = 0.2, blue] (3.28603,-.03) -- (3.28729, -.03) -- (3.28729, .03) -- (3.28603,.03) -- cycle;

% Interval 3.5
\fill[opacity = 0.2, blue] (3.29296,-.03) -- (3.29335, -.03) -- (3.29335, .03) -- (3.29296,.03) -- cycle;

% Interval 3.75
\fill[opacity = 0.2, blue] (3.33396,-.03) -- (3.33475, -.03) -- (3.33475, .03) -- (3.33396,.03) -- cycle;

% Interval 4
\fill[opacity = 0.2, blue] (3.359,-.03) -- (3.423, -.03) -- (3.423, .03) -- (3.359,.03) -- cycle;

% Interval 5
\fill[opacity = 0.2, blue] (3.464,-.03) -- (3.84, -.03) -- (3.84, .03) -- (3.464,.03) -- cycle;

% Interval 6
\fill[opacity = 0.2, blue] (3.873,-.03) -- (3.930691, -.03) -- (3.930691, .03) -- (3.873,.03) -- cycle;

% Interval 7
\fill[opacity = 0.2, blue] (3.93616,-.03) -- (3.943767, -.03) -- (3.943767, .03) -- (3.93616,.03) -- cycle;

% Interval 8
\fill[opacity = 0.2, blue] (3.944054,-.03) -- (3.971606, -.03) -- (3.971606, .03) -- (3.944054,.03) -- cycle;

% Interval 9
\fill[opacity = 0.2, blue] (3.97995,-.03) -- (3.9857, -.03) -- (3.9857, .03) -- (3.97995,.03) -- cycle;

%Berstein intervals
\foreach \x/\y in {3.001494/3.001889,
 3.011905/3.012844,
 3.026437/3.032218,
 3.043396/3.049177,
 3.049177/3.049296,
 3.072944/3.091490,
 3.116089/3.121575,
 3.122174/3.122890,
 3.123669/3.123758,
 3.128940/3.129192,
 3.166577/3.265590,
 3.280963/3.284892,
 3.284892/3.285169,
 3.289043/3.292022,
 3.325167/3.326240,
 3.332493/3.333233,
 3.333971/3.334734,
 3.341564/3.351944,
 3.405792/3.415638,
 3.432835/3.437413,
 3.440922/3.449666,
 3.459882/3.460606,
 3.462426/3.463064}
	\fill[opacity = 0.4, green] (\x,-.0225) -- (\y, -.0225) -- (\y, .0225) -- (\x,.0225) -- cycle;

% Line
\draw[->, thick] (2.95,0) -- (4.05,0);

% Cuts
\draw[thick] (3,0.02) -- (3,-0.02) node[below] {3};
\draw[thick] (3.334384,0.02) -- (3.334384,-0.02) node[below] {$t_1$};
\draw[thick] (4,0.02) -- (4,-0.02) node[below] {4};

\end{tikzpicture}
\caption{The good intervals produced in this paper before 4 are depicted in blue. Berstein's intervals are depicted in green.}
\end{figure}

\printbibliography

\end{document}